\documentclass[a4paper]{article}

\usepackage[utf8]{inputenc} 
\usepackage[T1]{fontenc} 
\usepackage{amsmath,amsthm,latexsym} 
\usepackage[fixamsmath]{mathtools}
\mathtoolsset{showmanualtags,mathic,centercolon} 
\usepackage{amssymb,amsfonts} 
\usepackage{booktabs} 
\usepackage[UKenglish]{babel} 
\usepackage[svgnames]{xcolor} 
\usepackage{algpseudocode} 
\usepackage{algorithm} 
\usepackage{url}

\usepackage{bibentry}
\newcommand{\ignore}[1]{}
\newcommand{\nobibentry}[1]{[{\let\nocite\ignore\bibentry{#1}}]}

\usepackage[margin=1.25in]{geometry} 
\usepackage[labelfont=sl,textfont=sl]{subcaption} 
\usepackage[font={small,sl},labelsep=period]{caption} 

\numberwithin{equation}{section} 
\theoremstyle{plain}
\newtheorem{theorem}{Theorem}[section] 

\theoremstyle{definition}
\newtheorem{definition}[theorem]{Definition}

\theoremstyle{remark}

\newcommand{\appref}[1]{Appendix~\ref{#1}}
\newcommand{\secref}[1]{Section~\ref{#1}}

\renewcommand{\algref}[1]{Algorithm~\ref{#1}}

\newcommand{\figref}[1]{Figure~\ref{#1}}

\newcommand{\tabref}[1]{Table~\ref{#1}}

\newcommand{\R}{\mathbb{R}} 
\newcommand{\bigO}{\mathcal{O}} 
\DeclareMathOperator*{\argmin}{arg\,min} 
\newcommand{\defeq}{:=} 

\def\be{\begin{equation}}
\def\ee{\end{equation}}

\renewcommand{\b}[1]{\mathbf{#1}} 
\renewcommand{\t}[1]{\widetilde{#1}} 

\newcommand{\bx}{\b{x}}
\newcommand{\by}{\b{y}}
\newcommand{\bs}{\b{s}}

\newcommand{\bg}{\b{g}}

\algrenewcommand\algorithmicrequire{\textbf{Input:}}
\algrenewcommand\algorithmicensure{\textbf{Output:}}

\begin{document}
\title{Escaping local minima with derivative-free methods:\\ a numerical investigation}
\author{
	Coralia Cartis\thanks{Mathematical Institute, University of Oxford, Radcliffe Observatory Quarter, Woodstock Road, Oxford, OX2 6GG, United Kingdom (\texttt{cartis, oliver.sheridan-methven@maths.ox.ac.uk}). This work was supported
          by the EPSRC Centre for Doctoral Training in Industrially
          Focused Mathematical Modelling (EP/L015803/1) in
          collaboration with the Numerical Algorithms Group Ltd.}
\and 
Lindon Roberts\thanks{Mathematical Sciences Institute, Building 145, Science Road, Australian National University, Canberra ACT 2601, Australia (\texttt{lindon.roberts@anu.edu.au}).}
\and
Oliver Sheridan-Methven\footnotemark[1]}

\date{\today}
\maketitle

\begin{abstract}
We apply a state-of-the-art, local derivative-free solver, Py-BOBYQA \cite{TR-DFOLS}, to global optimization problems, and propose an algorithmic improvement that is beneficial in this context. 
Our numerical findings are illustrated on a commonly-used but small-scale test set of global optimization problems \cite{Ali2005} and associated noisy variants, and on hyperparameter tuning for the machine learning test set MNIST \cite{Snoek2012}. 
As Py-BOBYQA is a model-based trust-region method, we compare mostly (but not exclusively) with other global optimization methods for which (global) models are important, such as Bayesian optimization and response surface methods; we also consider state-of-the-art representative deterministic and stochastic codes, such as DIRECT and CMA-ES. 
As a heuristic for escaping local minima, we find numerically that Py-BOBYQA is competitive with global optimization solvers for all accuracy/budget regimes, in both smooth and noisy settings.
In particular, Py-BOBYQA variants are best performing for smooth and multiplicative noise problems in high-accuracy regimes. 
As a by-product, some preliminary conclusions can be drawn on the relative performance of the global solvers we have tested with default settings.

\end{abstract}

\\


\section{Introduction}\label{Intro}

Derivative-Free Optimization (DFO) refers to a large class of methods that are designed to find local optima without needing derivative information  to  be provided  or to even exist, a ubiquitous practical situation arising in problems where the often black-box function evaluations may be expensive or inaccurate/noisy \cite{CSV}.   Global Optimization (GO) algorithms, whether derivative-based or derivative-free, tackle the (NP-hard) question of finding  global optima, a much more computationally intense effort but a crucial requirement for some applications \cite{Locatelli2013}. 
Both classes of methods aim to go beyond the remit of the successful and scalable  local derivative-based optimization solvers that, by construction, get easily trapped in/close to local optima; this larger applicability  comes at the cost of reduced scalability, especially in the case of generic GO algorithms.
 
Despite only local theoretical guarantees of convergence, 
DFO methods have long been regarded as a more suitable alternative for GO than local derivative-based methods.
Indeed, Rios and Sahinidis \cite{Rios2013} recently performed extensive benchmarking of various GO and DFO solvers, and found that in some instances, such as when close to a solution, state-of-the-art DFO methods such as NEWUOA \cite{Powell2007} can be comparable, or even more efficient than GO solvers (to reach the same accuracy), but that generally, GO is better than DFO. 
We have recently improved upon the performance of NEWUOA, with our solver Py-BOBYQA \cite{TR-DFOLS}, which is a Python-conversion and extension of  BOBYQA\footnote{Note that NEWUOA and BOBYQA are very similar codes, with the latter allowing bound constraints.}\cite{Powell2009}.
Our improvements to BOBYQA that are relevant here take the form of {\it restarts}; we have found in \cite{TR-DFOLS} that 
they significantly improve the performance of BOBYQA in the presence of noise, and also sometimes give improved accuracy in the noiseless case on standard nonlinear least-squares problems. 
Encouraged by this, in this paper, we further modify the restart procedure in a way that is (even more) advantageous to global optimization, and test both the usual Py-BOBYQA (with restarts) and its GO modification against state of the art GO codes---on standard GO test problems and on a parameter tuning problem from machine learning---with encouraging results, as we describe below.

\paragraph{Literature review}
There are several types of DFO algorithms: pattern or direct search \cite{Kolda2003}, the well-known simplex/Nelder-Mead algorithm \cite{MHWright}, linesearch methods \cite{Powell1998}, implicit filtering \cite{Kelley2011}, conjugate directions \cite{Powell1998},  model-based trust region methods \cite{CSV} and many more, including extensions to constrained problems;
Chapter 9 in \cite{Nocedal2006} gives a good survey and the reference textbook \cite{CSV} provides the details. It is well-known that model-based trust region methods are competitive within existing DFO methods and solvers \cite{Custodio2017, More2009}, 
and we have recently compared different methods with similar conclusions \cite{TR-Miniproject-Lindon, TR-DFOLS}. 
Thus, in this work, we focus on model-based trust-region DFO methods for unconstrained and bound-constrained optimization problems, in particular on  Py-BOBYQA variants \cite{TR-DFOLS}, which improve upon BOBYQA \cite{Powell2009}.  
These methods build, at each iteration, a {\it local quadratic model} of the objective function by interpolation of function values, that is minimized over a trust region ball, to find an improvement on the current, best-available, estimate of the solution. When no improvement can be achieved,
the size of the trust region and the geometry of the interpolation set are adjusted in the run of the algorithm. The evaluation cost of these methods is larger than the number of variables for the initial iteration (when the first local quadratic model is set up), but it is only order $1$ subsequently, as usually only one, or a small number of, points leaves and enters the interpolation set for the remaining iterations. Py-BOBYQA has additional features such as restarts, that we describe below.

GO methods typically fall in two broad categories: deterministic such as branch-and-bound, Lipschitz optimization or 
response surface/surrogate methods, and stochastic such as random search, evolutionary methods, simulated annealing and Bayesian methods; there are also mixed variants,
with both random and deterministic features (such as multi-start methods, as well as variants of surrogate methods). For our comparisons here, we focus mostly on GO methods that also use simplified models of the objective in their construction. Hence we consider state-of-the-art Bayesian optimization solvers, that use (global) Gaussian or other models, and response surface methods that use radial basis function surrogate models;
and also SNOBFIT, that is similar to Py-BOBYQA in using local quadratic models, combined with a branching strategy to enforce global searches. 
 For reference, we also select the well-known CMA-ES \cite{Hansen1996, HansenOstermeier2001} and  DIRECT \cite{Jones1993} methods, as representatives for 
the evolutionary and deterministic classes of methods, respectively\footnote{Due to our use of Python for the DFO code, for the sake of like-for-like comparisons, we considered GO solvers written in or interfacing to Python.}.

In terms of test problems, we have selected a class of 50 standard GO test problems \cite{Ali2005}, of dimensions up to $50$ variables,  and we also generate noisy variants of these problems by random perturbations, since the stochastic regime is very much within the design scope of both DFO, Bayesian and response surface methods. This test set is representative and common to usual GO solvers such as DIRECT, but probably less familiar to Bayesian solvers which are very much developed and used in the machine learning community.  
To play to the strengths of Bayesian solvers, we also choose a case study in machine learning, namely, an algorithm parameter tuning problem for training classifiers on the MNIST image data set. 

The benchmarking paper \cite{Rios2013}  compares a large number of DFO and GO solvers on a large test set, as mentioned above, but do not include  Bayesian solvers as extensively as we consider here (but do include the related response surface methods). The results in \cite{Ghanbari} compare model-based DFO with Bayesian solvers on machine learning test sets, so as to emphasize that there is no need to use GO solvers for parameter tuning of training algorithms since local DFO gives as good low-accuracy solutions as Bayesian solvers, with reduced computational and evaluation cost. Our aim here is to report preliminary testing of our recent and significantly improved model-based DFO solver (Py-BOBYQA) in similar (but reduced size) contexts as these two papers, so as to ascertain its suitability or otherwise for GO problems.




\paragraph{Our contributions and conclusions}
In some situations, local DFO methods may be sufficient for GO, as they may happen to sample points in a region of attraction of the global minimum, or may have been initialized at a starting point in this region.
Compared to GO solvers, local DFO methods benefit from being well-suited for exploiting extra problem information, such as a good starting point or problem structure, and are not as affected by the curse of dimensionality.
However, local DFO methods focus heavily on `exploitation' (i.e.~searching regions where the objective is known to be small) at the expense of `exploration' (i.e.~searching regions where the objective has not yet been evaluated), and so may fail to solve GO problems.
Here, we propose a modification of Py-BOBYQA to improve its `exploration' ability.
This technique is a heuristic, and we do not claim any guarantees that it can find global minima, but we show encouraging results that it improves the ability of local DFO methods to escape local minima (and possibly find the global minimum), while maintaining the benefits of local solvers.

Specifically, Py-BOBYQA has a restart mechanism that whenever the algorithm no longer makes much progress, it resets the (small) size of trust region radius to its initial (large) value. 
It also regenerates completely or partially the tightly-clustered interpolation set to be more spread out in the enlarged trust region. 
The restarts procedure was designed with noisy problems in mind, where we often observed that the trust region gets too small too quickly, before the solver reaches desired accuracy. 
In \cite{TR-DFOLS}, we found that for local optimization, restarts---and the associated auto-detection mechanism---help improve the accuracy of solutions for noisy problems, and do not adversely affect, and may even help, solver performance on smooth problems; our testing for this last point was limited to nonlinear least squares problems \cite{More2009}.
The modification of Py-BOBYQA proposed here is an additional improvement, {\it adaptive restarts}, with global optimization in mind: every time a restart is called, the trust region size is not just reset to its original value used at the start of the algorithm, but possibly increased, to allow a larger local region of search and optimization.
These larger search regions improve the `exploration' ability of Py-BOBYQA, without compromising its `exploitation' ability.

We test both the default restart configurations and the adaptive restart modification of Py-BOBYQA on the GO test problems and the machine learning test set. 
We find that for low accuracy solutions on the GO test sets, most solvers have similar performance on a given budget, for both smooth and noisy problems. 
Differences start appearing for longer budgets/high accuracy settings; here, adaptive Py-BOBYQA variants are the most efficient solvers for smooth  problems on a given evaluation budget, while for noisy problems, an adaptive Py-BOBYQA variant (with full quadratic models) is amongst the more efficient solvers, competitive with DIRECT and the surrogate solver PySOT \cite{pysot2015}. 
For the MNIST test set, due to the high cost of the evaluations, only low accuracy solutions were attempted. 
In this setting, PySOT and DIRECT are the most efficient solvers, while adaptive Py-BOBYQA (with full quadratic models) performs acceptably. 
Overall, we find that including the adaptive restart strategy in Py-BOBYQA is advantageous for escaping local minima, and that there is reasonably small variability in performance due to the choice of starting point.
In terms of the other solvers, DIRECT and PySOT are most competitive across all problems types.
 
\paragraph{Structure of the paper}
Section \ref{sec_algo_framework} describes the Py-BOBYQA algorithm and its improvements, Section \ref{sec_global} summarizes the main solvers we compare with, while Section \ref{sec_test_methodology} presents the testing methodology, including test problems and measuring performance of solvers. We show our numerical findings in Section \ref{sec_numerics}, comparing the different Py-BOBYQA variants, and then the best performing one with the remaining solvers; we give our conclusions in Section \ref{sec_conclusions}.



\section{Py-BOBYQA: a local DFO solver with improvements}
\label{sec_algo_framework}

The Py-BOBYQA algorithm is a model-based derivative-free optimization method, where, as is typical for such methods, a local model for the objective is constructed by interpolation and minimized on each iteration \cite{CSV}.
More specifically, Py-BOBYQA \cite{TR-DFOLS} is based on Powell's BOBYQA \cite{Powell2009}, which uses underdetermined quadratic interpolation for model construction.
At each iteration $k$, we have a collection of points $Y_k\subset\R^n$, where we allow $|Y_k| \in [n+1, (n+1)(n+2)/2]$, and we construct a local model for the objective
\be f(\bx_k+\bs) \approx m_k(\bs) = c_k + \bg_k^{\top}\bs + \frac{1}{2}\bs^{\top}H_k \bs, \label{eq_quadratic_models} \ee
satisfying the interpolation conditions
\be m_k(\by_t-\bx_k) = f(\by_t), \quad \text{for all $\by_t\in Y_k$.} \label{eq_bobyqa_interp_conditions} \ee
If $|Y_k|<(n+1)(n+2)/2$, then there are multiple solutions to \eqref{eq_bobyqa_interp_conditions}, so we use the remaining degrees of freedom by solving 
\be \min_{c_k,\bg_k,H_k} \|H_k-H_{k-1}\|_F^2 \quad \text{subject to \eqref{eq_bobyqa_interp_conditions}}, \label{eq_bobyqa_interp_problem} \ee
where $\|\cdot\|_F$ is the matrix Frobenius norm and we have the convention $H_{-1}=0$.
The value of $|Y_k|$ is a user-specified input; larger values of $|Y_k|$ capture more objective information, but require more initial objective evaluations and have a higher linear algebra cost for solving \eqref{eq_bobyqa_interp_problem}.
For smooth problems, we use the default value $|Y_k|=2n+1$, which provides a balance between capturing curvature information and low initialization cost.
However, for noisy problems, we use the default value $|Y_k|=(n+1)(n+2)/2$, which we have found yields a more robust solver \cite{TR-DFOLS}.

Once we have a method for constructing local models \eqref{eq_quadratic_models}, we calculate new steps using a trust-region method \cite{CSV}, which ensures global convergence (i.e.~convergence to a local minimum regardless of starting point).
That is, we maintain a parameter $\Delta_k>0$ and compute a new step by solving
\be \bs_k \approx \argmin_{\|\bs\|\leq\Delta_k} m_k(\bs). \label{eq_tr_subproblem} \ee
Efficient algorithms exist for solving \eqref{eq_tr_subproblem} approximately (e.g.~\cite{Conn2000}).
We then evaluate the objective $\bx_k+\bs_k$, and determine if the step produce a sufficient decrease in the objective, by computing:
\be r_k = \frac{\text{actual decrease}}{\text{expected decrease}} \defeq \frac{f(\bx_k) - f(\bx_k+\bs_k)}{m_k(\b{0}) - m_k(\bs_k)}. \label{eq_tr_ratio} \ee
If $r_k$ is sufficiently large, we accept the step ($\bx_{k+1}=\bx_k+\bs_k$) and increase $\Delta_k$, otherwise we reject the step ($\bx_{k+1}=\bx_k$) and decrease $\Delta_k$.

The last key component of the algorithm is the maintenance of the interpolation set $Y_k$.
We want $Y_k$ to consist of points close to $\bx_k$, with `good' geometry.
The quality of $Y_k$ is calculated based on Lagrange polynomials; that is, the polynomial functions $\{\ell_1(\by), \ldots, \ell_p(\by)\}$ for $p=|Y_k|$ defined by
\be \ell_t(\by_s) = \begin{cases}1, & t=s, \\ 0, & t\neq s, \end{cases} \qquad \forall s,t = 1,\ldots,p. \ee
Given these Lagrange polynomials, the measure of the quality of $Y_k$ is given by the following definition.

\begin{definition}[$\Lambda$-poised; Chapter 6, \cite{CSV}] \label{def_poised}
	For $B\subset\R^n$ and $\Lambda>0$, the set $Y_k$ is $\Lambda$-poised in $B$ if
	\be \max_{t=1,\ldots,p}\: \max_{\by\in B} |\ell_t(\by)| \leq \Lambda, \label{eq_poised} \ee
	for all $\by\in B$, where $\{\ell_1(\by), \ldots, \ell_p(\by)\}$ are the Lagrange polynomials for $Y_k$.
\end{definition}

Updating $Y_k$ proceeds as follows: at each iteration we add $\bx_k+\bs_k$ to $Y_k$ (and remove a point far from $\bx_{k+1}$), and at certain times we must also replace a point in $Y_k$ with a new point designed solely to improve the $\Lambda$-poisedness of $Y_k$ (rather than decreasing the objective).

\paragraph{The complete Py-BOBYQA algorithm}
The basic Py-BOBYQA algorithm is summarized in \algref{alg_pybobyqa}; a complete description can be found in \cite{Powell2009,TR-DFOLS}.
Its implementation also allows optional bound constraints $\b{a} \leq \bx \leq \b{b}$; the main impact of this is the need for a specialized routine for solving \eqref{eq_tr_subproblem} with bound constraints, which we take from \cite{Powell2009}.
When bound constraints are supplied, Py-BOBYQA can internally shift and scale the variables to lie in $[0,1]$.

\begin{algorithm}
	\small{
	\begin{algorithmic}[1]
		\Require Starting point $\bx_0\in\R^n$, initial trust region radius $\Delta_0>0$ and integer $n+1\leq p\leq (n+1)(n+2)/2$. 

		\State Build an initial interpolation set $Y_0\subset B(\bx_0,\Delta_0)$ of size $p$, with $\bx_0\in Y_0$.
		\For{$k=0,1,2,\ldots$} \label{ln_loop}
			\State Given $\bx_k$ and $Y_k$, construct the model $m_k(\bs)$ \eqref{eq_quadratic_models} satisfying \eqref{eq_bobyqa_interp_conditions}.
			\State Approximately solve the trust region subproblem \eqref{eq_tr_subproblem} to get a step $\bs_k$. \label{ln_trs}
			\If{$\|\bs_k\| < \Delta_k/2$}
				\State \underline{Safety Phase}: Set $\bx_{k+1}=\bx_k$, $\Delta_{k+1}<\Delta_k$, and form $Y_{k+1}$ by improving the geometry of $Y_k$.
				\State \textbf{goto} line \ref{ln_loop} (i.e.~next iteration).
			\EndIf
			\State Evaluate $f(\bx_k+\bs_k)$ and calculate ratio $r_k$ \eqref{eq_tr_ratio}.
			\If{$r_k \geq 0.1$}
				\State \underline{Successful Phase}: Accept step $\bx_{k+1}=\bx_k+\bs_k$ and set $\Delta_{k+1}>\Delta_k$.
			\ElsIf{\textit{geometry of $Y_k$ is not good}}
				\State \underline{Model Improvement Phase}: Reject step $\bx_{k+1}=\bx_k$, set $\Delta_{k+1}<\Delta_k$ and improve the geometry of $Y_{k+1}$.
			\Else
				\State \underline{Unsuccessful Phase}: Reject step $\bx_{k+1}=\bx_k$ and set $\Delta_{k+1}<\Delta_k$. \label{ln_rho_redn}
			\EndIf \label{ln_loop_end}
			\State Form $Y_{k+1}=Y_k\cup\{\bx_k+\bs_k\}\setminus\{\by\}$ for some $\by\in Y_k$.
		\EndFor
	\end{algorithmic}
	} 
	\caption{Py-BOBYQA for general minimization; see \cite{Powell2009,TR-DFOLS} for full details.}
	\label{alg_pybobyqa}
\end{algorithm}

\paragraph{Comparison with original BOBYQA}
The underlying Py-BOBYQA algorithm is the same as the original (Fortran) BOBYQA \cite{Powell2009}.
However, there are several ways in which Py-BOBYQA is an improvement over BOBYQA:
\begin{itemize}
	\item Py-BOBYQA allows the user to specify an interpolation set size of $|Y_k|=n+1$, compared to $|Y_k|\geq n+2$ as required by BOBYQA (which gives $H_k=0$ at each iteration, from \eqref{eq_bobyqa_interp_problem} and $H_{-1}=0$). However, we note that Py-BOBYQA does not use this value by default (see above);
	\item Py-BOBYQA has a larger choice of termination conditions than BOBYQA. Specifically, Py-BOBYQA has all the termination conditions of BOBYQA, but can also terminate when it has found a sufficiently small objective value, sees slow objective decrease over several iterations, or all interpolation points have objective value within a user-specified noise level; 
	\item Py-BOBYQA allows different default trust-region update parameters depending on whether or not the objective is noisy, and allows the use of flexible sample averaging strategies to improve performance for noisy problems; and
	\item Py-BOBYQA implements a multiple restarts mechanism, described in \secref{sec_restarts}. This mechanism was introduced in \cite{TR-DFOLS} to improve its performance for noisy problems, and is used in place of sample averaging (but used in conjunction with different default trust-region parameters).
\end{itemize}
In addition, simplify the linear algebra used in the model construction stage, to allow for a more streamlined code.
Solving the interpolation problem \eqref{eq_bobyqa_interp_problem} and building Lagrange polynomials both require the solution of a linear system.
In BOBYQA, this is done by maintaining the inverse of the associated matrix, which is updated via the Sherman-Morrison-Woodbury formula when interpolation points are changed.
In Py-BOBYQA, we solve the linear system directly.

Numerical testing in \cite{TR-DFOLS} demonstrates that Py-BOBYQA has similar performance to BOBYQA for local optimization of smooth objectives, and substantially improved performance for noisy objectives.

\subsection{Multiple restarts feature} \label{sec_restarts}
An important feature of Py-BOBYQA, originally designed to improve its robustness for noisy problems \cite{TR-DFOLS} is the ability to restart the solver.
When the objective is noisy, this stops Py-BOBYQA from stagnating in regions far from local minima; however, we will see that it is also useful for escaping from local minima.
A restart is triggered in four scenarios:
\begin{enumerate}
	\item The trust region radius $\Delta_k$ becomes sufficiently small;
	\item The decrease in the objective over consecutive iterations is too slow;
	\item If the objective values $f(\by)$ for $\by\in Y_k$ are all within some (user-specified) noise level; or
	\item Over consecutive iterations, the trust region radius is consistently being decreased, and the changes $\log \|\bg_k-\bg_{k-1}\|$ and $\log \|H_k-H_{k-1}\|_F$ are consistently increasing.
\end{enumerate}
The last two scenarios are designed primarily for noisy objectives.

When a restart is called, we increase the trust region radius $\Delta_{k+1} = \Delta_{\mathrm{reset}}$ and improve the geometry of $Y_k$ in this new trust region via two possible mechanisms:
\begin{description}
	\item[\normalfont\textit{Hard restart:}] Rebuild $Y_k$ in the new trust region from scratch using the same mechanism as how $Y_0$ was originally constructed; and,
	\item[\normalfont\textit{Soft restart (moving $\bx_k$):}] First, move $\bx_k$ to a geometry-improving point in the new trust region, shifting the trust region to this new point. 
	Then, move the $N-1<p$ points in $Y_k$ which were closest to the old value of $\bx_k$ to geometry-improving points in the new trust region centred at the (shifted) $\bx_k$. 
	Py-BOBYQA continues from whichever of these $N$ new points has the least objective value. 
\end{description}
The final solution returned by the solver takes the optimal value seen so far, including the saved endpoints from previous restarts.
The soft restart approach with $N=\min(3,p)$ is the default approach in Py-BOBYQA.
The value of $\Delta_{\mathrm{reset}}$ is $\Delta_0$ by default.

We note that multiple restarts are different to the interpolation set updating in the model-improvement phase (related to condition \eqref{eq_poised}).
Model-improving steps in \algref{alg_pybobyqa} are used to ensure that the interpolation points are sufficiently close to the current iterate $\bx_k$, and that they span all of $\R^n$ well; this guarantees that the interpolation model is a good approximation to the true objective within the current trust region.
Without restarts, we have theoretical guarantees that $\Delta_k\to 0$ as $k\to\infty$, and so our interpolation points will cluster together.
In the case of noisy problems, this can occur far from local minima, and the solver stagnates---in this situation, model-improving phases occur, but simply ensure that the points are well-spaced within the (very small) trust region.
Restarts couple model-improving steps with a very large increase in the trust-region radius, which allows Py-BOBYQA to explore larger regions and avoid this stagnation.
Here, we argue that the same mechanism is also useful for helping Py-BOBYQA to escape from local minima.

\paragraph{Adaptive restarts} The new restart set up in Py-BOBYQA, to help with (faster) escape from local minima proceeds as follows. Instead of setting $\Delta_{\mathrm{reset}}$ to $\Delta_0$ every time a restart is triggered (which we refer to as a fixed restart), 
Py-BOBYQA allows $\Delta_{\mathrm{reset}}$ to be increased by a constant factor  (default is 1.1) if the previous restart failed to produce a reduction in the objective. Thus $\Delta_{\mathrm{reset}}$ is still set to $\Delta_0$ for the first restart, but will be progressively increased on restarts that do not make progress. 

\paragraph{Connection to multistart GO}
Our restarts mechanism has some similarity to multistart methods for GO, in which several runs of a local optimization algorithm are initialized from different starting points \cite{Locatelli2013,Hare2018}.
Multistart methods can be enhanced by techniques for identifying the basins of attraction for different local minima, such as clustering \cite{RinnooyKan1987,RinnooyKan1987a} and topographs \cite{Torn1994}---these two approaches are combined in \cite{Ali1994}.
Multistart approaches have also been used for local DFO (e.g.~\cite{Custodio2011} or \cite[Example 3.5]{Audet2017}.
However, in multiple restarts, we only ever use a single initialization point for a run, and instead use the final iterate of a run as the starting point for the restarted run (this may also be viewed simply as a continuation of the initial run).
As mentioned above, our goal here is to develop a method that enhances local DFO with greater exploration of the feasible region---and so retaining the benefits of local methods, such as scalability---rather than to build a method with guarantees for GO.

\begin{figure}[t]
	\centering
	\begin{subfigure}[b]{0.48\textwidth}
		\includegraphics[width=\textwidth]{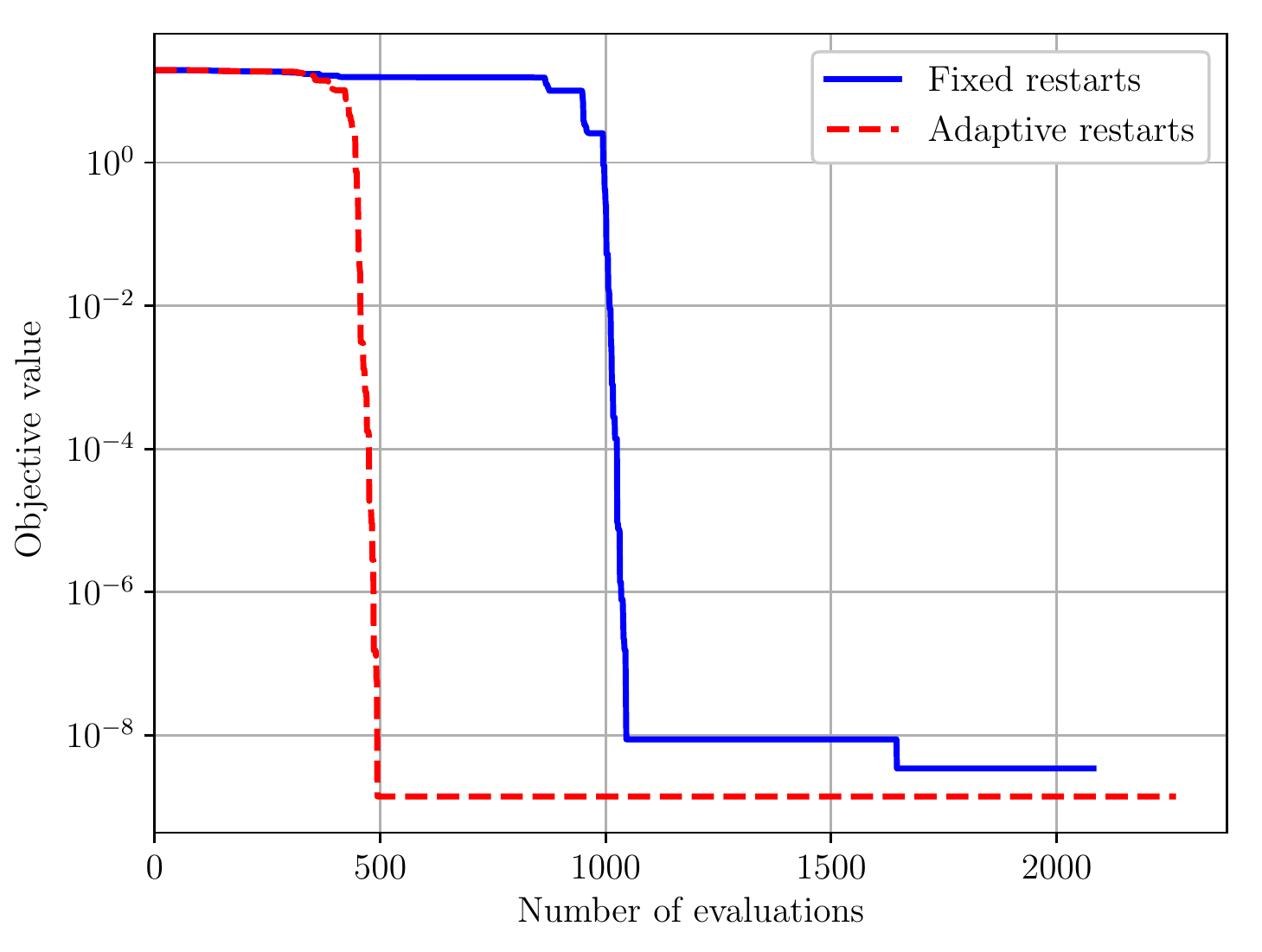}
		\caption{Objective decrease versus evaluation budget.}
	\end{subfigure}
	~
	\begin{subfigure}[b]{0.48\textwidth}
		\includegraphics[width=\textwidth]{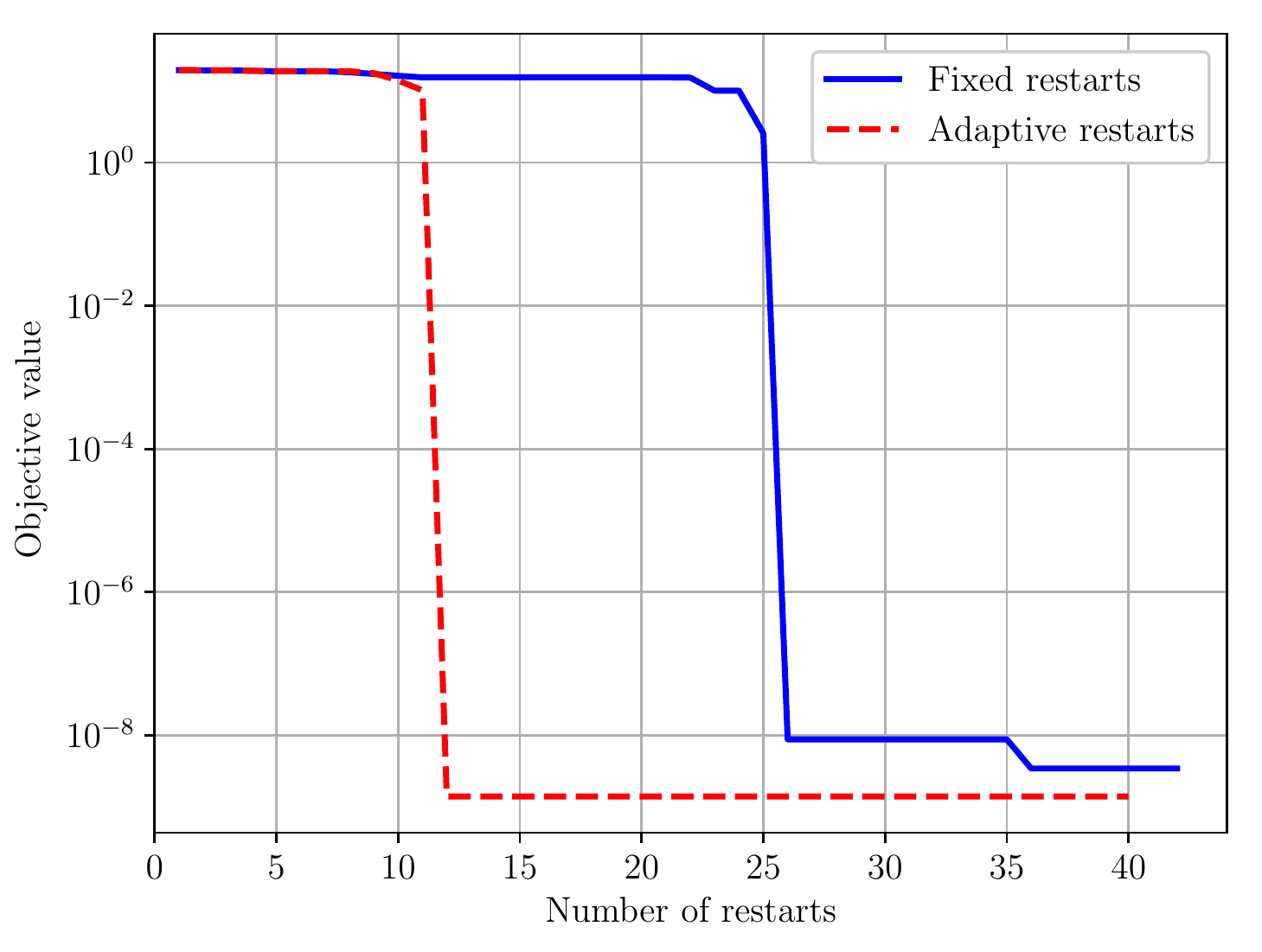}
		\caption{Objective decrease versus number of restarts.}
	\end{subfigure}
	\\
	\begin{subfigure}[b]{0.48\textwidth}
		\includegraphics[width=\textwidth]{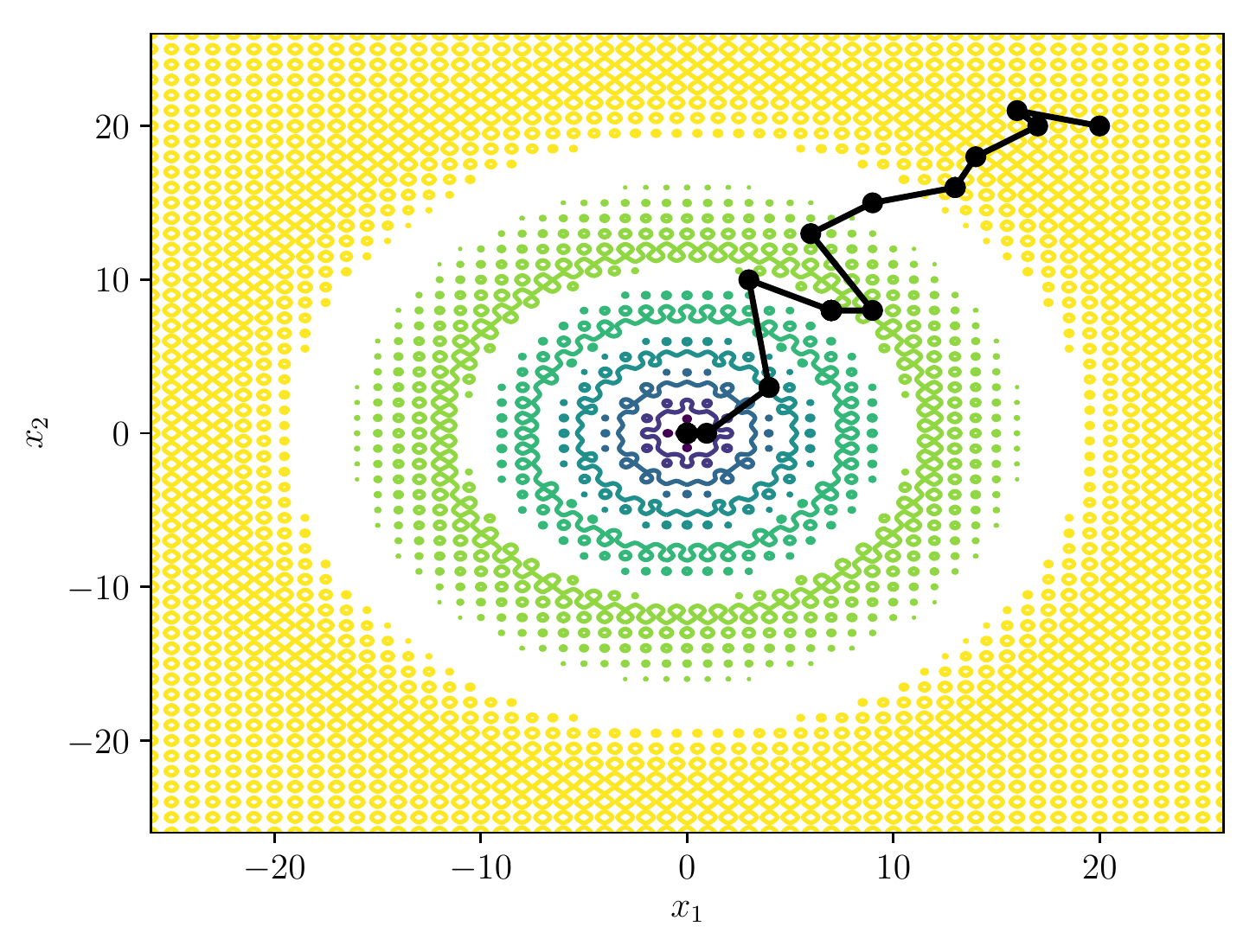}
		\caption{Fixed restarts, path of iterates after each restart.}
	\end{subfigure}
	~
	\begin{subfigure}[b]{0.48\textwidth}
		\includegraphics[width=\textwidth]{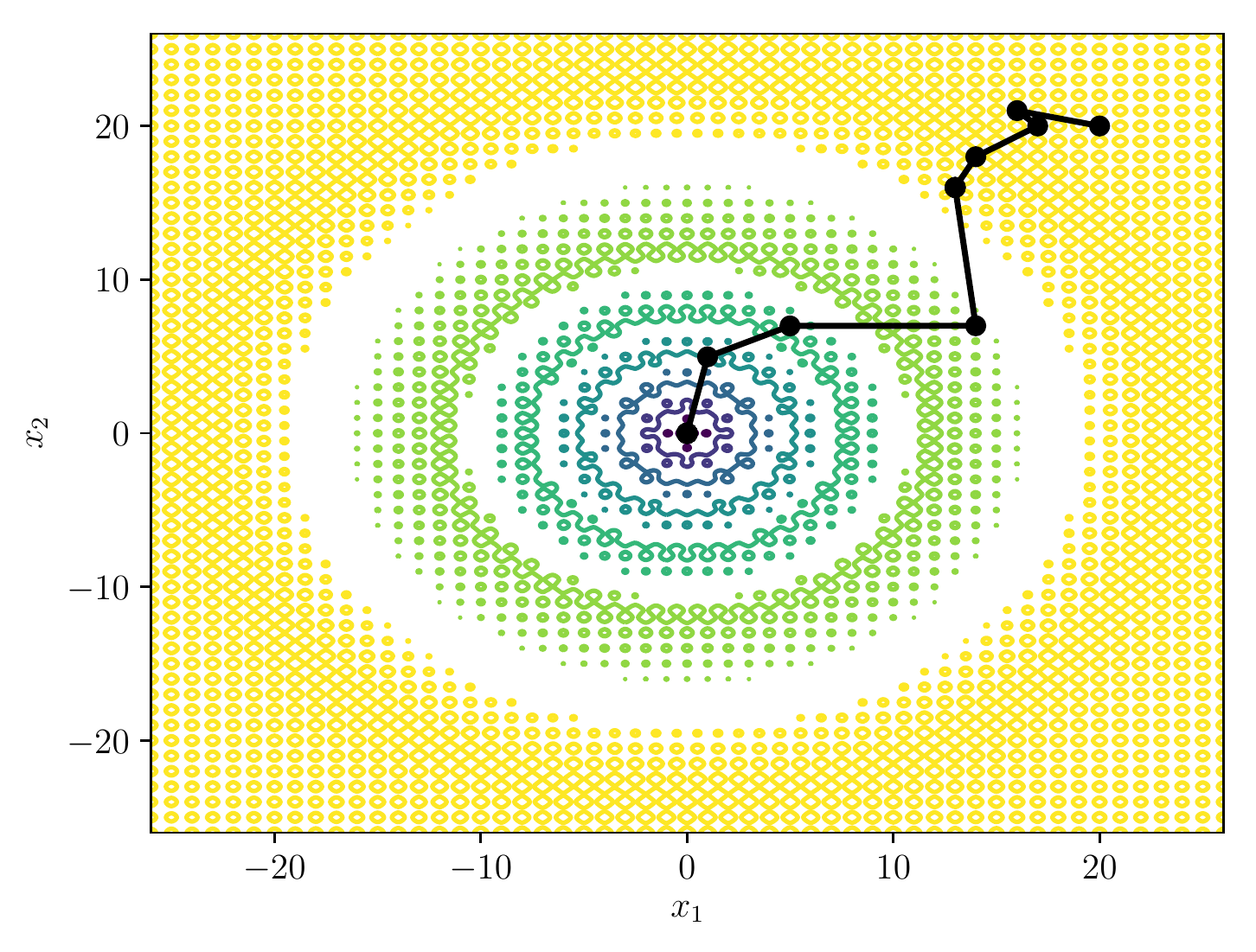}
		\caption{Adaptive restarts, path of iterates after each restart.}
	\end{subfigure}
	\caption{Demonstration of Py-BOBYQA with multiple restarts on the Ackley function.}
	\label{fig_Ackley}
\end{figure}

\subsection{Examples of multiple restarts functionality} \label{sec_restarts_demo}
We now provide several simple examples to demonstrate the benefit of multiple restarts for escaping local minima.

\paragraph{Ackley function} 
We use the well-known Ackley function, that has a global minimum surrounded by a high number of local minima, and illustrate the impact of using fixed and adaptive restarts in Py-BOBYQA on finding the global minimum.
In  \figref{fig_Ackley},  Py-BOBYQA with $2n+1$ initialization points is applied to
 the two-dimensional Ackley function with bound constraints
$\bx\in [-26,26]^2$ and no noise, from starting point $\bx_0=[20,\,20]^{\top}$ and with a maximum of $10^3 (n+1)$ evaluations.
Both usual (fixed) restarts and adaptive restarts are tested.
From this starting point, Py-BOBYQA without restarts converges to a nearby local minimum, as expected for a local solver---it returns the point $\bx=[20-\delta,\,20-\delta]^{\top}$ with $\delta\approx 6.878\times 10^{-4}$---and so needs multiple restarts to escape from this point.

In \figref{fig_Ackley}, we show the results of augmenting Py-BOBYQA with fixed and adaptive restarts.
We first show the objective value achieved by both restart variants, measured by number of evaluations, and by number of restarts (where we plot the best value found so far, just before each time a restart is called).
We then show the location of the iterates from each restart variant, again just before each time a restart is called.

Both fixed and adaptive restarts enable Py-BOBYQA to reach the global minimum to high accuracy---unlike Py-BOBYQA without restarts.
Furthermore, by comparison to fixed restarts, adaptive restarts both achieve a larger reduction in the cost function within fewer evaluations and require fewer restarts to achieve this improvement. 
This is because the adaptivity allows us to more quickly explore larger areas, and so we are more likely to take larger steps (which increases our ability to find better local minima).

\begin{figure}[t]
	\centering
	\begin{subfigure}[b]{0.48\textwidth}
		\includegraphics[width=\textwidth]{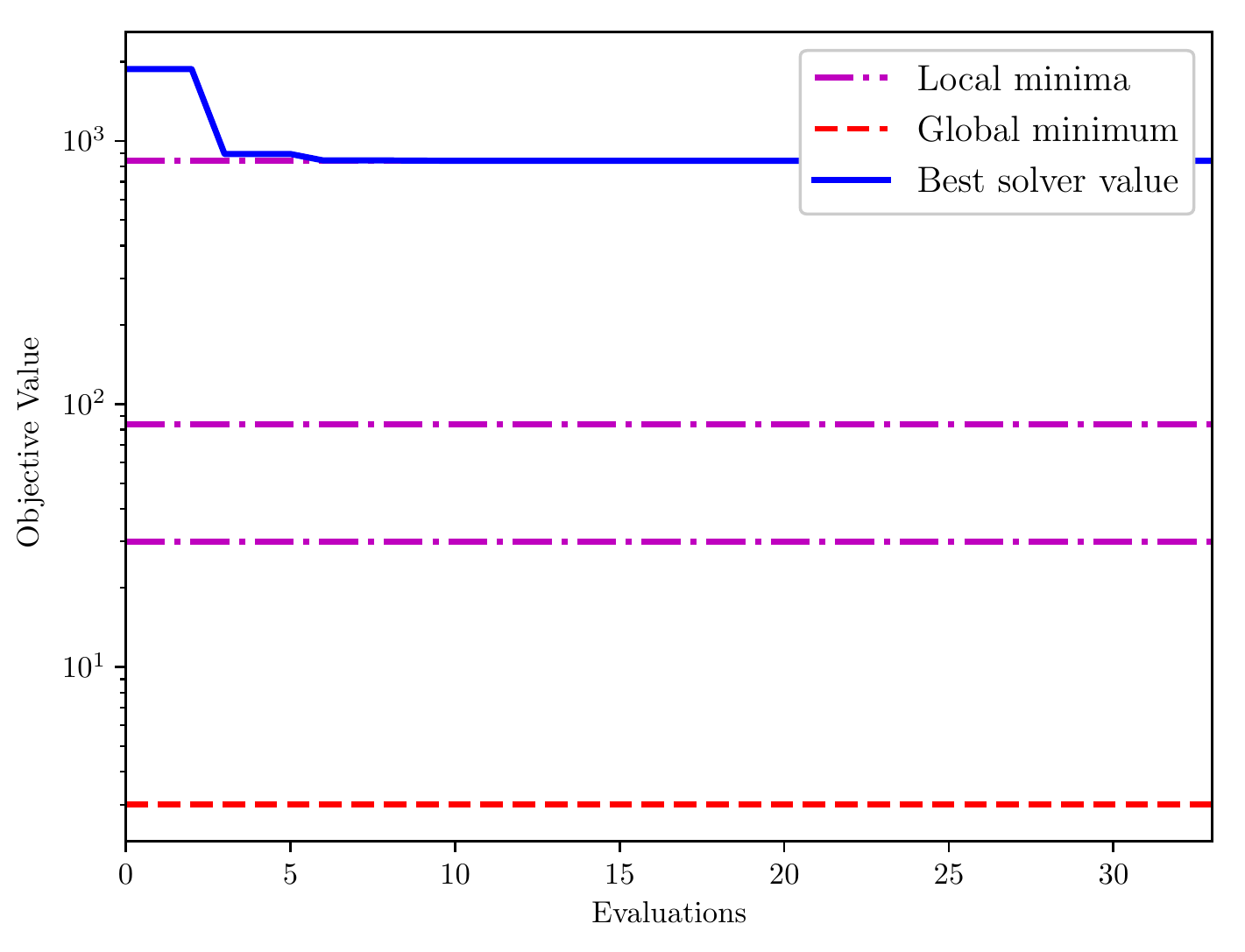}
		\caption{No restarts: objective value.}
	\end{subfigure}
	~
	\begin{subfigure}[b]{0.48\textwidth}
		\includegraphics[width=\textwidth]{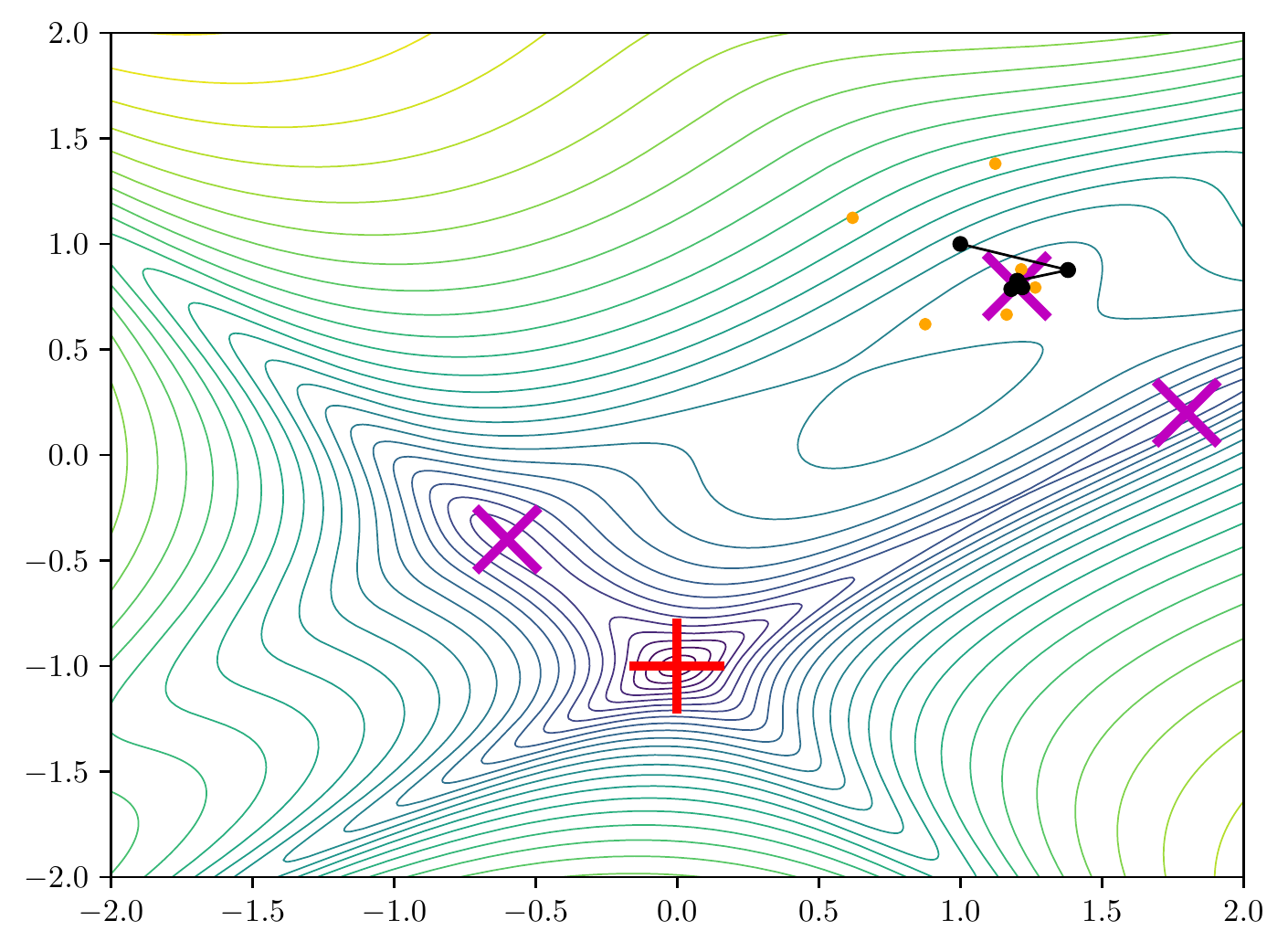}
		\caption{No restarts: evaluated points.}
	\end{subfigure}
	\\
	\begin{subfigure}[b]{0.48\textwidth}
		\includegraphics[width=\textwidth]{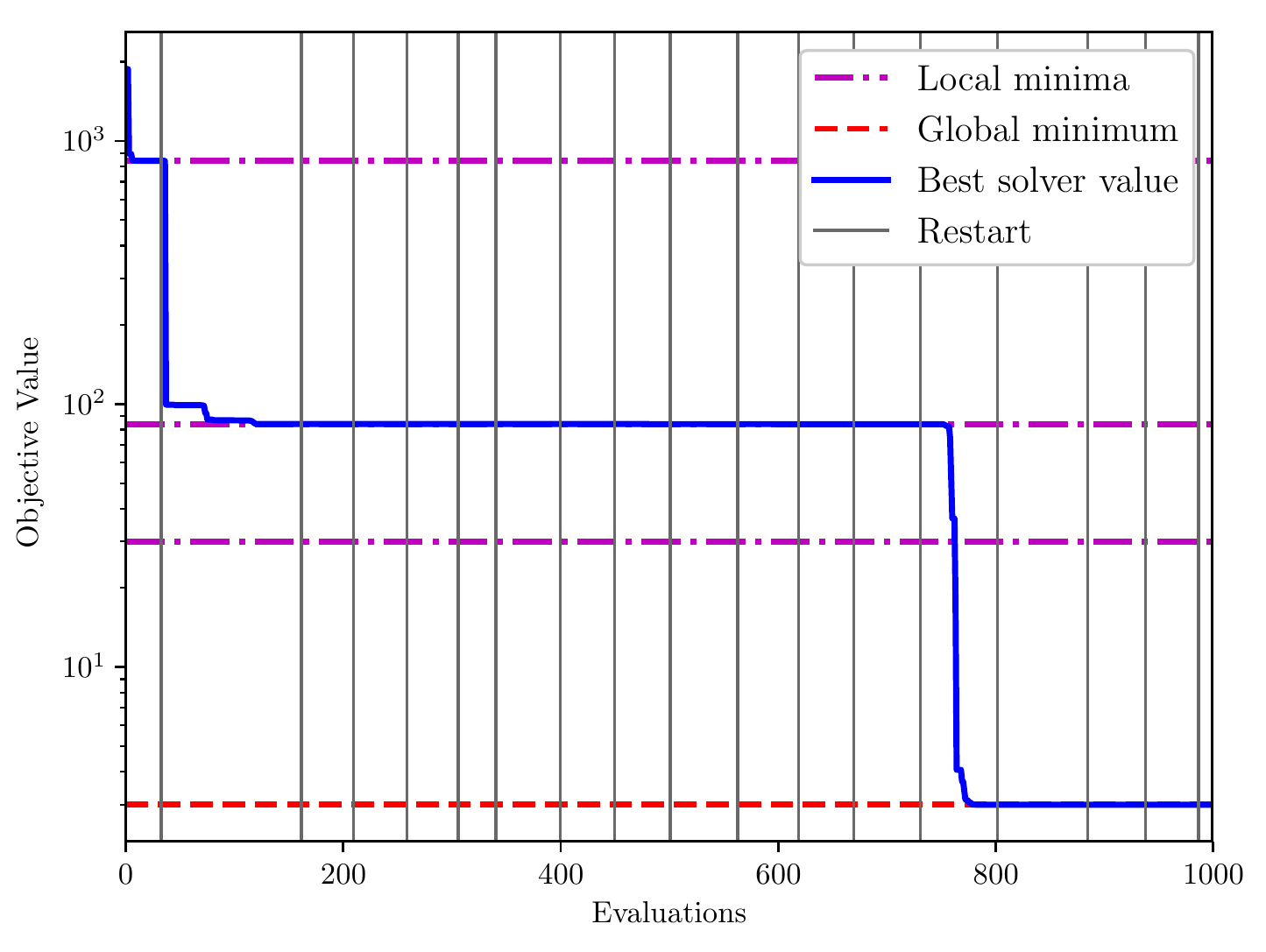}
		\caption{Fixed restarts: objective value.}
	\end{subfigure}
	~
	\begin{subfigure}[b]{0.48\textwidth}
		\includegraphics[width=\textwidth]{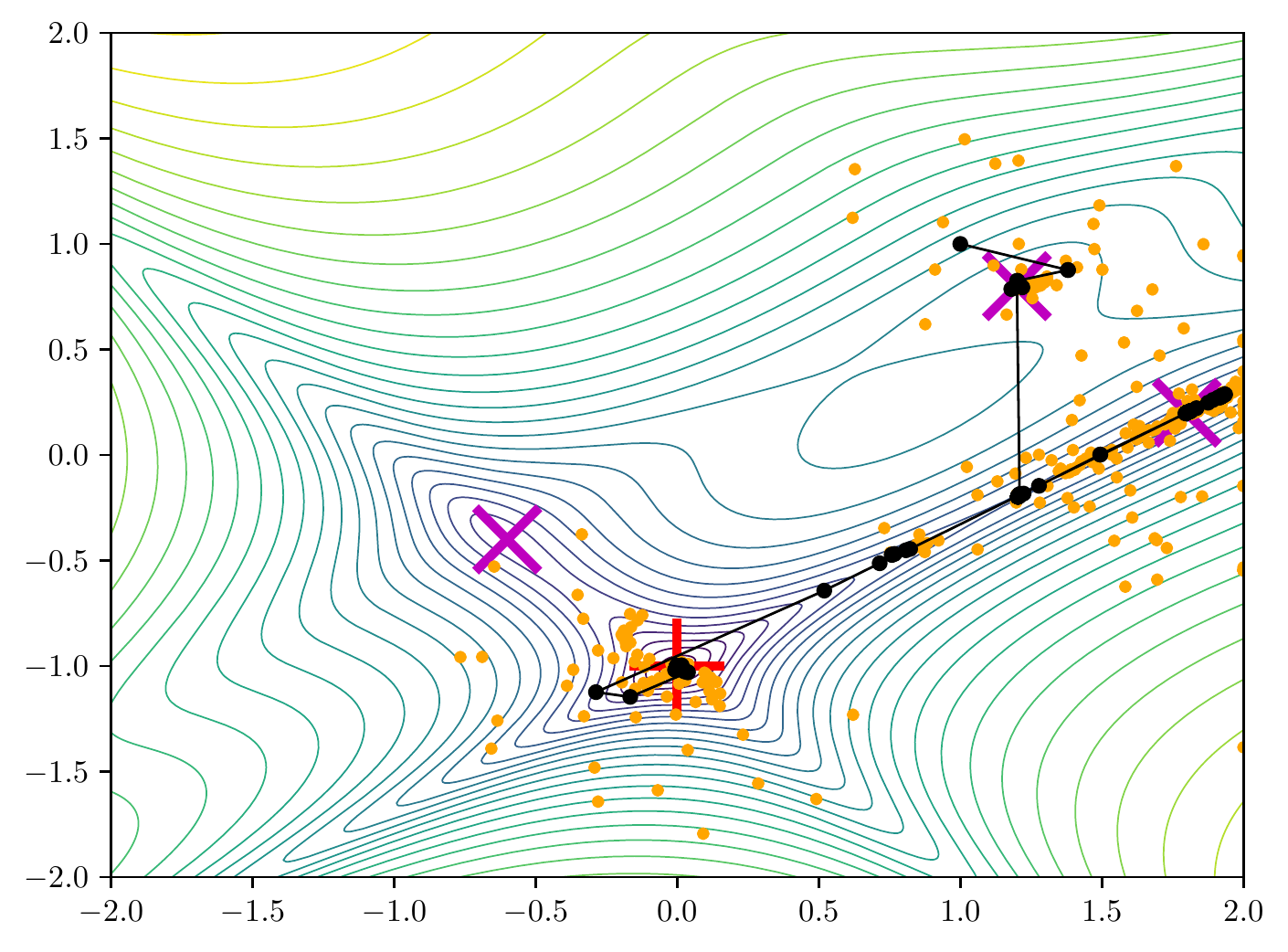}
		\caption{Fixed restarts: evaluated points.}
	\end{subfigure}
	\\
	\begin{subfigure}[b]{0.48\textwidth}
		\includegraphics[width=\textwidth]{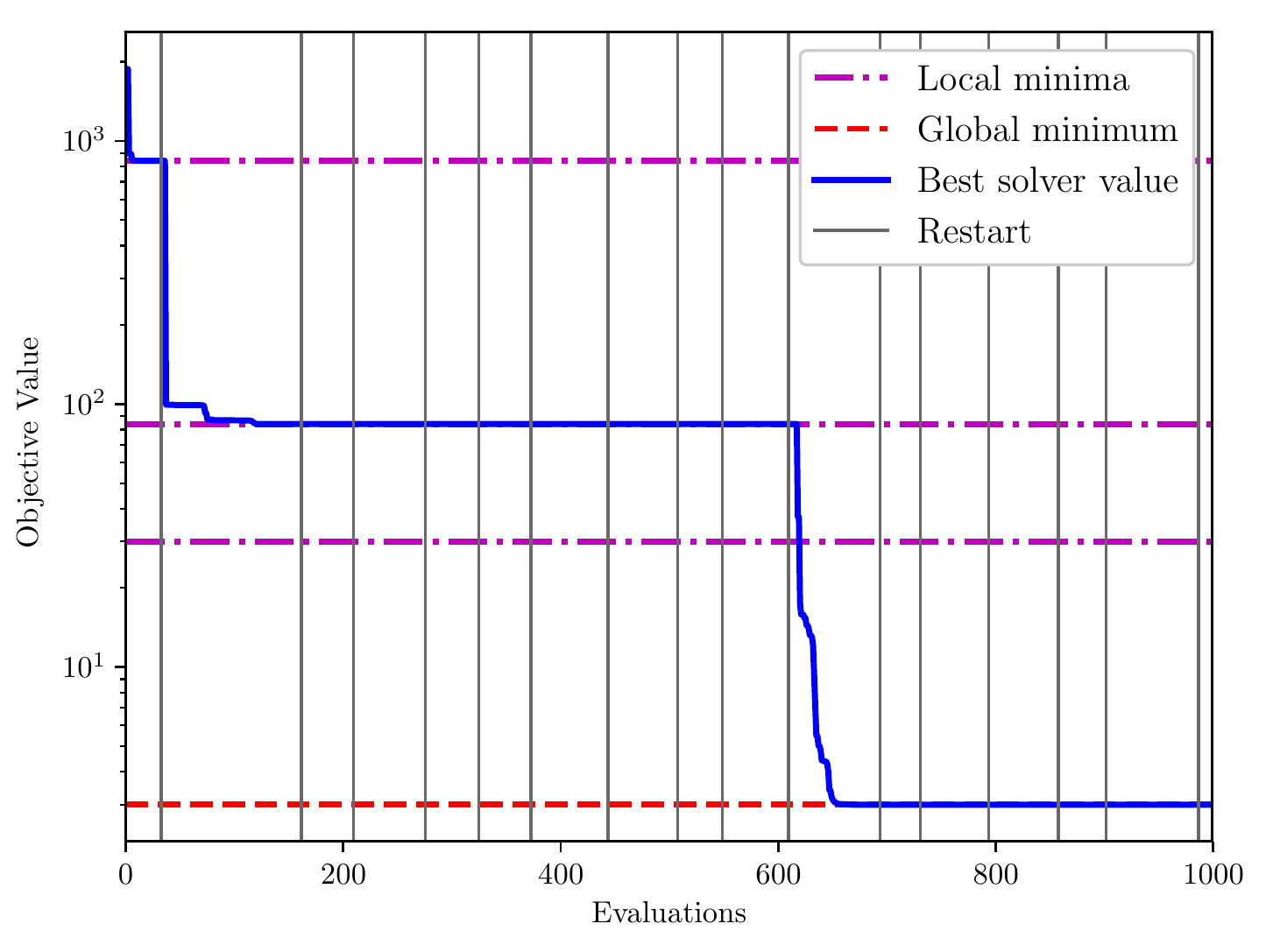}
		\caption{Adaptive restarts: objective value.}
	\end{subfigure}
	~
	\begin{subfigure}[b]{0.48\textwidth}
		\includegraphics[width=\textwidth]{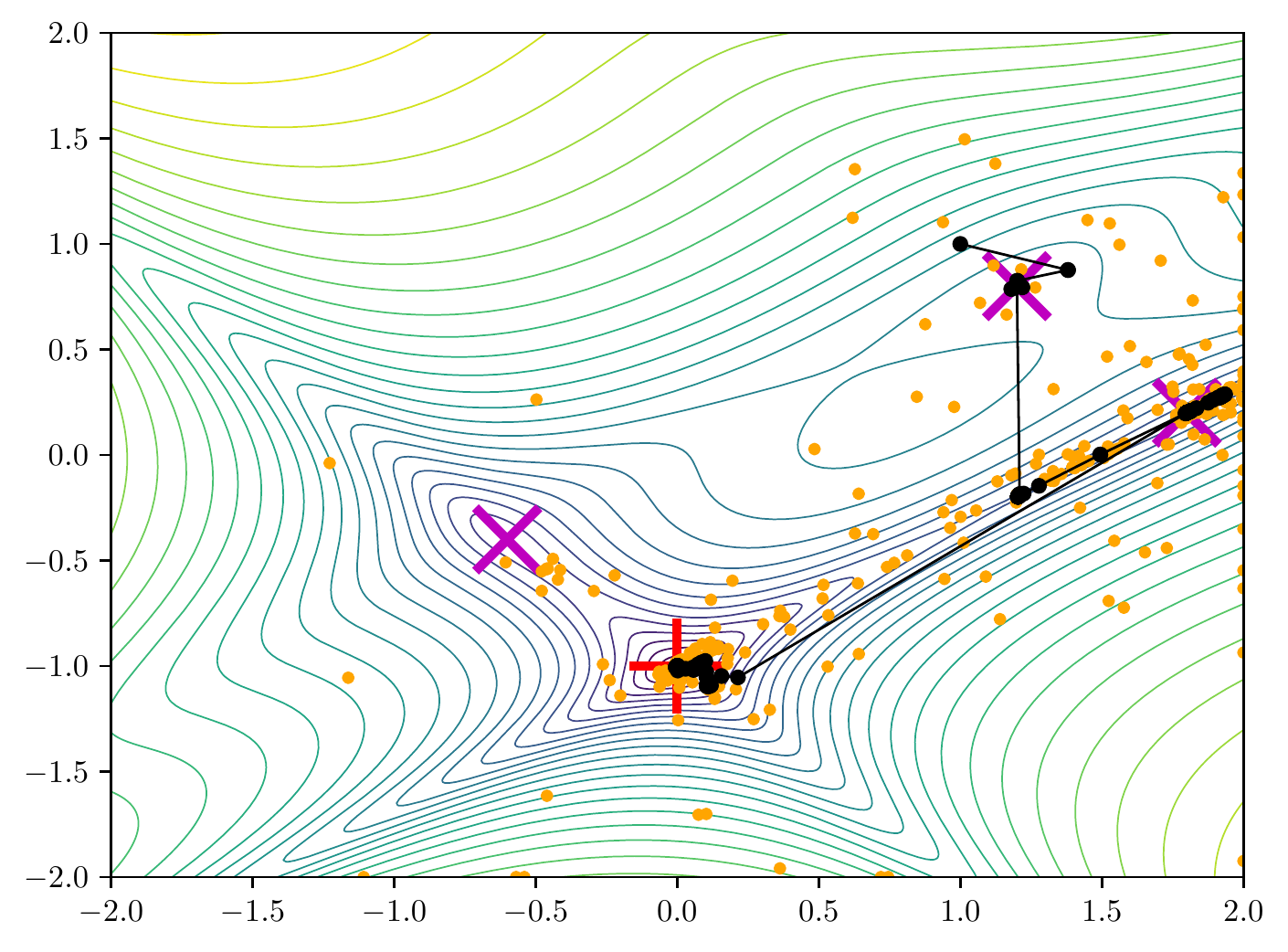}
		\caption{Adaptive restarts: evaluated points.}
	\end{subfigure}
	\caption{Demonstration of Py-BOBYQA with multiple restarts on the Goldstein-Price function. For the evaluated points plot, the pink and red crosses are the local and global minima of the objective, the black points are the iterates of Py-BOBYQA, and the orange points are all other points evaluated by Py-BOBYQA.}
	\label{fig_goldstein_price}
\end{figure}

\paragraph{Goldstein-Price function}
Next, we perform a similar study for the Goldstein-Price function.
This is another well-known GO test problem with $n=2$ defined on $[-2,2]^2$, with four minima (of which one is global).
We start Py-BOBYQA without restarts at $\bx_0=[1,\,1]^{\top}$ with a budget of 1000 evaluations, with the results illustrated in \figref{fig_goldstein_price}.
For each restarts variant---none, fixed and adaptive---we plot the objective value versus number of evaluations (compared to the objective value at the four minima) and the path of the iterates through the feasible set.

Without restarts, Py-BOBYQA converges to the nearest local minimum to $\bx_0$, $\bx=[1.2,\,0.8]^{\top}$, which is also the one with the worst objective value.
If we allow one restart, Py-BOBYQA manages to find the next-best local minimum $\bx=[1.8,\,0.2]^{\top}$.
After more restarts, Py-BOBYQA eventually converges to the global minimum $\bx=[0,\,-1]^{\top}$.
Again, both fixed and adaptive restarts achieve this, but Py-BOBYQA with adaptive restarts explores the feasible region more quickly, and hence requires fewer restarts and evaluations.

\paragraph{High-dimensional Ackley function}
Finally, we demonstrate one advantage of using local solvers over GO solvers, namely improved scalability.
This benefit comes from local solvers not attempting to explore the whole space, and so avoiding the curse of dimensionality.
To show the benefit, we consider Ackley's test function in $\R^n$,
\be f(\bx) = -a\exp\left(-\frac{b\|\bx\|_2}{\sqrt{n}}\right) - \exp\left(\frac{1}{n}\sum_{i=1}^{n}\cos(c x_i)\right) + a + \exp(1). \label{eq_ackley} \ee
As we will be increasing the value of $n$, we make the function less oscillatory by taking $(a,b,c)=(20,0.2,\pi)$ rather than the usual $c=2\pi$, and we use a smaller feasible domain, $[-3,9]^{n}$ and (where necessary) a starting point $\bx_0=[3,\ldots,3]^{\top}$.
The global minimum is still $f(\b{0})=0$.

For this comparison, we take Py-BOBYQA (with adaptive restarts), together with GO solvers DIRECT, CMA-ES, PySOT and SNOBFIT---these are described in \secref{sec_global} below.
We run all solvers for 2000 evaluations, and in \figref{fig_ackley_scaling} we plot the best objective value achieved by each solver for \eqref{eq_ackley} with $n\in[5,10,20,50]$.
We see that Py-BOBYQA achieves the fastest objective decrease of all the solvers, and the relative benefit of Py-BOBYQA compared to the GO solvers is greater as $n$ increases.


\section{Global solvers: a brief description}
\label{sec_global}
\enlargethispage{2em}

Here we briefly describe the global solvers we test against local DFO. We chose a number of solvers that use global models, also constructed from just black-box, possibly noisy and expensive, function values, compared to local models as in Py-BOBYQA, and then a small number of other state-of-the-art deterministic and stochastic global solvers. 
Our aim here is not to provide a comprehensive benchmarking exercise for global optimization solvers, but to make sure we choose state-of-the-art global solvers with a similar remit, to compare with Py-BOBYQA. 
We refer the interested reader to the extensive comparisons and clear descriptions of state-of-the-art derivative-free local and global solvers given in \cite{Rios2013}, and earlier in \cite{Ali2005}, which we also used as a guide for some of our solver choices.

\clearpage

\begin{figure}[t]
	\centering
	\begin{subfigure}[b]{0.48\textwidth}
		\includegraphics[width=\textwidth]{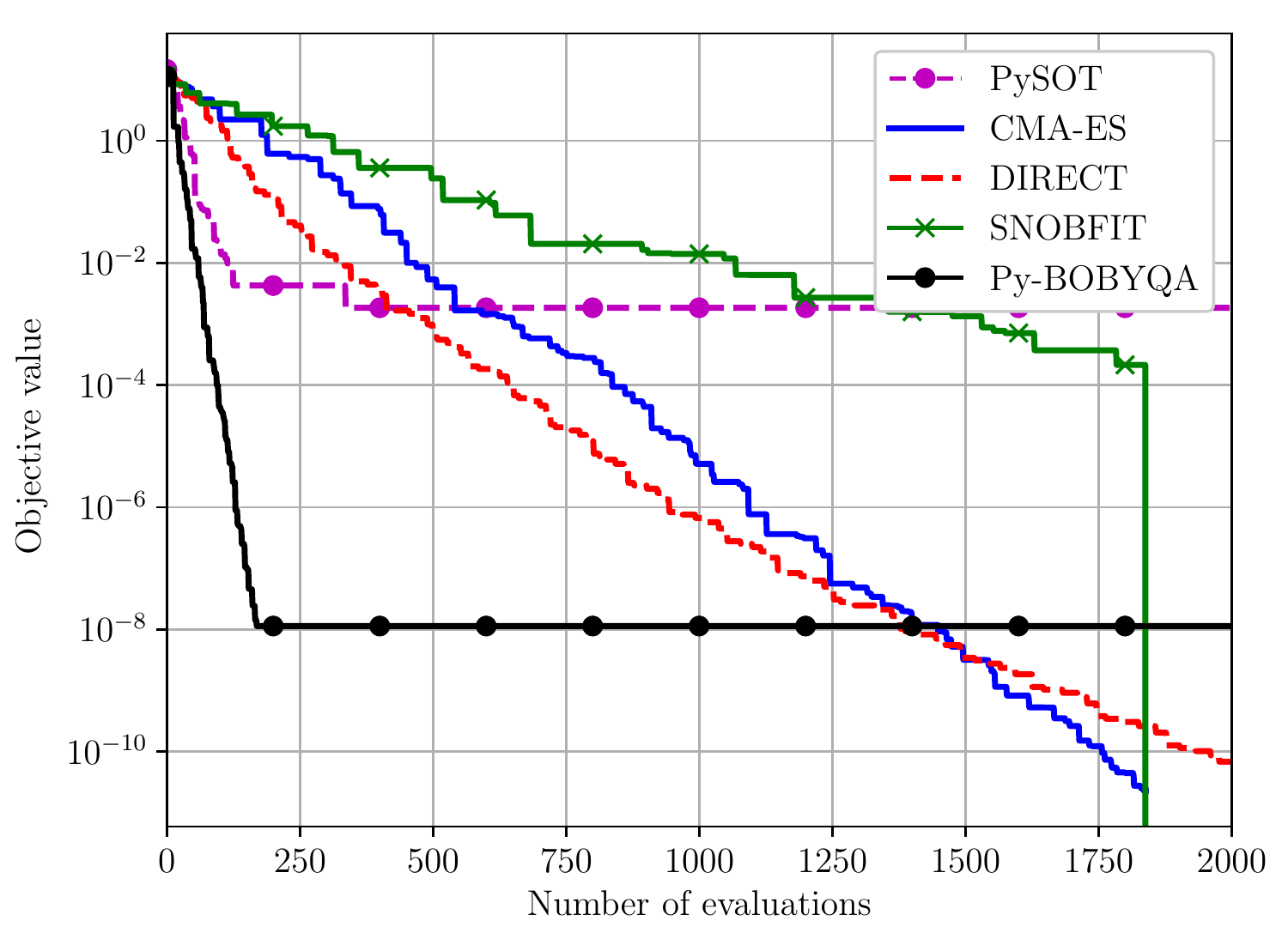}
		\caption{$n=5$}
	\end{subfigure}
	~
	\begin{subfigure}[b]{0.48\textwidth}
		\includegraphics[width=\textwidth]{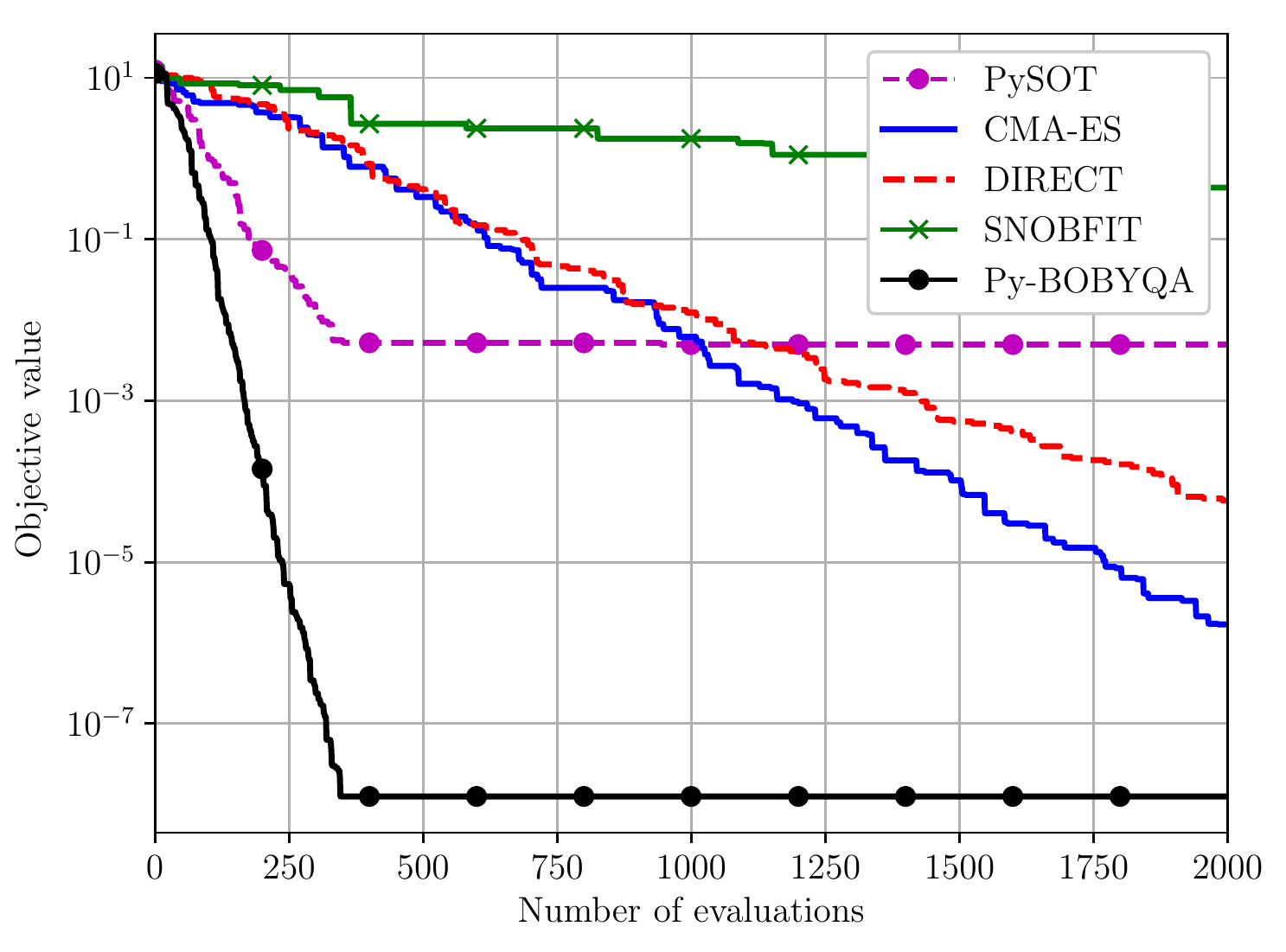}
		\caption{$n=10$}
	\end{subfigure}
	\\
	\begin{subfigure}[b]{0.48\textwidth}
		\includegraphics[width=\textwidth]{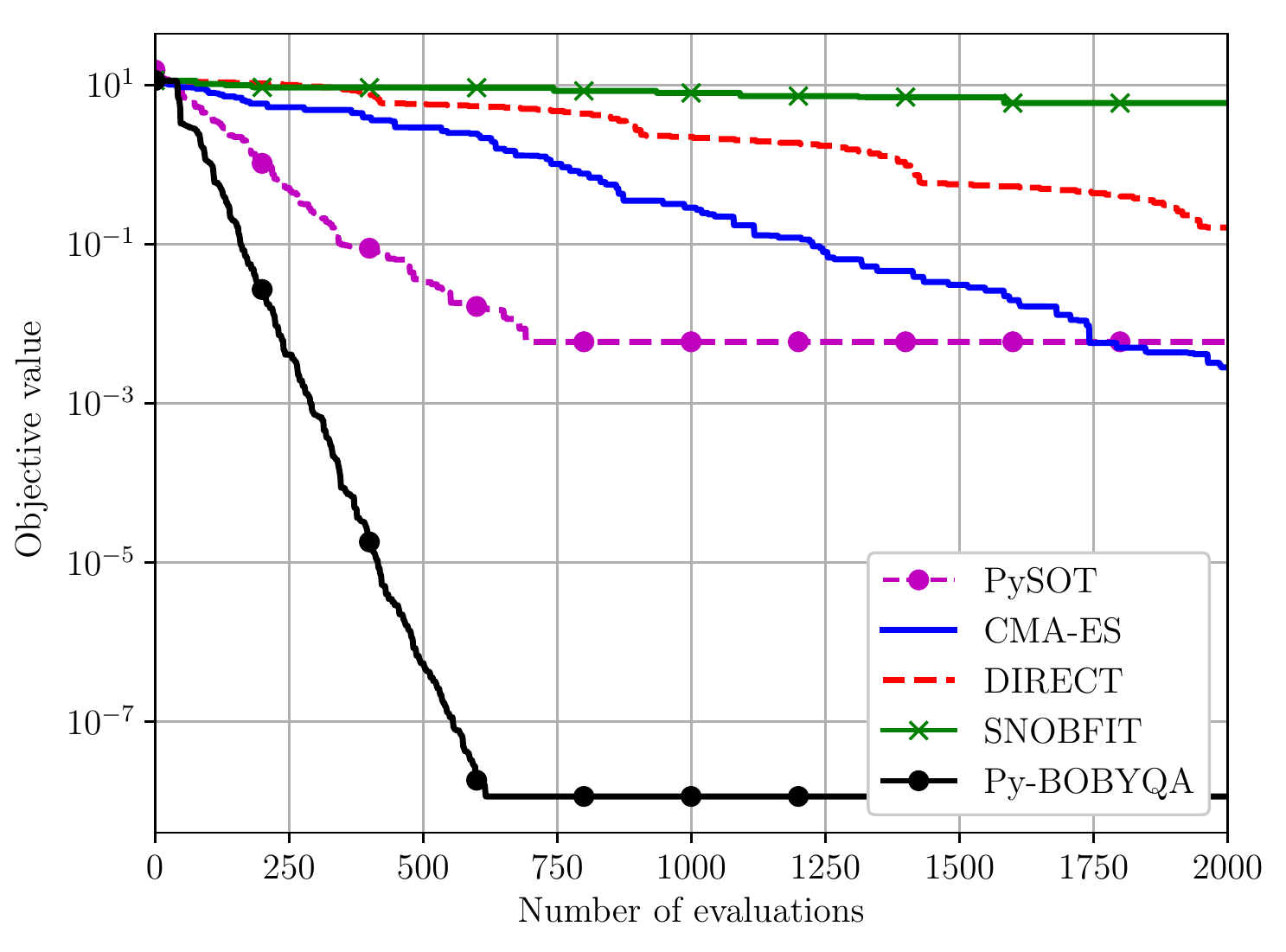}
		\caption{$n=20$}
	\end{subfigure}
	~
	\begin{subfigure}[b]{0.48\textwidth}
		\includegraphics[width=\textwidth]{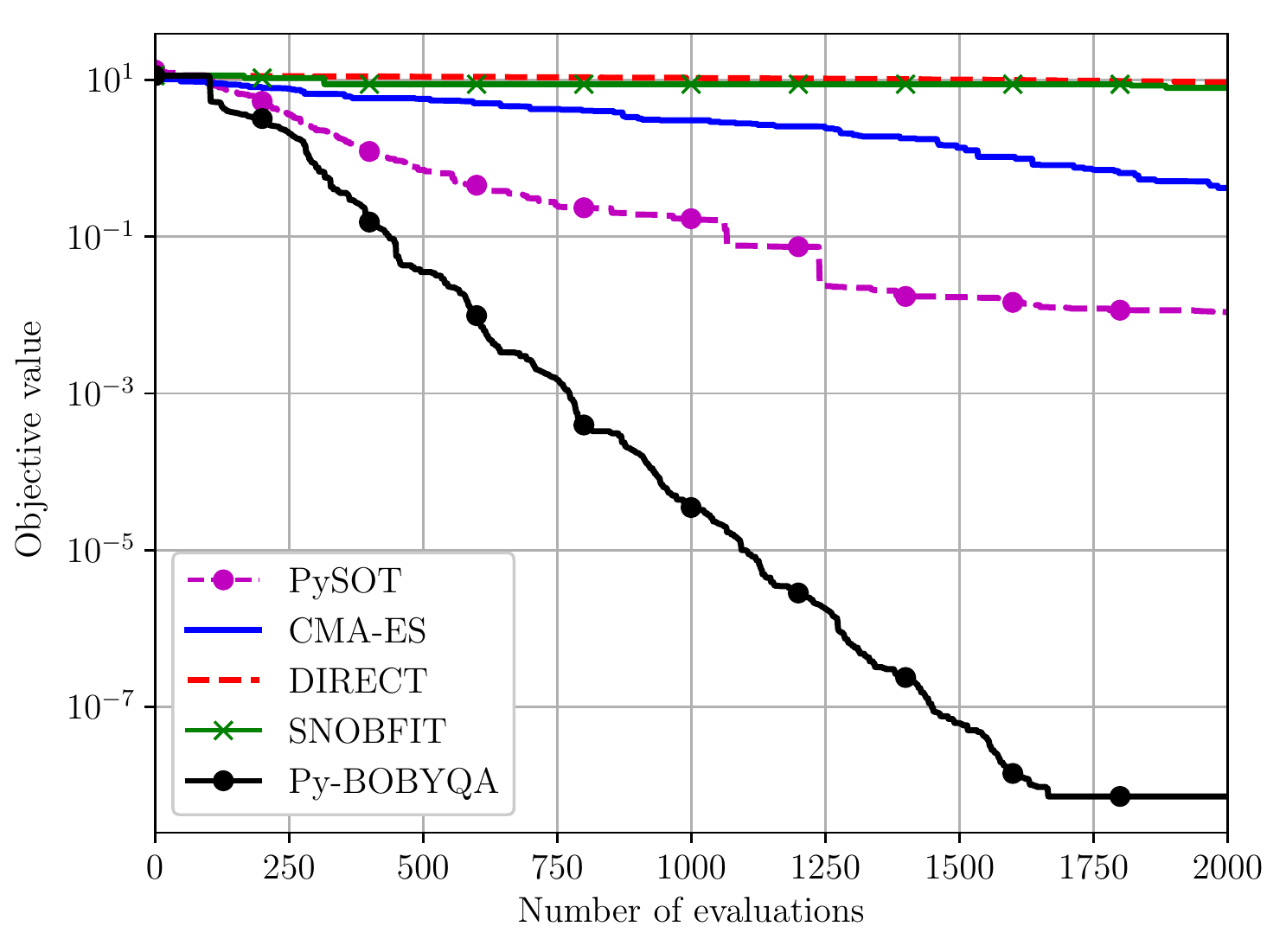}
		\caption{$n=50$}
	\end{subfigure}
	\caption{Comparison of Py-BOBYQA with multiple (adaptive) restarts with GO solvers on the variable-dimension Ackley function.}
	\label{fig_ackley_scaling}
\end{figure}

\paragraph{Bayesian optimization solvers}
Ever since \cite{mockus1978toward, jones1998efficient}, Bayesian methods have been  designed for the global optimization of black-box functions, which are typically  expensive to compute. 
This remit implies that the methods spend computational time and resources to make informed decisions about where to sample the function in the course of their runs. 
As such, Bayesian algorithms attempt to strike a balance between \textit{exploring} unsampled regions, and \textit{exploiting} regions with a promising objective value; see \cite{brochu2010tutorial, shahriari2016taking, fowkes2011bayesian} and the references therein for detailed surveys of this area and connections to global optimization methods and applications. 

As the name indicates, these methods rely on Bayes' theorem and related Bayesian processes in probability, namely, they assume a (global) model (a prior distribution) of the objective $f$ is available together with some point samples  from this distribution, say $\mathcal{D}_k$
at iteration $k$ of the algorithm, which  are iteratively used to determine an a posteriori distribution for $f$ and evolve the current model for $f$. Typical choices of prior distributions  for the objective include Gaussian processes (used in the GPyOpt solver \cite{gpyopt2016}), random forests (PySMAC \cite{Hutter2011}), regression trees,  and tree Parzen estimators (HyperOpt \cite{Bergstra2013}); see 
\cite{rasmussen2006gaussian} for more details and descriptions. The prior/current model of $f$ implies (often closed-form) expressions or choices  for the  sampled 
mean $\mu(\cdot|\mathcal{D}_k)$ and covariance/noise kernel $\sigma(\cdot|\mathcal{D}_k)$.
The noise kernel allows a Gaussian process to capture a wide range of functional/distributional constructions;
various options have been proposed, such as the well-known squared exponential/Gaussian kernel, Mat\' ern $\nu-$kernel, radial basis functions (used in PySOT \cite{pysot2015}), and more \cite{rasmussen2006gaussian}.
 
Using and combining the expectation and covariance (noise kernel) of the distribution, an acquisition 
function  $a(\cdot)$ is formed, which is then employed in order to find the next sampling point, $\bx_{k+1}=\argmin_{\bx}a(\bx)$.  Instead of minimizing directly the model of the objective, as it is done in local DFO for
example, Bayesian optimization methods minimize some measure of `error' or `discrepancy' between the objective and its model.
The latter optimization problem  is  simpler,  and it is inexpensive with respect to evaluations of the objective. The choice and optimization of the acquisition function determines whether the next point $\bx_{k+1}$/model leads to exploring areas with higher uncertainty, or exploiting regions with more promising expected values. 
Various choices of acquisition functions are in use, such as the {\it expected improvement}, {\it probability of improvement},   {\it upper confidence bounds}, and more.
The acquisition function  is updated with the new point that was calculated and its associated function value, iteratively, until convergence or the given budget are reached. Then, in addition to optimization solution information, the acquisition function can be extracted, which contains a distribution and confidence levels for the solution. 

Bayesian optimization is frequently used in machine learning applications, such as for hyperparameter tuning of optimization algorithms for training classifiers \cite{hutter2007automatic, Snoek2012, li2016hyperband}. Ghanbari \& Scheinberg \cite{Ghanbari} compared derivative-free local solvers against Bayesian global solvers on such hyperparameter
tuning problems, with encouraging results. 

\paragraph{A  surrogate optimization algorithm}
Surrogate or response surface optimization methods are similar to Bayesian algorithms \cite{gutmann, fowkes2011bayesian, jones2001taxonomy}, in that they use a global model/surrogate (such as a radial basis function)  to approximate the black-box objective function and then iteratively refine this model. This model is deterministic and constructed by interpolation of function values; a `bumpiness' function measures the discrepancy between the model and true function, and its minimization is used to find the next interpolation point (similarly to the acquisition function in Bayesian methods). Gaussian priors can be viewed as particular radial basis functions, and so the two classes of methods are closely connected \cite{fowkes2011bayesian}.
PySOT \cite{pysot2015} is a state-of-the-art surrogate optimization solver that uses radial basis functions for the response surface; we have included it here due to its similarities to Bayesian optimization, in order to assess the potential benefits of different global models.






\paragraph{A branch and fit method} 

SNOBFIT \cite{Huyer2008} constructs local quadratic models from interpolation of function values, and uses a branching strategy to enhance their global properties. It connects well with Py-BOBYQA in that it is a local method, with improvements that make it suitable for global optimization, hence our choice to include it here. 

\paragraph{An evolutionary algorithm}
Covariance Matrix Adaptation evolutionary strategy (CMA-ES) \cite{Hansen1996, HansenOstermeier2001} is a genetic algorithm, that attempts to find the basin of a global minimizer by iteratively moving and reshaping a Gaussian distribution using samples of objective values (populations).  The point with the best objective value is labelled as the `best point' and the mean of the distribution is updated,
which shifts the mean of the distribution in the direction of the best point. Updating the covariance of the distribution determines its shape and spread, with the algorithm deciding how much old information to keep and how much to update from current best point obtained. In CMA-ES, the updating of the kernel/covariance is a weighted formula allowing three types of contributions (past, rank-one updates along a direction of progress or higher rank updates); the weights in the covariance updating formula are also adjusted in the algorithm. The sophisticated updates are key to the success of this evolutionary strategy and can be viewed as providing some curvature information for the objective. 

\paragraph{The DIRECT algorithm \cite{Jones1993}} 
The well-known DIvide a hyperRECTangle (DIRECT) method is an extension of Lipschitz global optimization strategies  that 
does not require the Lipschitz constant to be known/provided, but uses trial values of this constant chosen heuristically \cite{gutmann}.
In particular, the domain is subdivided into hyper-rectangles and the function is evaluated at the midpoint of each box, and this evaluation is then used to construct Lipschitz-like lower bounds. 
The splitting rules for subdomains are carefully chosen to focus on those that yield the largest potential decrease in the objective, so as to encourage the algorithm towards global searches.

\section{Testing methodology}\label{sec_test_methodology}

\paragraph{Solvers tested}
We compared Py-BOBYQA with a selection of global optimization routines, all implemented or wrapped in Python\footnote{We used Python 2.7 with NumPy 1.12.1 and SciPy 1.0.1.}.
Here, we list each solver tested and any non-standard parameters set---in some cases, these choices depended on whether we were testing an objective with or without noise:
\begin{description}
	\item[\normalfont\textit{CMA-ES v2.5.7 \cite{Hansen1996,cmaes_URL}:}] Used an initial standard deviation $0.2\min_i (b_i-a_i)$, where the feasible region is $\b{a} \leq \bx \leq \b{b}$;
	\item[\normalfont\textit{DIRECT v1.0.1 \cite{Jones1993,pydirect_URL}:}] Used all defaults, except for increasing the maximum number of objective evaluations to be sufficiently large for our testing;
	\item[\normalfont\textit{GPyOpt v1.2.1 \cite{gpyopt2016}:}] Used an initial design of $2n$ random points, and set the \texttt{exact\_feval} flag depending on whether noisy is objective or not; 
	\item[\normalfont\textit{HyperOpt v0.1 \cite{Bergstra2013,HyperOpt_URL}:}] Used a uniform search space over the feasible region, Tree of Parzen Estimators algorithm, with and without providing noise variance for the objective evaluation; 
	\item[\normalfont\textit{PySMAC v0.9.1 \cite{Hutter2011,pysmac_URL}:}] We set the \texttt{deterministic} option according to whether the objective had noise or not; 
	\item[\normalfont\textit{PySOT v0.1.36 \cite{pysot2015}:}] Used a radial basis function surrogate with an initial design of $2n+1$ Latin Hypercube points and a DYCORS \cite{Regis2013} strategy with $100n$ points for selecting new points; 
	\item[\normalfont\textit{SNOBFIT v1.0.0 \cite{Huyer2008,snobfit_URL}:}] Used all default parameters. 
\end{description}
For Py-BOBYQA, we used version 1.1\footnote{Available from \url{https://github.com/numericalalgorithmsgroup/pybobyqa}}, with options: scale feasible region  to $[0,1]^n$; initial trust region radius 0.1 and final trust region radius $10^{-8}$ (in scaled variables); on unsuccessful restarts, increase trust region radius $\Delta_{\mathrm{reset}}$ by a factor\footnote{The constant factor $1.1$ was not decided based on testing or tuning; it was simply considered to be a sensible choice.} of $1.1$; terminate after 10 consecutive unsuccessful restarts or 20 total unsuccessful restarts.

All solvers were allowed a budget of at most $10^4 (n+1)$ objective evaluations for an $n$-dimensional problem, and a runtime of at most 12 hours per problem instance.
In addition, all solvers except DIRECT require a starting point; for each solver, problem, and instance, this was chosen independently and uniformly from the entire feasible region.
Given this, we ran 10 instances of each solver on each problem (where, if applicable, each instance had a different starting point, and realizations of the stochastic noise in the objective).

\paragraph{Performance measurement}
Following the approach from Mor\'e \& Wild \cite{More2009}, we measure solver performance by counting the number of evaluations required to first achieve a sufficient objective reduction; that is, for a solver $\mathcal{S}$, problem $p$ (with true objective $f$ and noisy realized objective $\t{f}$), and accuracy level $\tau\in(0,1)$, we define the number of evaluations required to solve the problem as:
\begin{align}
	N_p(\mathcal{S}; \tau) &\defeq \text{\# objective evaluations taken to find a point $\bx$ satisfying} \nonumber \\
	&\qquad \t{f}(\bx) \leq \mathbb{E}[\t{f}(\bx^*) + \tau(\t{f}(\bx_0) - \t{f}(\bx^*))], \label{eq_reduction_measure}
\end{align}
where $\bx^*$ is the global minimizer of $f$ (also the minimizer of $\mathbb{E}[\t{f}]$).
If \eqref{eq_reduction_measure} was not achieved at all by the solver, we set $N_p(\mathcal{S}; \tau)=\infty$.
We also note that $\t{f}(\bx_0)$ is different for each instance of each problem (and each solver).


To compare solvers, we will use data  profiles \cite{More2009}\footnote{Another possibility is to use performance profiles \cite{More2009},
where the 'difficulty' of the problem is also taken into account.  We have generated these as well, and compared them with data profiles. We found that  they did not bring new/different information in our tests compared to data profiles and so we have not included them here for space reasons.}.
Data profiles measure the proportion of problems solved within a fixed number of evaluations (normalized by problem dimension); that is, for each solver $\mathcal{S}$ on a set of test problems $\mathcal{P}$, we plot the curve
\be d_{\mathcal{S}}(\alpha) \defeq \frac{1}{|\mathcal{P}|} \cdot \left|\left\{p\in\mathcal{P} : N_p(\mathcal{S}; \tau_p) \leq \alpha(n_p+1)\right\}\right|, \qquad \text{for $\alpha\geq 0$,} \label{eq_data_profile} \ee
where $n_p$ is the dimension of problem $p$. 
In our testing, we ran each solver on 10 instances of each problem (i.e.~independent realizations of the starting point and/or random noise in objective evaluations).
For data profiles, we take each problem instance as a separate problem to be `solved'; i.e.~we plot the proportion of the 500 problem instances solved within a given computational budget.

For noisy problems, some of the objective reduction in \eqref{eq_reduction_measure} may be due to a `lucky' realization of the noise, rather than genuine progress towards the objective.
To ensure that our results only capture genuine progress in the optimization, we follow \cite{TR-DFOLS} and floor $\tau$ on a per-problem basis; specifically, for each solve instance of each problem we define:
\be \tau_{\mathrm{crit}}(p) \defeq 10^{\lceil \log_{10}\widehat{\tau}(p)\rceil}, \qquad \text{where} \qquad \widehat{\tau}(p) \defeq \frac{\sigma(\bx^*)}{\mathbb{E}\left[\t{f}(\bx_0)-\t{f}(\bx^*)\right]}, \label{eq_tau_thresh_used} \ee
where $\sigma(\bx^*)$ is the standard deviation of $\t{f}(\bx^*)$.
We then replace $\tau$ in \eqref{eq_reduction_measure} with 
\be \tau_p \defeq \min(\tau_{\mathrm{max}}, \max(\tau_{\mathrm{crit}}(p), \tau)), \ee
where $\tau_{\mathrm{max}}=10^{-1}$ and $\tau$ is the actual accuracy level desired (e.g.~$\tau=10^{-5}$).

Lastly, we note that DIRECT does not allow the specification of a starting point, and instead always selects $\bx_0$ to be the center of the feasible region. 
For 13 problems\footnote{These are: Ackley, Bohachevsky 1, Bohachevsky 2, Camel 3, Cosine Mixture, Exponential, Griewank, Periodic, Powell Quadratic, Rastrigin, Salomon, Schaffer 1, and Schaffer 2.}, this point is the global minimizer, and so DIRECT appears to solve these problems to global optimality after exactly one objective evaluation.
To provide a more reasonable measure of DIRECT's performance, in these instances, we replace $\bx_0$ by $\bx_1$ in \eqref{eq_reduction_measure} and look for the second time this reduction is achieved (i.e.~excluding $\bx_0$ from the count).



\subsection{Test problems}
\paragraph{GO test set} For our comparison, we used the collection of 50 problems from Ali, Khompatraporn, and Zabinsky \cite{Ali2005}. 
To make the testing more difficult, we increased the dimension of 6 problems to $n\in[25,50]$, where they had $n=10$ in \cite{Ali2005}.
In addition, we reduce the dimension of the Griewank test function from $n=10$ in \cite{Ali2005} to $n=5$, which again makes it more difficult to solve (see the discussion in \cite{Locatelli2003}).
The collection of test problems is given in \tabref{tab_final_global_problems} in \appref{app_test_problems}.

In addition, we optionally modified the problems by the addition of stochastic noise; that is, each solver has access to $\t{f}(\bx)\approx f(\bx)$ for each $\bx$.
We tested the following noise models:
\begin{itemize}
	\item Smooth (no noise): $\t{f}(\bx) = f(\bx)$
	\item Multiplicative Gaussian noise: $\t{f}(\bx) = f(\bx) (1 + \sigma \epsilon)$, where $\epsilon\sim N(0,1)$ i.i.d.~for each $\bx$ and each solver; and
	\item Additive Gaussian noise: $\t{f}(\bx) = f(\bx) + \sigma\epsilon$, where $\epsilon\sim N(0,1)$ i.i.d.~for each $\bx$ and each solver.
\end{itemize}
For our testing, we used a noise level of $\sigma=10^{-2}$.

\paragraph{Hyperparameter tuning (MNIST)}
We also considered a case study of hyperparameter tuning for machine learning, for which developers have often used Bayesian optimization solvers.
Here, we followed the approach from Snoek, Larochelle, and Adams \cite{Snoek2012}, and consider the problem of training a multi-class logistic classifier for MNIST\footnote{The multi-class classifier came from by training 10 one-versus-all classifiers, with prediction by selecting the class with highest probability. The training set had 60,000 images, the test set had 10,000 images.} using stochastic gradient descent with a fixed learning rate.
The 4 hyperparameters we optimized in our testing were:
\begin{itemize}
	\item $\log_{10}$ of learning rate, in $[-10, 0]$;
	\item $\log_{10}$ of $\ell_2$ regularization parameter, in $[-10, 0]$;
	\item Number of epochs, in $[1, 100]$; and
	\item Batch size, in $[1, 2000]$.
\end{itemize}
We converted real-valued inputs for the last two parameters to integer values by rounding up to the nearest integer.
The objective function we selected was the (top-1) error rate on the test set.
We ran the solvers as per \secref{sec_test_methodology} (i.e.~selecting a random feasible $\bx_0$ if needed), for a maximum of $100(n+1)=500$ objective evaluations or 36 hours.
For each solver, we ran 30 instances of the problem, and compared performance using data and performance profiles, in the same way as above.

{\it Selecting the `optimum' value for MNIST:}
To construct data and performance profiles, we need estimates of $\mathbb{E}[\t{f}(\bx_0)]$, $\mathbb{E}[\t{f}(\bx^*)]$, and $\sigma(\bx^*)$ in order to evaluate \eqref{eq_reduction_measure}.
However, unlike before, we do not have access to the true values for our objective, and so we need to provide sensible estimates of these values.
For $\mathbb{E}[\t{f}(\bx_0)]$, we used the value of $\t{f}(\bx_0)$ from each instance of each solver.
To determine $\bx^*$, and hence $\mathbb{E}[\t{f}(\bx^*)]$ and $\sigma(\bx^*)$, we did the following:
\begin{enumerate}
	\item Take the 20 points found by all instances of all solvers with the smallest (noisy) objective value;
	\item Evaluate the (noisy) objective $\t{f}$ independently 30 times at each of these points $\bx$, and use this to estimate $\mathbb{E}[\t{f}(\bx)]$ and $\sigma(\bx)$;
	\item Select as $\bx^*$ the point which gave the smallest value of $\mathbb{E}[\t{f}(\bx)]$.
\end{enumerate}
This process gave the point $\bx^*=[-5.2706, -2.2505, 100, 411.02]^{\top}$, with estimates $\mathbb{E}[\t{f}(\bx^*)] = 0.0818067$ and $\sigma(\bx^*)=1.01388\times 10^{-3}$.
This equates to a learning rate of $5.36\times 10^{-6}$, $\ell_2$ regularization parameter $5.62\times 10^{-3}$, 100 epochs and batch size 412, yielding a test set error rate of $8.1\%\pm0.1\%$.

\section{Numerical results}
\label{sec_numerics}

\subsection{Numerical results on the GO test set}

\paragraph{Selecting the best-performing Py-BOBYQA variant}

We compare various Py-BOBYQA variants on the GO test set, to find which variant is most efficient for different accuracy/budget 
regimes.  We divided our testing based on the size of the interpolation set, namely, $p=2n+1$ and $p=(n+1)(n+2)/2$, which corresponds to the case of an inexact local quadratic (as if having a diagonal Hessian) and a full quadratic model $m_k$ in \eqref{eq_quadratic_models}, respectively\footnote{The case $p=n+1$ when $m_k$ is linear (and hence contains no curvature information) does not yield competitive variants for general objectives \cite{CSV, Cartis2017a}.}, and we will compare these two sets of results to justify our default choices of $p$ for Py-BOBYQA given in \secref{sec_algo_framework}. Clearly, building and using an inexact quadratic local model has a reduced evaluation and computational cost (compared to a full quadratic model), but we need to assess whether the local efficiency translates into globally good performance, across iterations; hence the need to look at the two types of local models.
For both choices of interpolation models, we compare the different types of restarts allowed in Py-BOBYQA---{\it hard} and {\it soft} (choice of interpolation set at a restart); {\it fixed} and {\it adaptive} (size of trust region on restart). 



{\it Low accuracy/budget} ($\tau=10^{-2}$ in \eqref{eq_reduction_measure} and budget of $100(n+1)$ evaluations) Our tests show that in this accuracy/budget regime, there is no substantial difference between Py-BOBYQA variants, for the different interpolation models and restart strategies. We note a (small but clear) reduction in performance when no restarts are performed, and a (slightly) superior performance for soft restarts 
with fixed and adaptive choices of trust region radius on restarts, 
for both values of the interpolation set size $p$. We also note that the evaluation cost of building  full quadratic local models
can be seen in the early iterations of the algorithms, when the variants that use reduced quadratics can already make progress, while the full size models are still being built. 
We delegate the data profiles for this case to \appref{app_pybobyqa_extra_results}, where we illustrate also the different noise models and compare the two interpolation cases.

\begin{figure}[t]
	\centering
	\begin{subfigure}[b]{0.48\textwidth}
		\includegraphics[width=\textwidth]{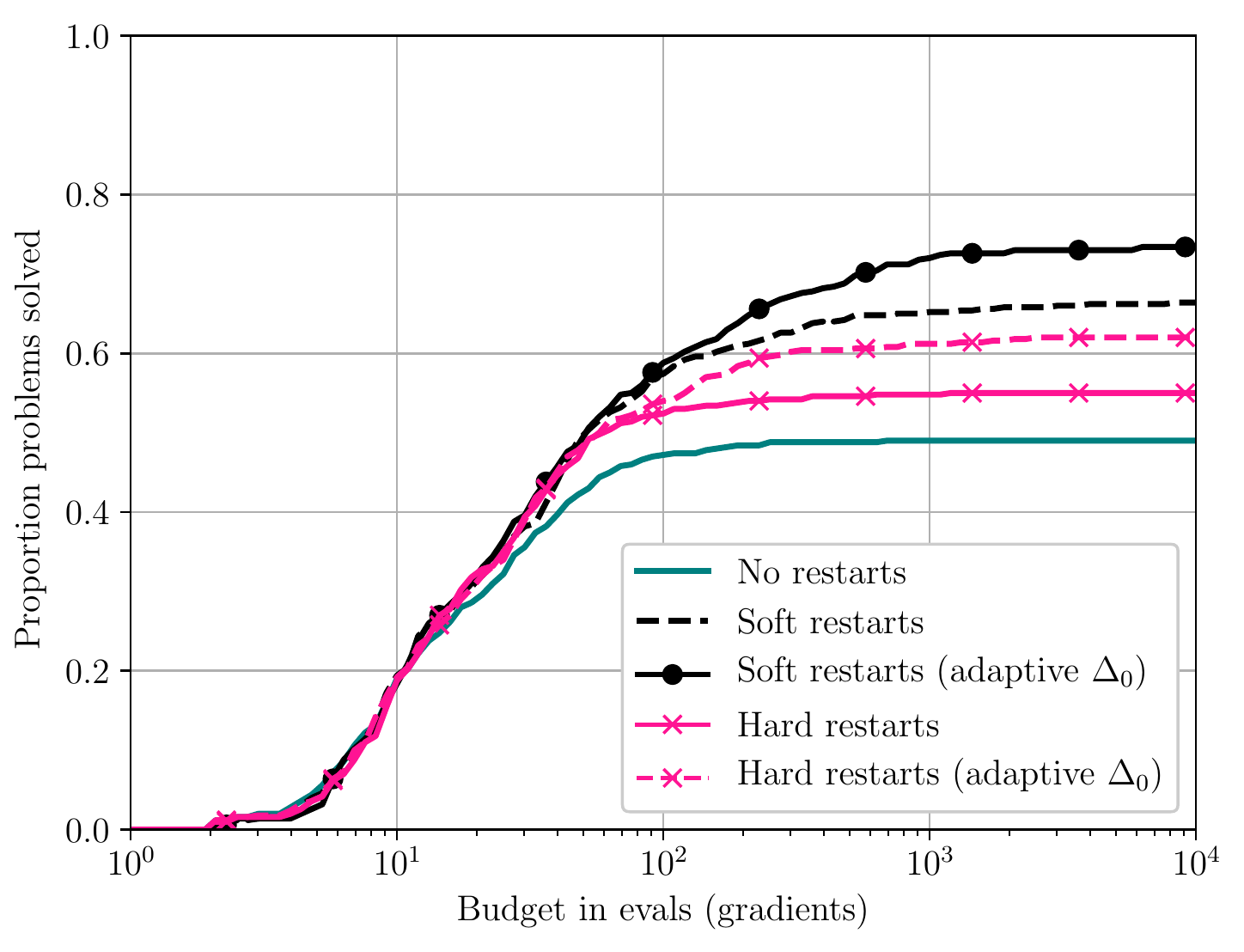}
		\caption{Partial interpolation models  ($p=2n+1$).}
	\end{subfigure}
	~
	\begin{subfigure}[b]{0.48\textwidth}
		\includegraphics[width=\textwidth]{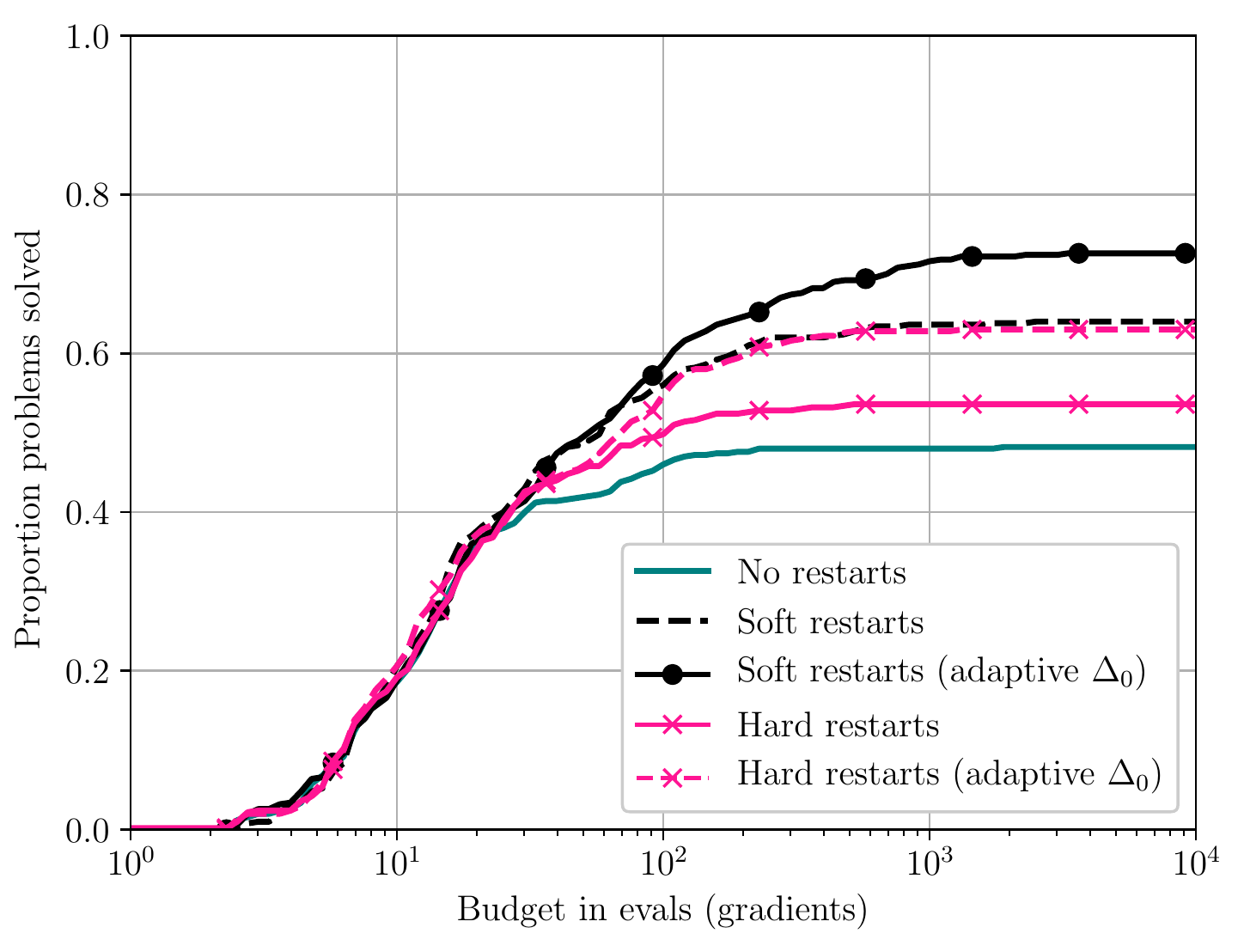}
		\caption{Full  interpolation models  ($p=\mathcal{O}(n^2)$).}
	\end{subfigure}
	\caption{Comparison of Py-BOBYQA variants on the smooth GO test set.}
	\label{pybobyqa-smooth-high}
\end{figure}

\begin{figure}[t]
	\begin{subfigure}[b]{0.48\textwidth}
		\includegraphics[width=\textwidth]{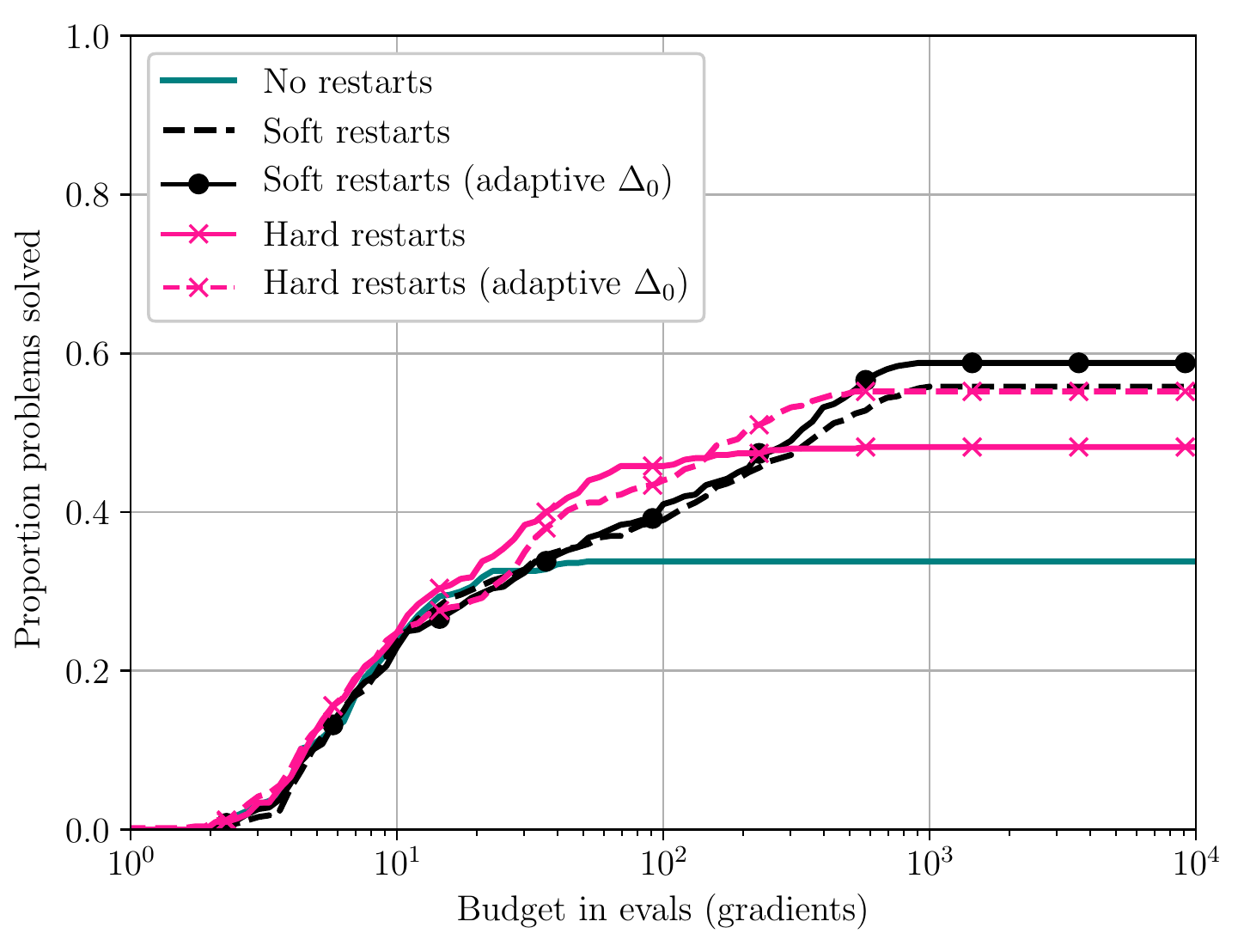}
		\caption{Partial  interpolation models  ($p=2n+1$).}
	\end{subfigure}
	~
	\begin{subfigure}[b]{0.48\textwidth}
		\includegraphics[width=\textwidth]{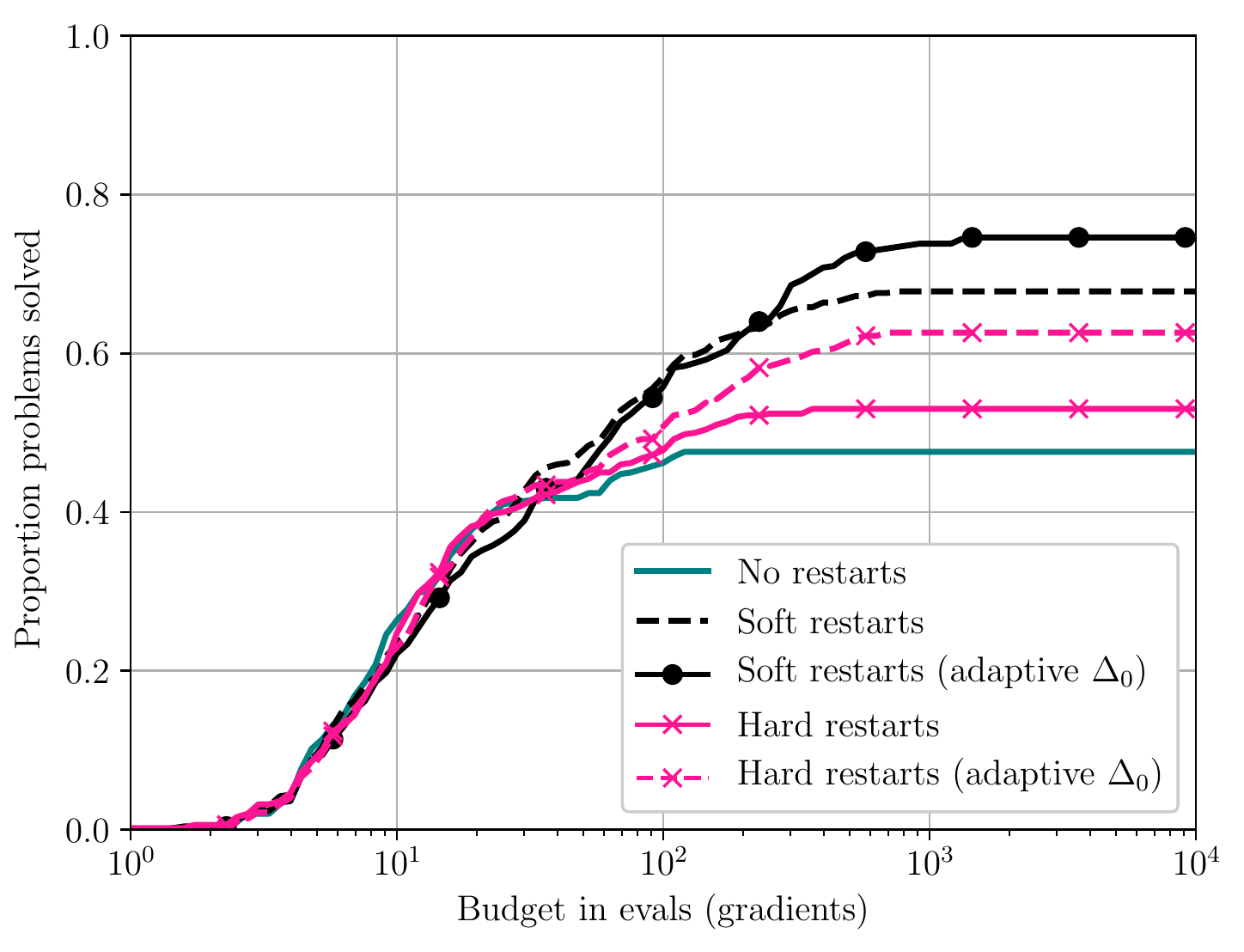}
		\caption{Full  interpolation models  ($p=\mathcal{O}(n^2)$).}	
	\end{subfigure}
	\caption{Comparison of Py-BOBYQA variants on the GO test set with multiplicative noise.}
	\label{pybobyqa-mult-high}
\end{figure}

\begin{figure}[t]
	\begin{subfigure}[b]{0.48\textwidth}
		\includegraphics[width=\textwidth]{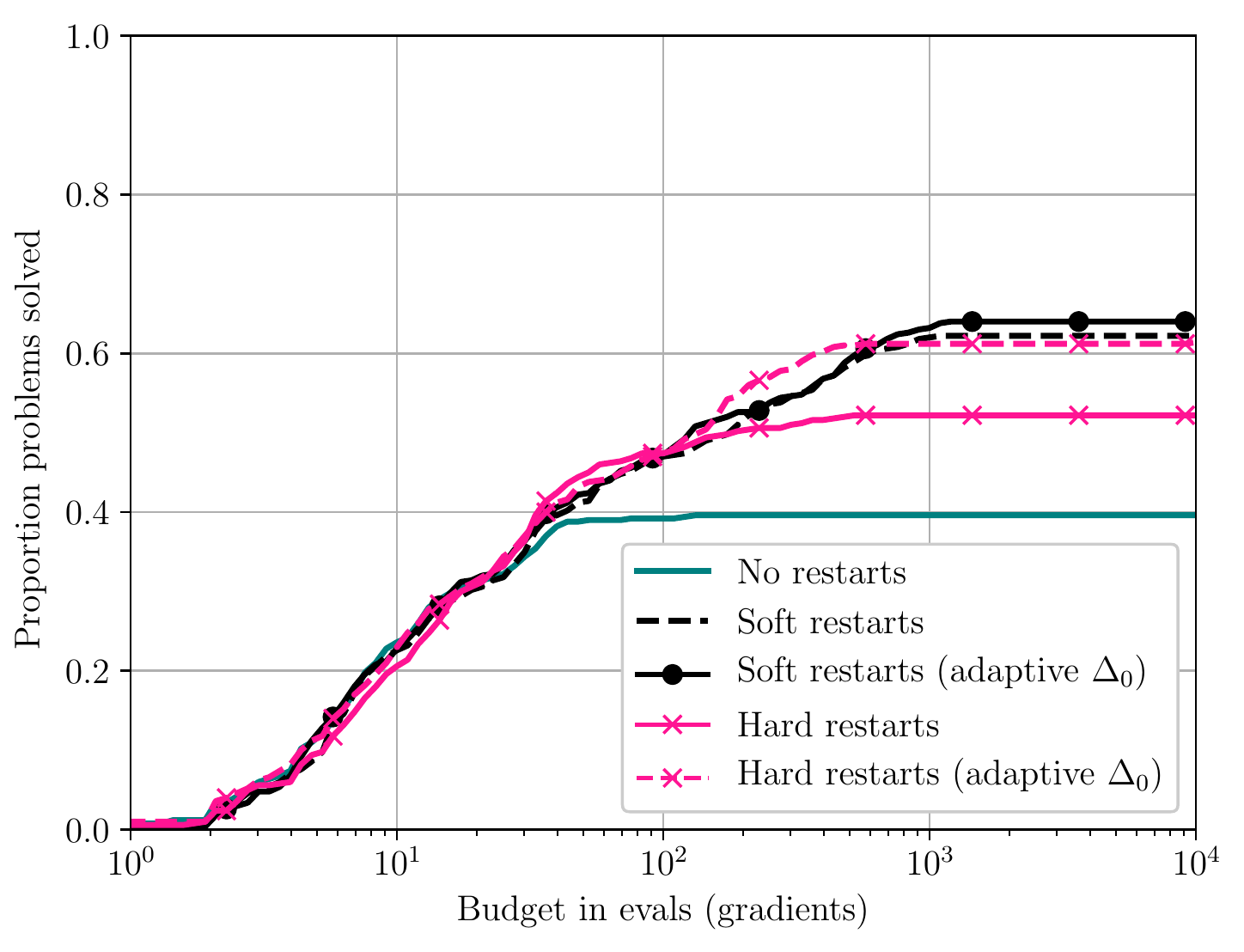}
		\caption{Partial  interpolation models  ($p=2n+1$).}
	\end{subfigure}
	~
	\begin{subfigure}[b]{0.48\textwidth}
		\includegraphics[width=\textwidth]{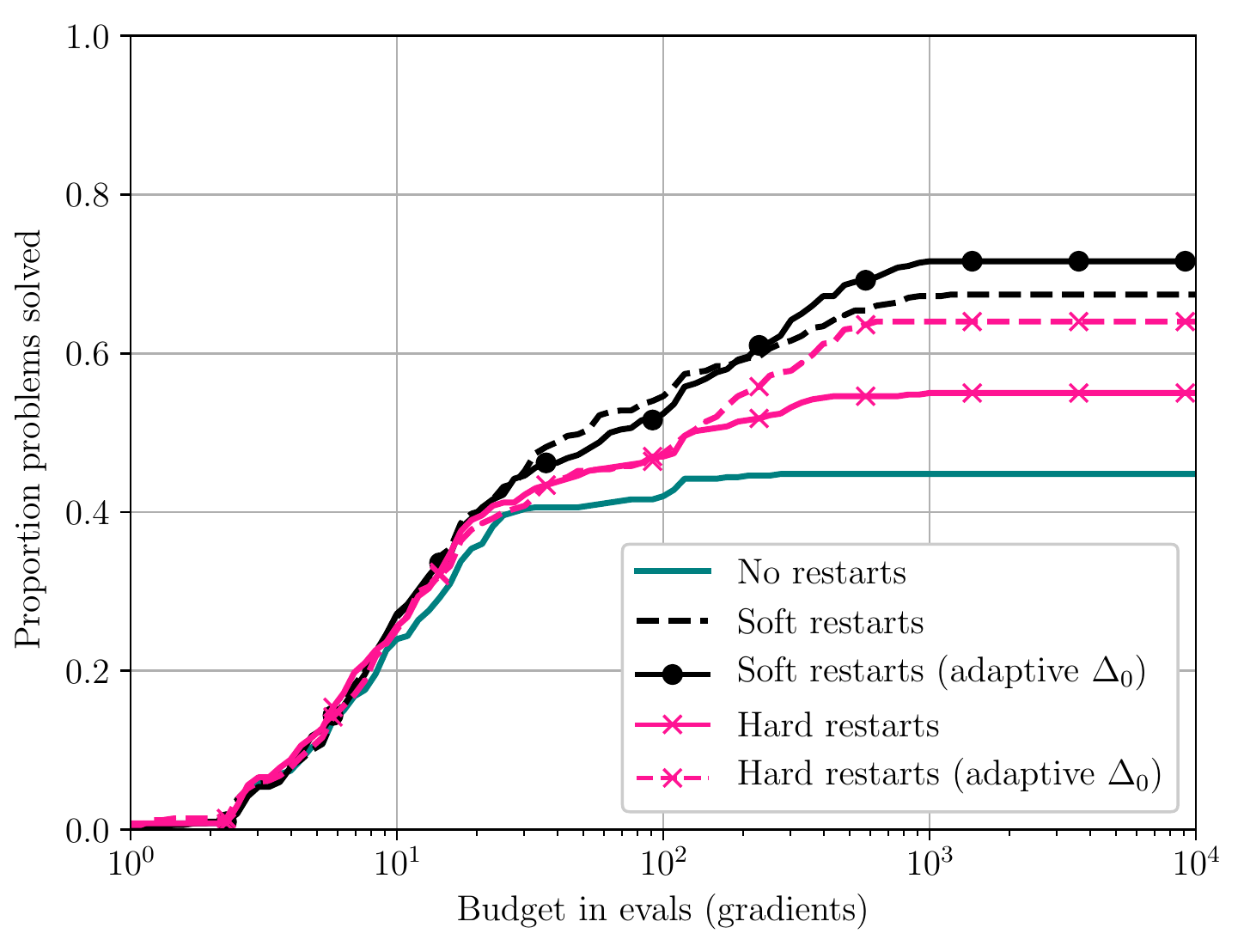}
		\caption{Full  interpolation models  ($p=\mathcal{O}(n^2)$).}	
	\end{subfigure}
	\caption{Comparison of Py-BOBYQA variants  on GO test set with additive noise.}
	\label{pybobyqa-add-high}
\end{figure}

{\it High accuracy/high budget} ($\tau=10^{-5}$ in \eqref{eq_reduction_measure} and full budget of $10^4 (n+1)$ evaluations)
For both interpolation models, and different noise models (including no noise), Figures \ref{pybobyqa-smooth-high}--\ref{pybobyqa-add-high} show the continuously improving performance of some Py-BOBYQA variants with restarts, while the no-restarts variant for example, can no longer improve its performance beyond the achievements in the low accuracy regime, despite the substantially larger budget. 
In particular, the no-restarts Py-BOBYQA is the worst performing variant, failing to solve the largest amount of problems  within the given (large) budget. 
In terms of restarts, the soft restart variant with adaptive restart radius (adaptive $\Delta_0$) is the most efficient and economical across all models and restart types, followed by standard soft restarts (with fixed radius size) and matched sometimes by hard restarts with adaptive radius size (such as for $p=2n+1$ and noisy problems; and $p=(n+1)(n+2)/2$ and no noise). 
To contrast the two interpolation model types, we see that Py-BOBYQA with full quadratic models ($p=(n+1)(n+2)/2$), soft restarts and adaptive radius on restarts gives the best performance overall in our GO tests, and is particularly effective on noisy problems. 
We note that the superior performance of this Py-BOBYQA variant is matched by Py-BOBYQA with inexact quadratic models ($p=2n+1$), also with soft and adaptive restarts on smooth problems, which is important as the computational cost of Py-BOBYQA ($p=2n+1$) is less than that of Py-BOBYQA with full quadratic models, as the linear interpolation system that is solved to find the local model is much smaller in dimension. 
We also note that in the large budget regime, we no longer `pay the price' of the more expensive model set up. 




Overall, we conclude that Py-BOBYQA with adaptive restart radius outperforms the other variants (no restarts and fixed restart radius), and using $p=2n+1$ interpolation points for smooth problems and $p=(n+1)(n+2)/2$ for noisy problems are suitable choices.
Although Py-BOBYQA with no restarts is able to solve some problems to high accuracy---possibly because having well-spaced interpolation points means the region of attraction of the global minimum is found in the main iteration---but, reflecting the examples in \secref{sec_restarts_demo}, many more problems can be solved to high accuracy with adaptive restarts.

\begin{figure}[t]
	\centering
	\begin{subfigure}[b]{0.48\textwidth}
		\includegraphics[width=\textwidth]{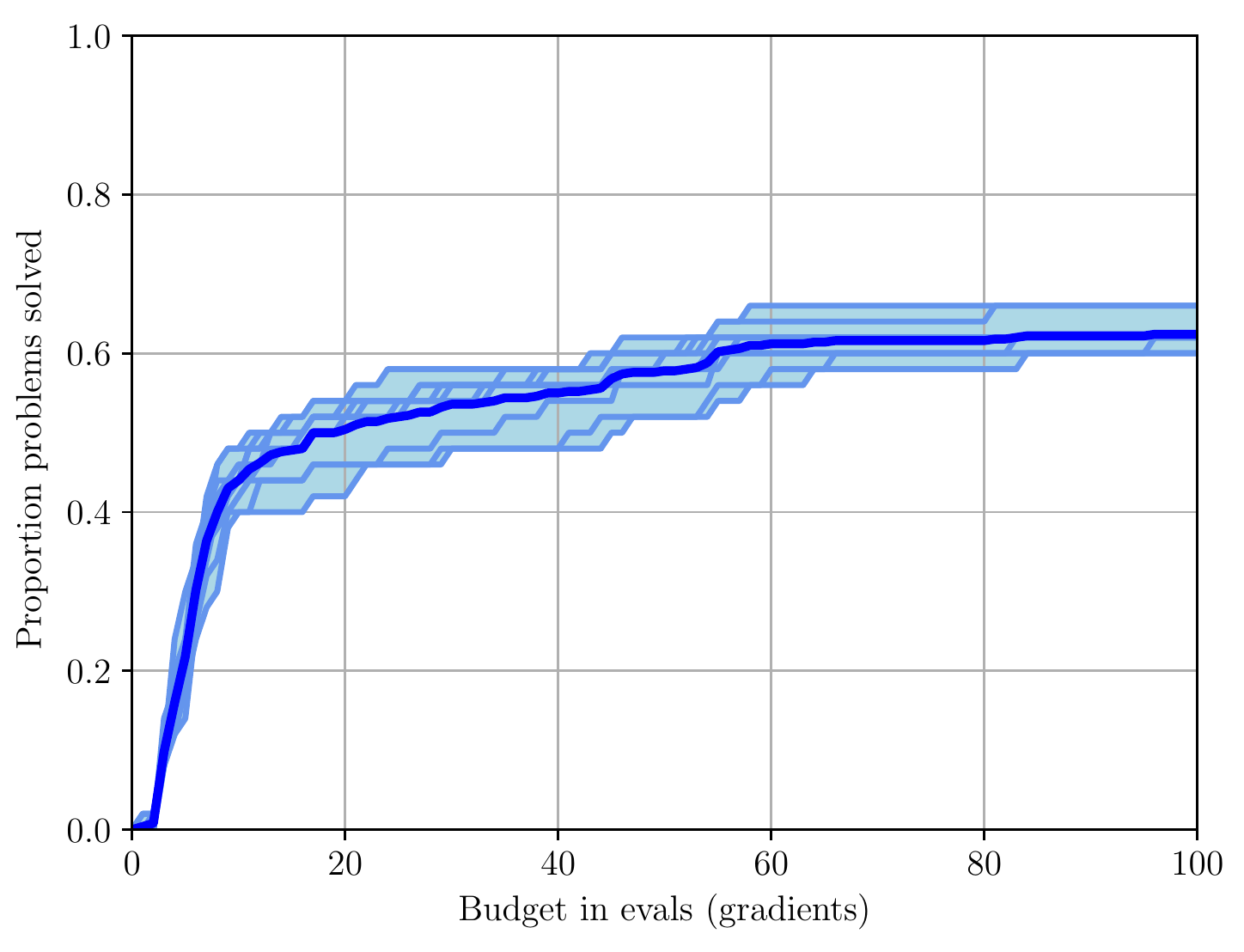}
		\caption{Fixed restart radius, low accuracy/budget.}
	\end{subfigure}
	~
	\begin{subfigure}[b]{0.48\textwidth}
		\includegraphics[width=\textwidth]{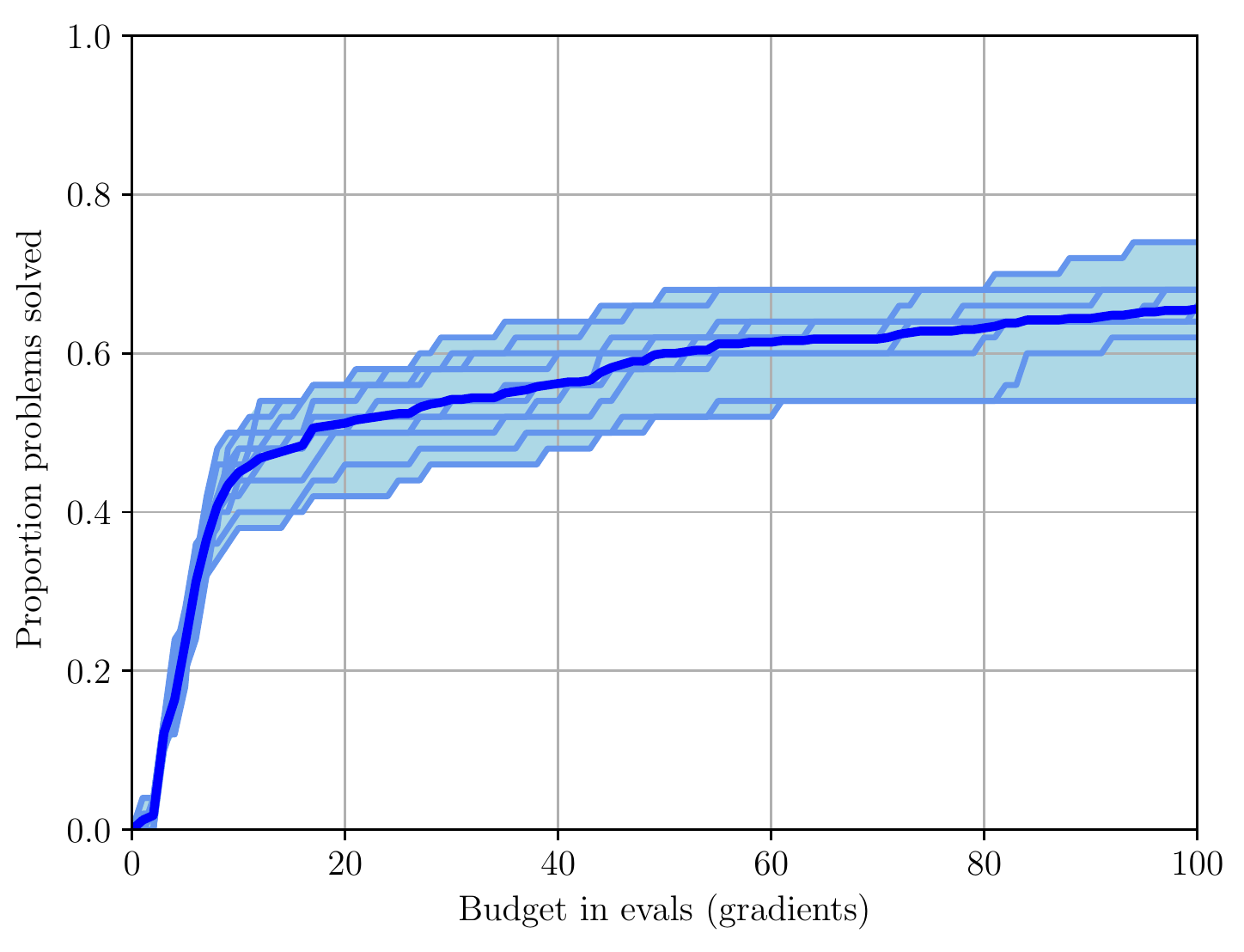}
		\caption{Adaptive restart radius, low accuracy/budget.}
	\end{subfigure}
	\\
	\begin{subfigure}[b]{0.48\textwidth}
		\includegraphics[width=\textwidth]{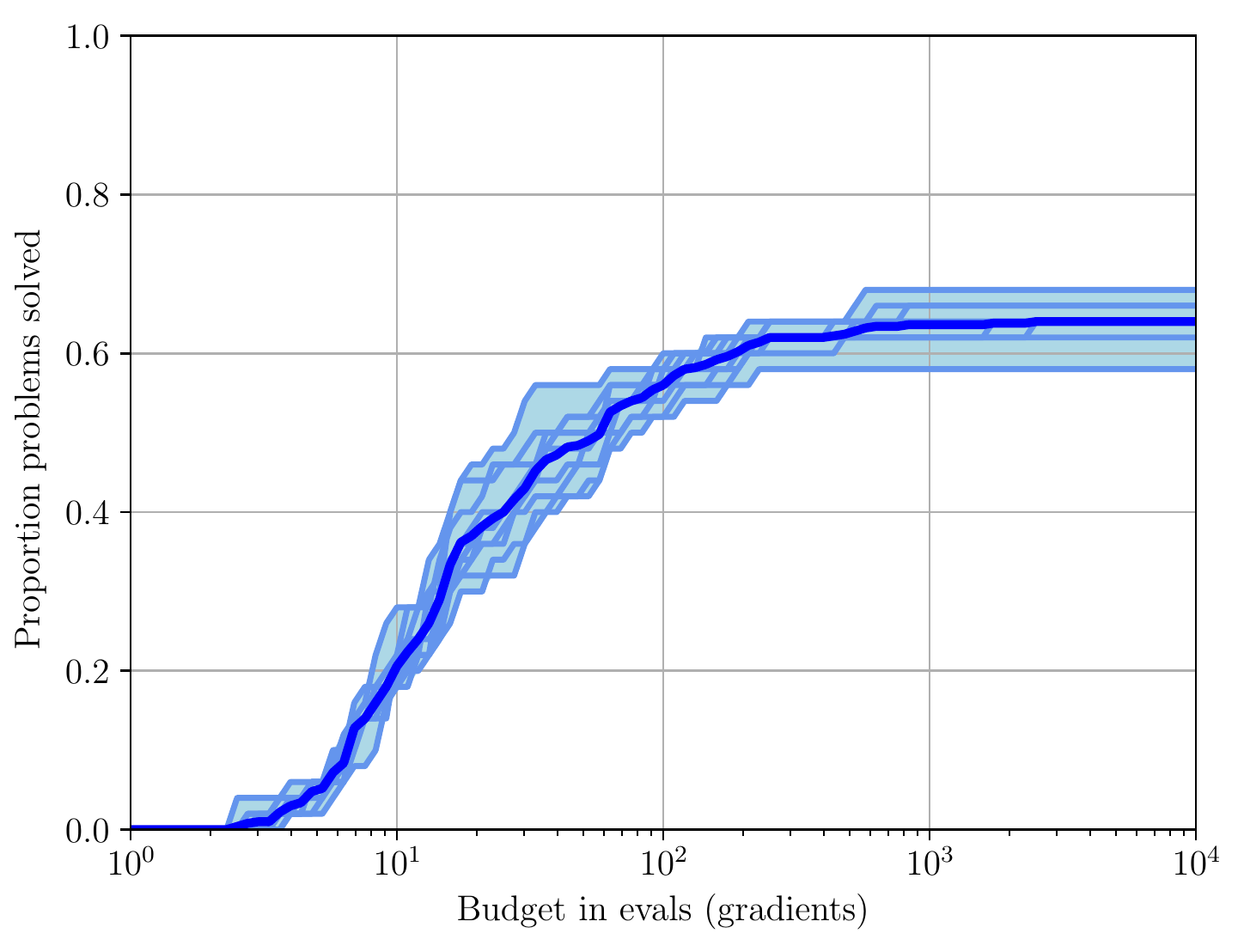}
		\caption{Fixed restart radius, high accuracy/budget.}
	\end{subfigure}
	~
	\begin{subfigure}[b]{0.48\textwidth}
		\includegraphics[width=\textwidth]{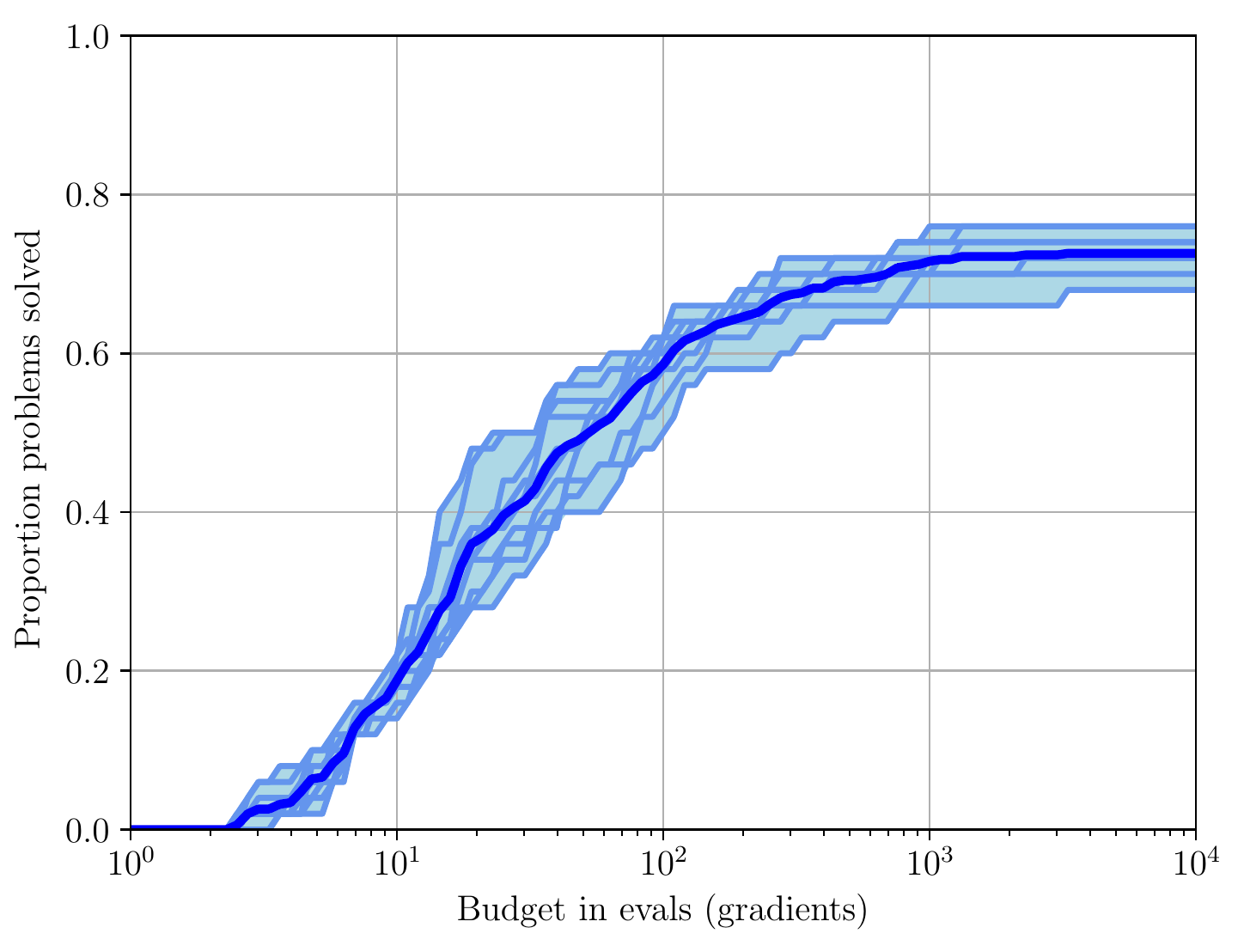}
		\caption{Adaptive restart radius, high accuracy/budget.}
	\end{subfigure}
	\caption{Demonstration of impact of random starting point for Py-BOBYQA ($p=\bigO(n^2)$ with soft restarts). Low accuracy and budget are $\tau=10^{-2}$ and $10^2 (n+1)$ evaluations (for an $n$-dimensional problem) respectively; high accuracy and budget are $\tau=10^{-5}$ and $10^4 (n+1)$ evaluations respectively. The dark lines are the average of all runs.}
	\label{fig_pybobyqa_variance_nsq}
\end{figure}

\paragraph{Py-BOBYQA variance test} \label{sec_pybobyqa_variance}
As Py-BOBYQA is a local solver, and as we are using (ten) randomly chosen starting points per problem instance, we investigate here the variation in performance over  ten different runs of Py-BOBYQA over the test set with different starting points each time. 
This illustrates the variance in our previous results, and how much the performance is influenced by a particular starting point, a key question for local solvers. 
Given the good performance of Py-BOBYQA with full quadratic models ($p=\bigO(n^2)$) and soft restarts across all test problem variants, we are focusing on this variant, with fixed and adaptive restart trust-region radius on restarts. 
We are only looking at smooth problems (with no noise), to make sure the only randomness present is due to choice of starting point. 
In \figref{fig_pybobyqa_variance_nsq}, we show separate data profiles for the ten runs of Py-BOBYQA (i.e.~each profile shows one instance of each problem with a random starting point), and the average data profile.


The variation in starting point usually leads to between 5\%--10\% variation in performance in the second half/end of given budget, 
except in the low accuracy regime and adaptive trust region radius where we see a variation in performance of about $20\%$
at the end of the budget. For high accuracy/large budget, the variation in performance is such that if we compare with \figref{pybobyqa-smooth-high}, we see that all/most adaptive radius runs are better in performance than the fixed radius variant. For low accuracy,
the situation is less clear and performances of different variants are less distinctive, which is no surprise given the average plots in \figref{pybobyqa-low}. 

We also note that these conclusions are similar for Py-BOBYQA with $p=2n+1$ interpolation points; see \figref{fig_pybobyqa_variance_2np1} in \appref{app_pybobyqa_extra_results}.


\paragraph{Selecting the best-performing Bayesian and surrogate solver}

We first compare the Bayesian solvers (GPyOpt, HyperOpt, and PySMAC)
and surrogate algorithm PySOT for smooth and noisy GO test set.

\begin{figure}[t]
	\centering
	\begin{subfigure}[b]{0.48\textwidth}
		\includegraphics[width=\textwidth]{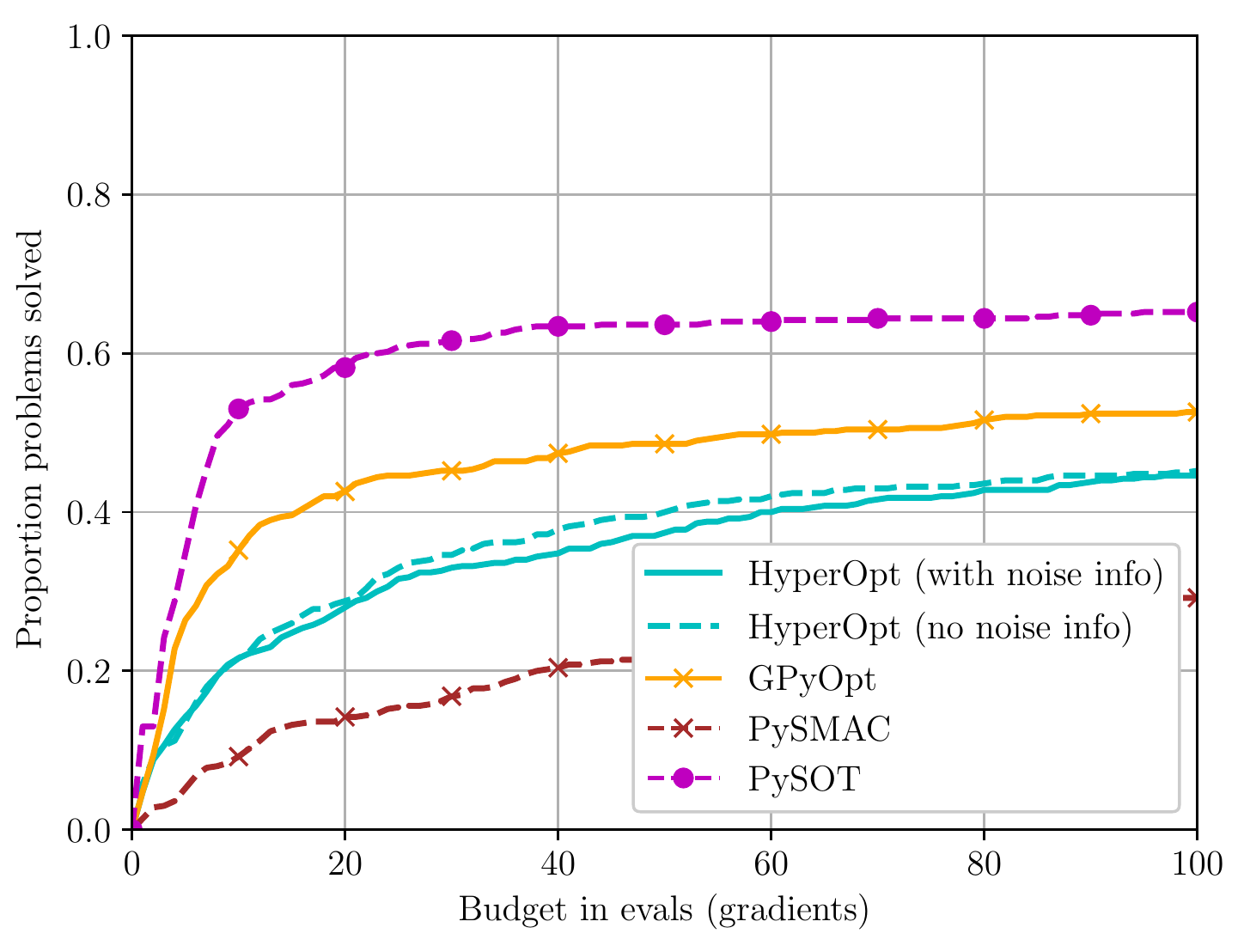}
		\caption{Low accuracy ($\tau=10^{-2}$), smooth.}
	\end{subfigure}
	~
	\begin{subfigure}[b]{0.48\textwidth}
		\includegraphics[width=\textwidth]{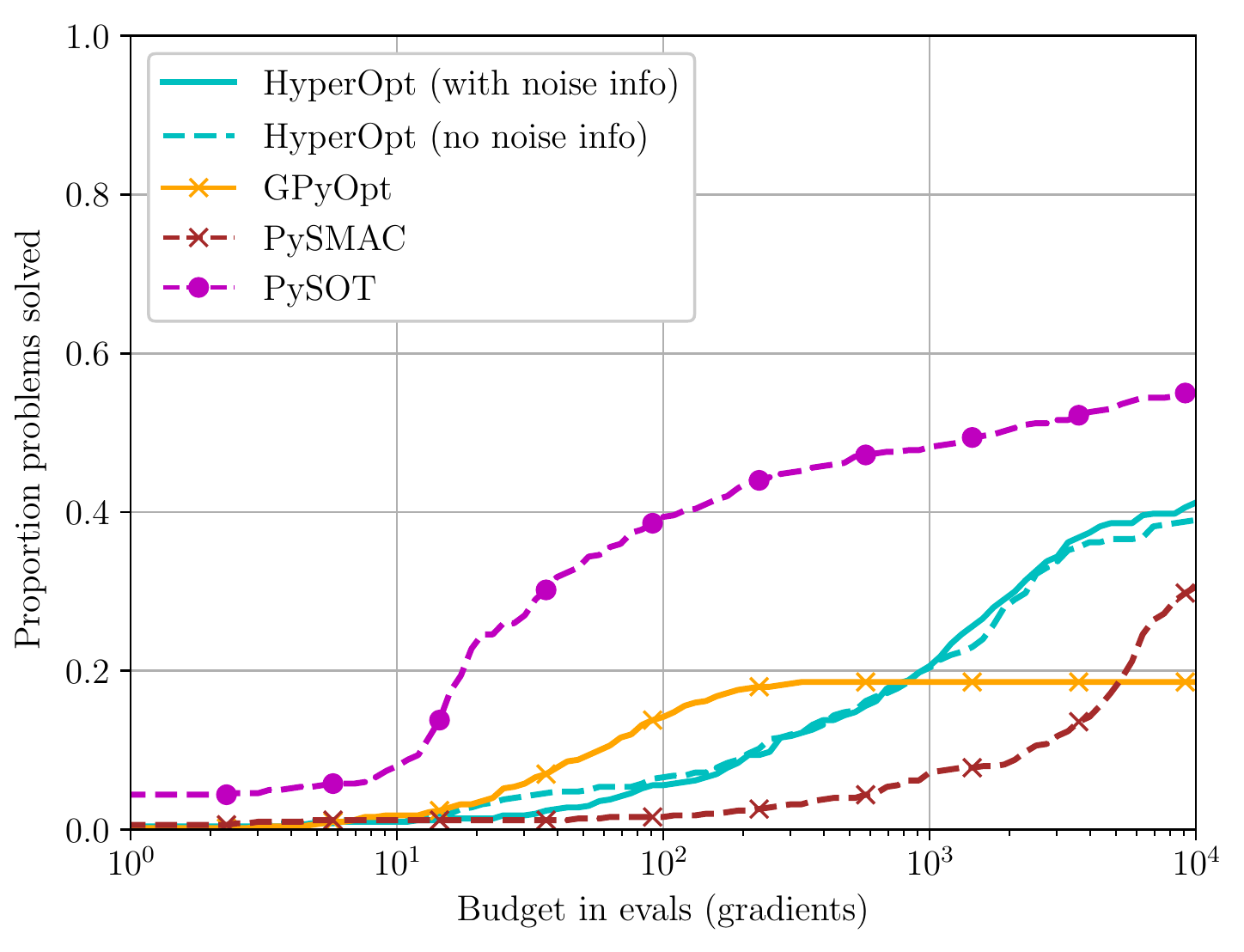}
		\caption{High accuracy ($\tau=10^{-5}$), smooth.}
	\end{subfigure}
	\caption{Comparison of Bayesian and surrogate solvers on the GO test set (smooth case; see \figref{BayesianComparison-GO-noisy} in \appref{app_pybobyqa_extra_results} for noisy results).}
	\label{BayesianComparison-GO}
\end{figure}

\figref{BayesianComparison-GO} shows that, for smooth problems, PySOT is consistently the best performing solver compared to all Bayesian solvers that we tested,
in both low and high accuracy regimes; the next best solver is GPyOpt for low accuracy regime, but when large budgets are available, HyperOpt and PySMAC improve their performance (compared to GPyOpt). 
The same results apply for noisy problems; these results are shown in \figref{BayesianComparison-GO-noisy} (\appref{app_pybobyqa_extra_results}).
Thus in our next set of comparisons between all solvers, we select only PySOT to `represent' the best behaviour of Bayesian and surrogate solvers.


\begin{figure}[t]
	\centering
	\begin{subfigure}[b]{0.48\textwidth}
		\includegraphics[width=\textwidth]{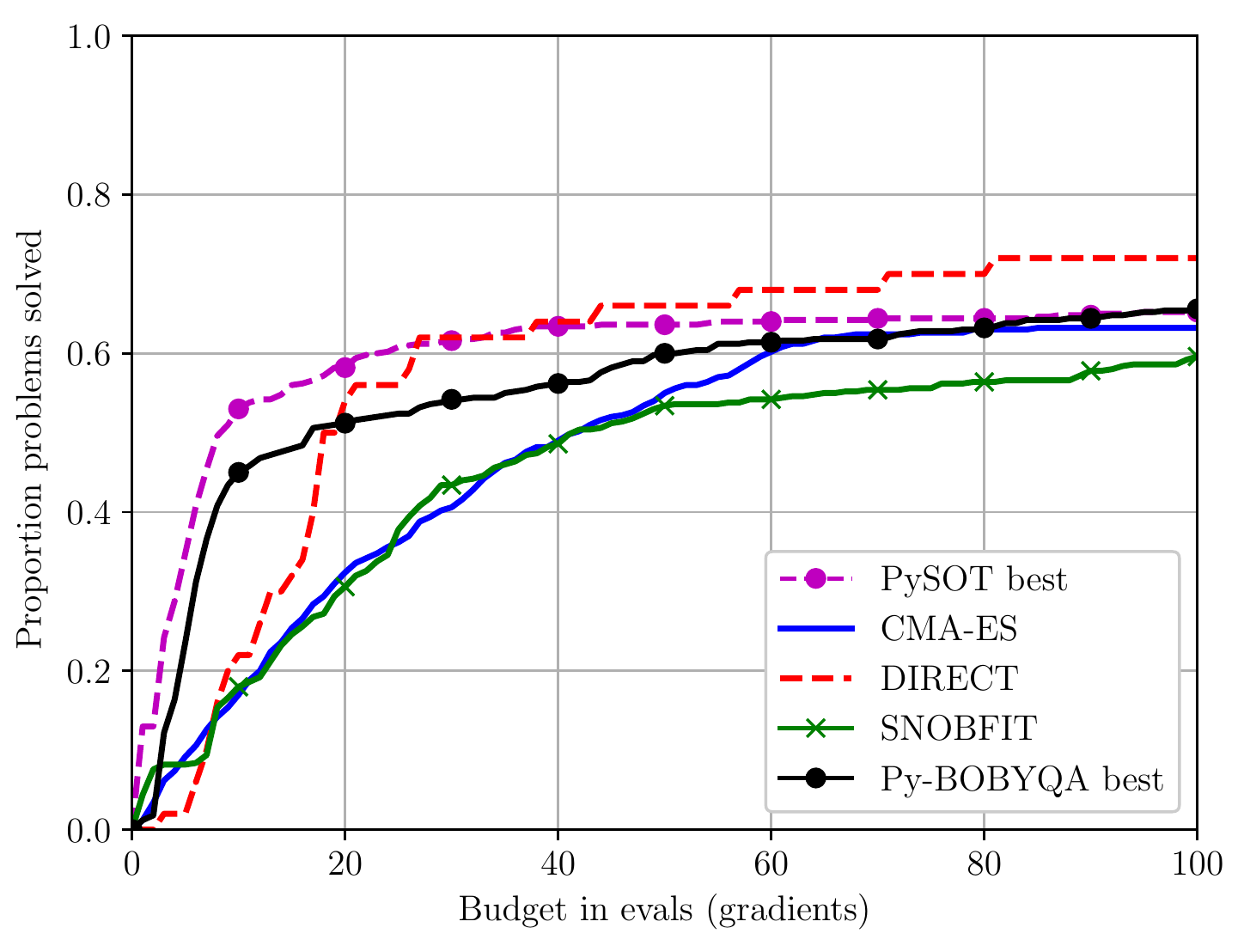}
		\caption{Low accuracy ($\tau=10^{-2}$), smooth.}
	\end{subfigure}
	~
	\begin{subfigure}[b]{0.48\textwidth}
		\includegraphics[width=\textwidth]{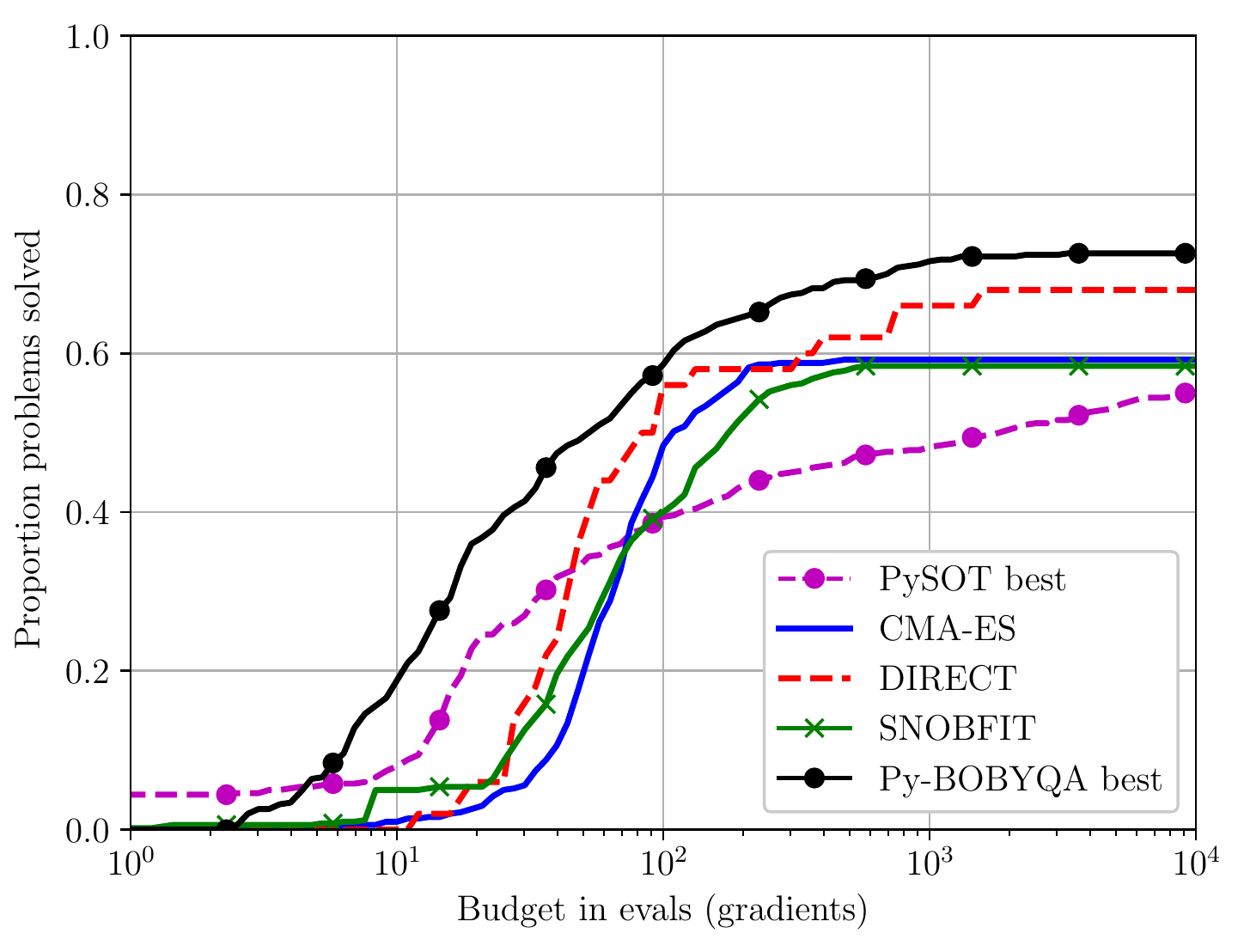}
		\caption{High accuracy ($\tau=10^{-5}$), smooth.}
	\end{subfigure}
	\\
	\begin{subfigure}[b]{0.48\textwidth}
		\includegraphics[width=\textwidth]{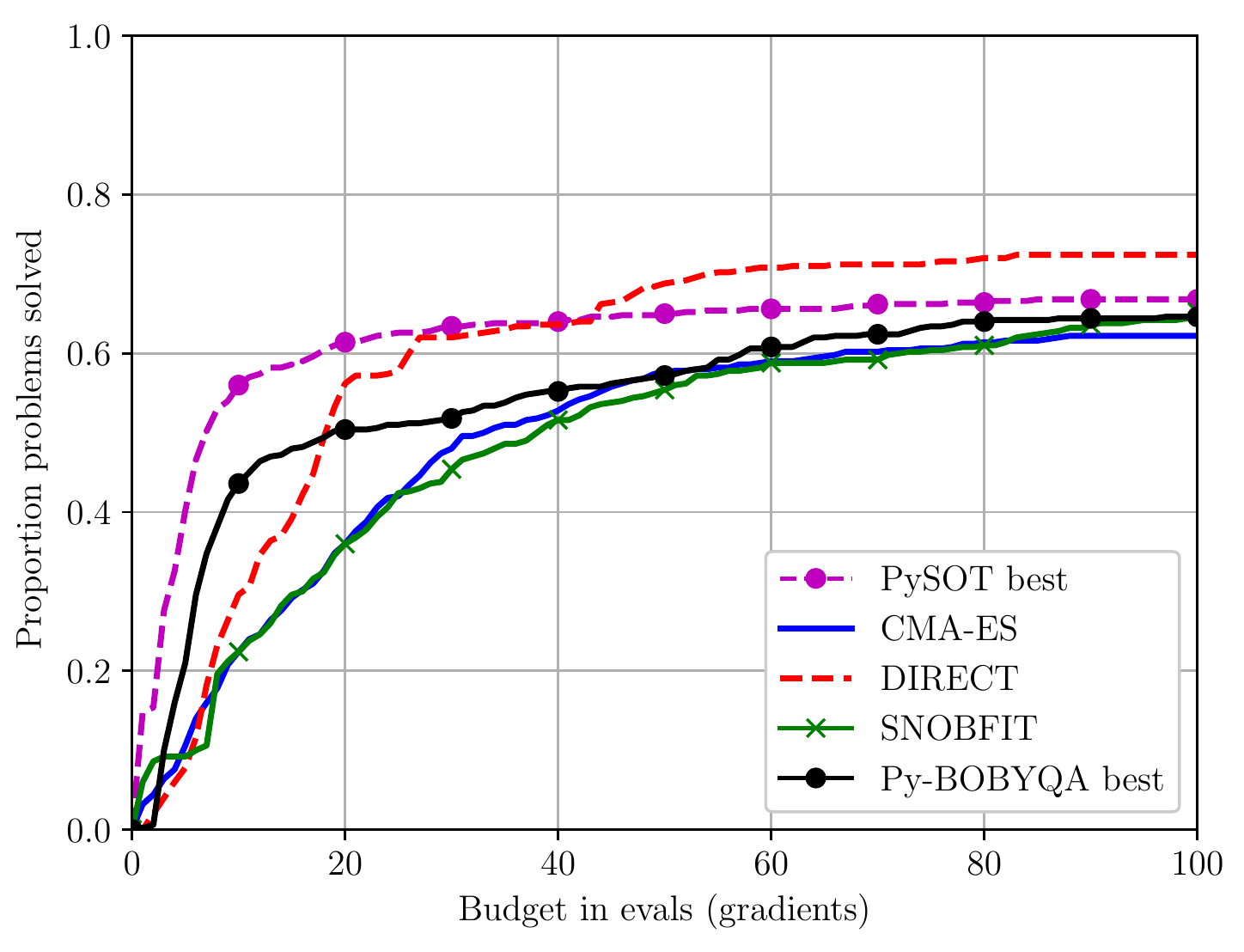}
		\caption{Low accuracy ($\tau=10^{-2}$), multiplicative noise.}
	\end{subfigure}
	~
	\begin{subfigure}[b]{0.48\textwidth}
		\includegraphics[width=\textwidth]{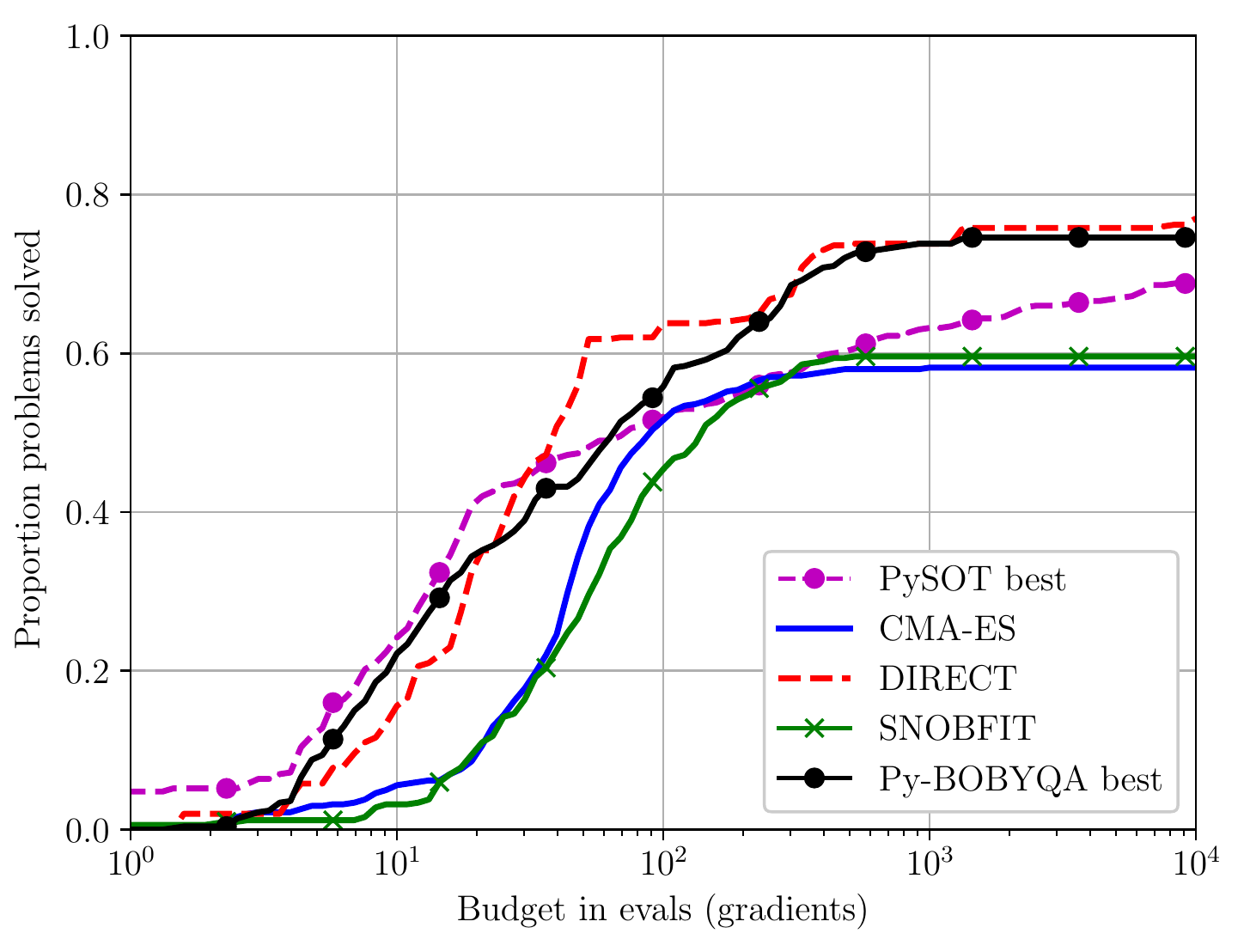}
		\caption{High accuracy ($\tau=10^{-5}$), multiplicative noise.}
	\end{subfigure}
	\\
	\begin{subfigure}[b]{0.48\textwidth}
		\includegraphics[width=\textwidth]{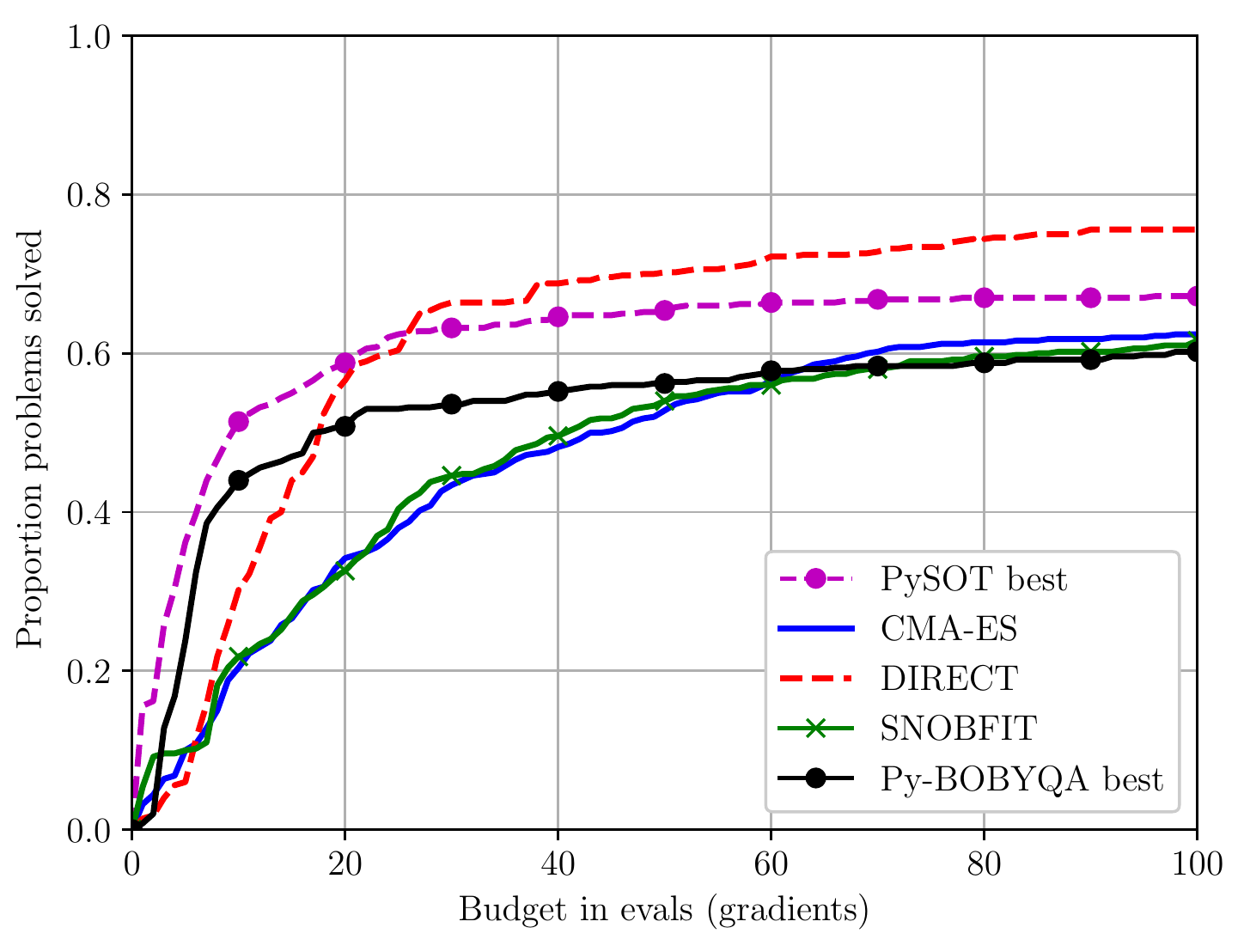}
		\caption{Low accuracy ($\tau=10^{-2}$), additive noise.}
	\end{subfigure}
	~
	\begin{subfigure}[b]{0.48\textwidth}
		\includegraphics[width=\textwidth]{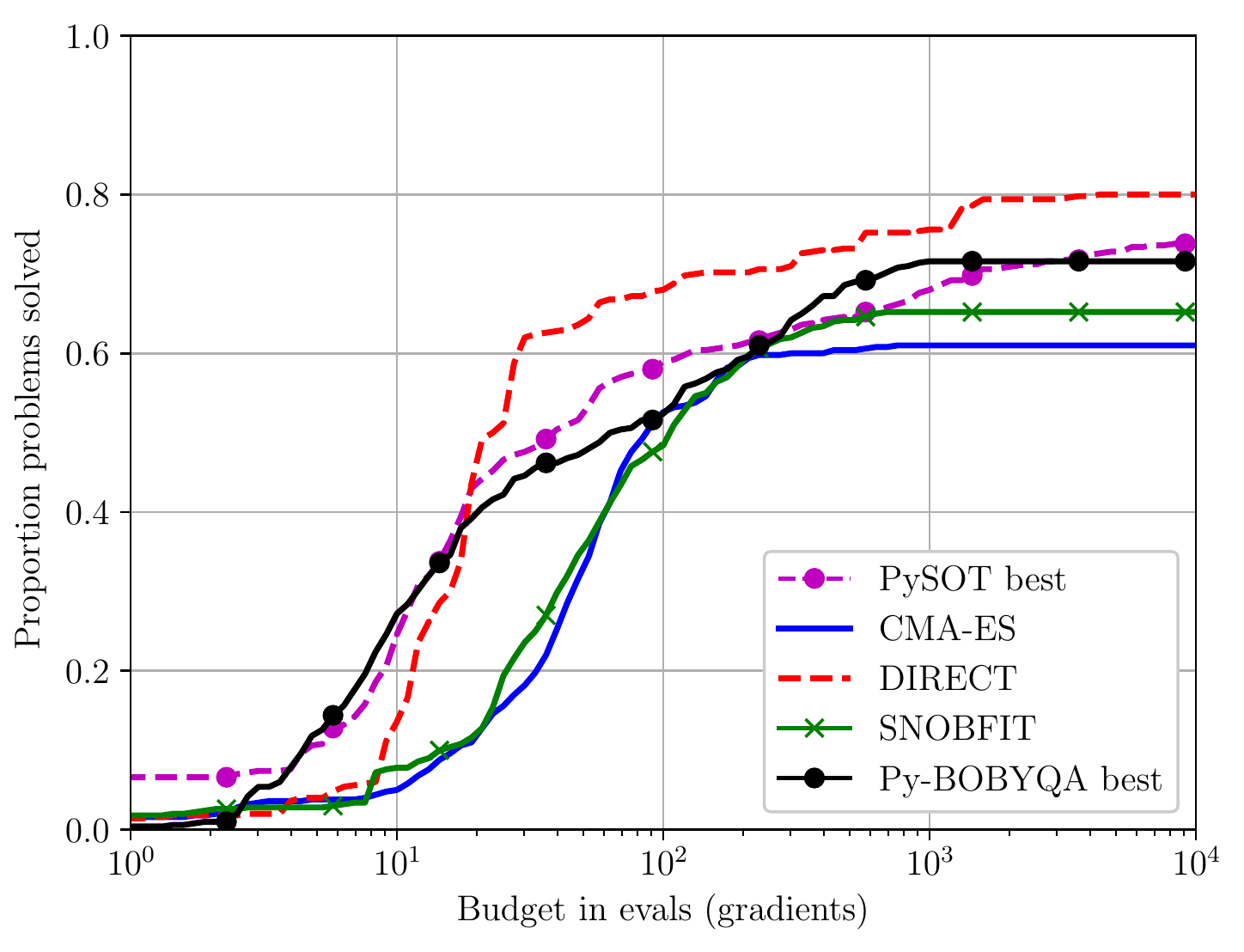}
		\caption{High accuracy ($\tau=10^{-5}$), additive noise.}
	\end{subfigure}
	\caption{Comparison of all and best-performing solvers on the GO test set.}
	\label{AllSolversComparison-GO}
\end{figure}

\paragraph{Comparisons of all solvers}


In  \figref{AllSolversComparison-GO}, we compare Py-BOBYQA against the global solvers; ``Py-BOBYQA best'' represents the variant with full quadratic interpolation model ($p=\mathcal{O}(n^2)$), soft restarts and adaptive trust region radius on restarts, that we found earlier
 to be performing best on the GO test set amongst all Py-BOBYQA variants, in all regimes of accuracy and noise; see Figures \ref{pybobyqa-smooth-high}--\ref{pybobyqa-add-high}, \figref{pybobyqa-low}, and our corresponding comments above. Also, ``PySOT best'' represents the best performing solver on the GO test set amongst the chosen Bayesian and surrogate solvers. 

{\it Low accuracy/budget} ($\tau=10^{-2}$ in \eqref{eq_reduction_measure} and budget $10^2 (n+1)$ evaluations)
The first column in \figref{AllSolversComparison-GO} shows that on average, DIRECT performs best in this regime, and the remaining solvers are essentially similar at the end of the budget. In the earlier phases of budget, Py-SOT and then the best 
Py-BOBYQA variant  have best performance,
being able to solve more problems on small budgets than the remaining solvers, while CMA-ES and SNOBFIT seem to need longer budgets to reach a similar amount of solved problems at the end of the low accuracy budget.

{\it High accuracy/budget} ($\tau=10^{-5}$ in \eqref{eq_reduction_measure} and budget $10^4 (n+1)$ evaluations)
For smooth problems, the best Py-BOBYQA variant is the best performing solver followed by DIRECT, while for multiplicative noise, these two solvers have comparably good performance. For additive noise, DIRECT has superior performance, followed by similarly good performance from PySOT and Py-BOBYQA. 




\subsection{Numerical results on MNIST}

\begin{figure}[t]
	\centering
	\begin{subfigure}[b]{0.48\textwidth}
		\includegraphics[width=\textwidth]{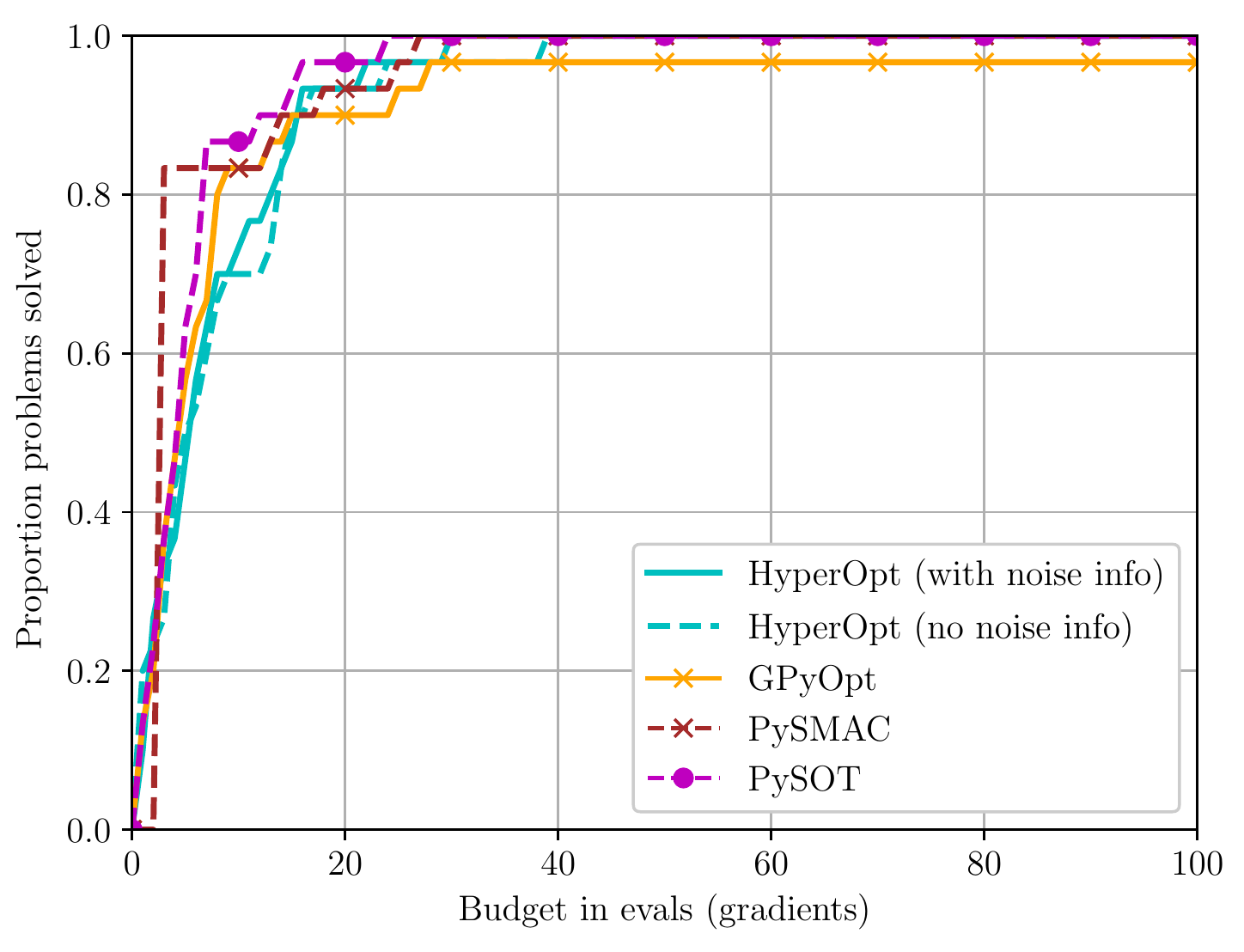}
		\caption{Bayesian and surrogate solvers.}
		\label{fig-mnist-bayesian}
	\end{subfigure}
	~
	\begin{subfigure}[b]{0.48\textwidth}
		\includegraphics[width=\textwidth]{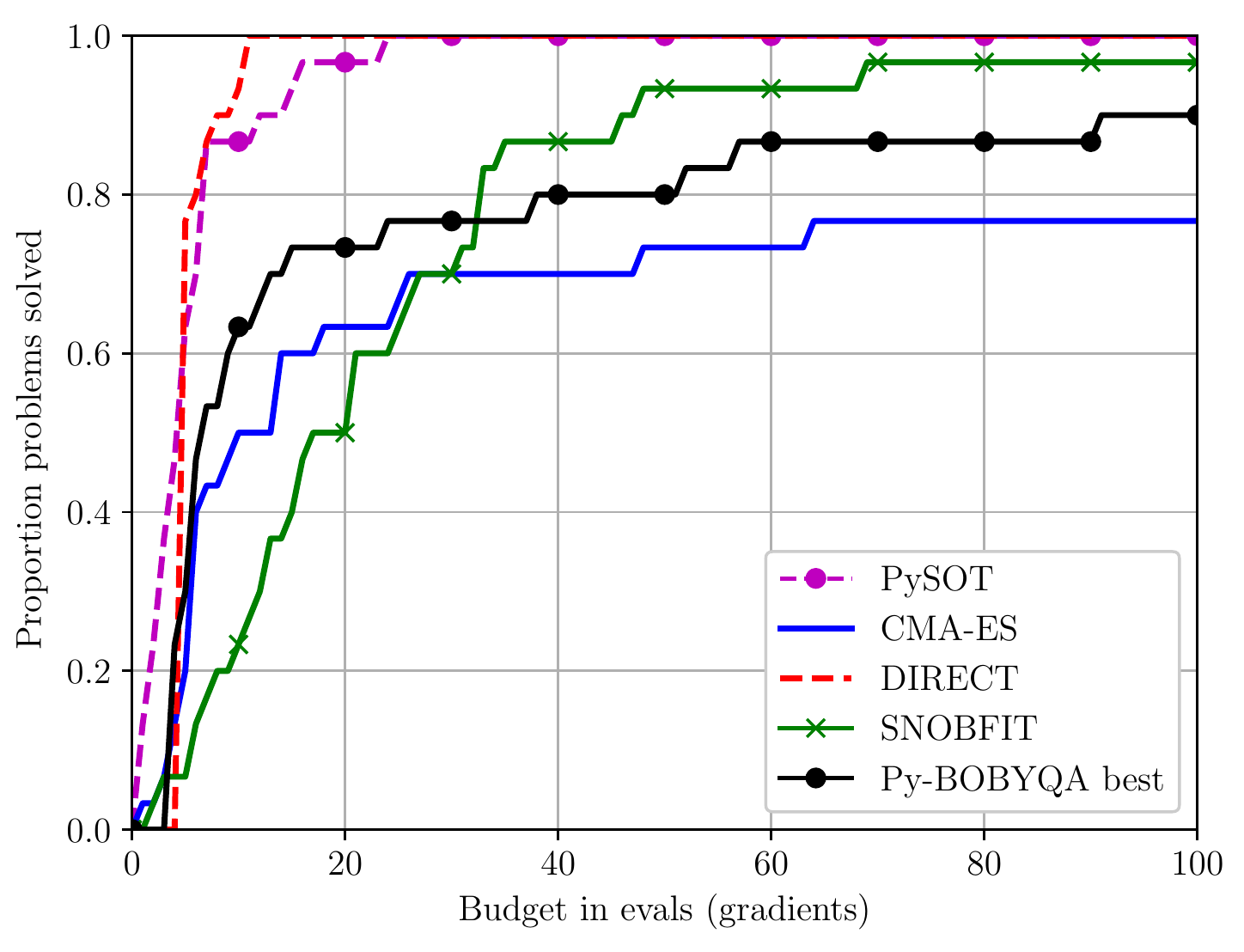}
		\caption{All/best performing solvers.}
		\label{fig-mnist-all}
	\end{subfigure}
	\caption{Comparison of solvers on the MNIST problem, low accuracy ($\tau=10^{-2}$).}
	\label{fig-mnist}
\end{figure}

We recall the description of the hyper-parameter tuning problem for the MNIST test set in Section \ref{sec_test_methodology},
where we described the four hyper-parameters, the objective value and the low accuracy level that is needed here.
Again, we compare first the Bayesian and surrogate optimization solvers in \figref{fig-mnist-bayesian}; we find similar performance
to the low accuracy regime for the GO test set (\figref{BayesianComparison-GO}), namely, the solvers are comparable with PySOT being slightly better. When comparing the remaining solvers with the best PyBOBYQA and PySOT in \figref{fig-mnist-all}, we find that DIRECT and PySOT are best performing, and the best performing Py-BOBYQA variant is not much worse. 

\section{Conclusions and future directions}
\label{sec_conclusions}

We investigated numerically the behaviour of a local derivative-free model-based trust region method Py-BOBYQA, and a new variant,
on global optimization problems. We found encouraging performance  of the improved variant when compared to Bayesian and surrogate algorithms, with DIRECT, SNOBFIT, and CMA-ES solvers, for both smooth and noisy problems solved to low or high accuracy requirements. Our experiments illustrate that local derivative free methods may be effective for global optimization, in the regime typical for global solvers, namely, when the objective is black-box, possibly noisy and computationally expensive to evaluate. This finding is useful as these local methods have greater   flexibility (such as allowing the use of a good initial guess of the solution), guarantee improvement in current best guess of the solution,  and have better scalability. 
However, we delegate to future directions the crucial task of estimating the performance of the local DFO solver on large-scale problems, both in terms of computational time  and evaluations. Theoretical developments that would show that the improved Py-BOBYQA can escape local minima would also be desirable.

\clearpage

\subsection{Acknowledgements}
We acknowledge the use of the University of Oxford Advanced Research Computing (ARC) facility\footnote{\url{http://dx.doi.org/10.5281/zenodo.22558}} in carrying out this work.








\addcontentsline{toc}{section}{References} 
\bibliographystyle{siam}
\enlargethispage{2em}
\bibliography{global_refs} 

\begin{thebibliography}{10}

\bibitem{Ali2005}
{\sc M.~M. Ali, C.~Khompatraporn, and Z.~B. Zabinsky}, {\em {A numerical
  evaluation of several stochastic algorithms on selected continuous global
  optimization test problems}}, J. Glob. Optim., 31 (2005), pp.~635--672.

\bibitem{Ali1994}
{\sc M.~M. Ali and C.~Storey}, {\em {Topographical Multilevel Single Linkage}},
  Journal of Global Optimization, 5 (1994), pp.~349--358.

\bibitem{Audet2017}
{\sc C.~Audet and W.~Hare}, {\em Derivative-Free and Blackbox Optimization},
  Springer Series in Operations Research and Financial Engineering, Springer,
  Cham, Switzerland, 2017.

\bibitem{pysmac_URL}
{\sc {AutoML}}, {\em {pySMAC}}.
\newblock \url{https://github.com/automl/pysmac}, 2018.

\bibitem{Bergstra2013}
{\sc J.~Bergstra, D.~L.~K. Yamins, and D.~D. Cox}, {\em {Making a Science of
  Model Search: Hyperparameter Optimization in Hundreds of Dimensions for
  Vision Architectures}}, in 30th Int. Conf. Mach. Learn., Atlanta, 2013,
  pp.~115--123.

\bibitem{HyperOpt_URL}
\leavevmode\vrule height 2pt depth -1.6pt width 23pt, {\em {Hyperopt:
  Distributed Asynchronous Hyper-parameter Optimization}}.
\newblock \url{https://github.com/hyperopt/hyperopt}, 2018.

\bibitem{brochu2010tutorial}
{\sc E.~Brochu, V.~M. Cora, and N.~De~Freitas}, {\em A tutorial on {B}ayesian
  optimization of expensive cost functions, with application to active user
  modeling and hierarchical reinforcement learning}, arXiv preprint
  arXiv:1012.2599,  (2010).

\bibitem{TR-DFOLS}
{\sc C.~Cartis, J.~Fiala, B.~Marteau, and L.~Roberts}, {\em Improving the
  flexibility and robustness of model-based derivative-free optimization
  solvers}, ACM Transactions on Mathematical Software, 45 (2019),
  pp.~32:1--32:41.

\bibitem{Cartis2017a}
{\sc C.~Cartis and L.~Roberts}, {\em A derivative-free {G}auss-{N}ewton
  method}, Mathematical Programming Computation,  (2019).

\bibitem{Conn2000}
{\sc A.~R. Conn, N.~I.~M. Gould, and P.~L. Toint}, {\em {Trust-Region
  Methods}}, MPS-SIAM Series on Optimization, MPS/SIAM, Philadelphia, 2000.

\bibitem{CSV}
{\sc A.~R. Conn, K.~Scheinberg, and L.~N. Vicente}, {\em {Introduction to
  Derivative-Free Optimization}}, vol.~8 of MPS-SIAM Series on Optimization,
  MPS/SIAM, Philadelphia, 2009.

\bibitem{Custodio2011}
{\sc A.~L. Cust{\'{o}}dio, J.~F.~A. Madeira, A.~I.~F. Vaz, and L.~N. Vicente},
  {\em Direct multisearch for multiobjective optimization}, SIAM Journal on
  Optimization, 21 (2011), pp.~1109--1140.

\bibitem{Custodio2017}
{\sc A.~L. Cust{\'{o}}dio, K.~Scheinberg, and L.~N. Vicente}, {\em
  {Methodologies and Software for Derivative-free Optimization}}, in Adv.
  Trends Optim. with Eng. Appl., T.~Terlaky, M.~F. Anjos, and S.~Ahmed, eds.,
  MOS-SIAM Book Series on Optimization, SIAM, Philadelphia, 2017.

\bibitem{pysot2015}
{\sc D.~Eriksson, D.~Bindel, and C.~Shoemaker}, {\em Surrogate optimization
  toolbox ({pySOT})}.
\newblock \url{https://github.com/dme65/pySOT}, 2015.

\bibitem{fowkes2011bayesian}
{\sc J.~Fowkes}, {\em Bayesian Numerical Analysis: Global Optimization and
  Other Applications}, PhD thesis, Oxford University, 2011.

\bibitem{snobfit_URL}
{\sc Z.~Fu}, {\em {Package snobfit}}.
\newblock \url{http://reflectometry.org/danse/docs/snobfit}, 2009.

\bibitem{Ghanbari}
{\sc H.~Ghanbari and K.~Scheinberg}, {\em Black-box optimization in machine
  learning with trust region based derivative free algorithm}, arXiv preprint
  arXiv:1703.06925,  (2017).

\bibitem{gutmann}
{\sc H.-M. Gutmann}, {\em Radial Basis Function Methods for Global
  Optimization}, PhD thesis, Cambridge University, 2001.

\bibitem{cmaes_URL}
{\sc N.~Hansen}, {\em {pycma}}.
\newblock \url{https://github.com/CMA-ES/pycma}, 2018.

\bibitem{Hansen1996}
{\sc N.~Hansen and A.~Ostermeier}, {\em {Adapting arbitrary normal mutation
  distributions in evolution strategies: the covariance matrix adaptation}}, in
  Proc. IEEE Int. Conf. Evol. Comput., Nagoya, 1996, pp.~312--317.

\bibitem{HansenOstermeier2001}
\leavevmode\vrule height 2pt depth -1.6pt width 23pt, {\em Completely
  derandomized self-adaptation in evolution strategies}, Evolutionary
  Computation, 9 (2001), pp.~159--195.

\bibitem{Hare2018}
{\sc W.~Hare, J.~Loeppky, and S.~Xie}, {\em Methods to compare expensive
  stochastic optimization algorithms with random restarts}, Journal of Global
  Optimization, 72 (2018), pp.~781--801.

\bibitem{Hutter2011}
{\sc F.~Hutter, H.~H. Hoos, and K.~Leyton-Brown}, {\em {Sequential Model-Based
  Optimization for General Algorithm Configuration}}, in Int. Conf. Learn.
  Intell. Optim., Rome, 2011, pp.~507--223.

\bibitem{hutter2007automatic}
{\sc F.~Hutter, H.~H. Hoos, and T.~St{\"u}tzle}, {\em Automatic algorithm
  configuration based on local search}, in AAAI, vol.~7, 2007, pp.~1152--1157.

\bibitem{Huyer2008}
{\sc W.~Huyer and A.~Neumaier}, {\em {SNOBFIT -- Stable Noisy Optimization by
  Branch and Fit}}, ACM Trans. Math. Softw., 35 (2008), pp.~1--25.

\bibitem{pydirect_URL}
{\sc {J. M. Gablonsky}}, {\em {pydirect}}.
\newblock \url{https://code.google.com/archive/p/pydirect}, 2012.

\bibitem{jones2001taxonomy}
{\sc D.~R. Jones}, {\em A taxonomy of global optimization methods based on
  response surfaces}, Journal of Global Optimization, 21 (2001), pp.~345--383.

\bibitem{Jones1993}
{\sc D.~R. Jones, C.~D. Perttunen, and B.~E. Stuckman}, {\em {Lipschitzian
  optimization without the Lipschitz constant}}, J. Optim. Theory Appl., 79
  (1993), pp.~157--181.

\bibitem{jones1998efficient}
{\sc D.~R. Jones, M.~Schonlau, and W.~J. Welch}, {\em Efficient global
  optimization of expensive black-box functions}, Journal of Global
  Optimization, 13 (1998), pp.~455--492.

\bibitem{Kelley2011}
{\sc C.~T. Kelley}, {\em Implicit Filtering}, Software, Environments and Tools,
  SIAM, Philadelphia, 2011.

\bibitem{Kolda2003}
{\sc T.~Kolda, R.~Lewis, and V.~Torczon}, {\em Optimization by direct search:
  new perspectives on some classical and modern methods}, SIAM Review, 45
  (2003), pp.~385--482.

\bibitem{li2016hyperband}
{\sc L.~Li, K.~Jamieson, G.~DeSalvo, A.~Rostamizadeh, and A.~Talwalkar}, {\em
  Hyperband: A novel bandit-based approach to hyperparameter optimization},
  arXiv preprint arXiv:1603.06560,  (2016).

\bibitem{Locatelli2003}
{\sc M.~Locatelli}, {\em A note on the {G}riewank test function}, Journal of
  Global Optimization, 25 (2003), pp.~169--174.

\bibitem{Locatelli2013}
{\sc M.~Locatelli and F.~Schoen}, {\em Global Optimization: Theory, Algorithms,
  and Applications}, vol.~15 of MOS-SIAM Series on Optimization, MOS/SIAM,
  Philadelphia, 2013.

\bibitem{mockus1978toward}
{\sc J.~Mockus, V.~Tiesis, and A.~Zilinskas}, {\em Bayesian methods for seeking
  the extremum}, Toward Global Optimization, 2 (1978), pp.~117--129.

\bibitem{More2009}
{\sc J.~J. Mor{\'e} and S.~M. Wild}, {\em Benchmarking derivative-free
  optimization algorithms}, SIAM Journal on Optimization, 20 (2009),
  pp.~172--191.

\bibitem{Nocedal2006}
{\sc J.~Nocedal and S.~J. Wright}, {\em {Numerical Optimization}}, Springer
  Series in Operations Research and Financial Engineering, Springer, New York,
  2nd~ed., 2006.

\bibitem{Powell1998}
{\sc M.~J.~D. Powell}, {\em Direct search algorithms for optimization
  calculations}, Acta Numerica, 7 (1998), pp.~287--336.

\bibitem{Powell2007}
\leavevmode\vrule height 2pt depth -1.6pt width 23pt, {\em {Developments of
  NEWUOA for minimization without derivatives}}, Tech. Rep. DAMTP 2007/NA05,
  University of Cambridge, 2007.

\bibitem{Powell2009}
\leavevmode\vrule height 2pt depth -1.6pt width 23pt, {\em {The BOBYQA
  algorithm for bound constrained optimization without derivatives}}, Tech.
  Rep. DAMTP 2009/NA06, University of Cambridge, 2009.

\bibitem{rasmussen2006gaussian}
{\sc C.~E. Rasmussen and C.~K. Williams}, {\em Gaussian processes for machine
  learning. 2006}, The MIT Press, Cambridge, MA, USA, 38 (2006), pp.~715--719.

\bibitem{Regis2013}
{\sc R.~G. Regis and C.~A. Shoemaker}, {\em {Combining radial basis function
  surrogates and dynamic coordinate search in high-dimensional expensive
  black-box optimization}}, Eng. Optim., 45 (2013), pp.~529--555.

\bibitem{RinnooyKan1987}
{\sc A.~H.~G. {Rinnooy Kan} and G.~T. Timmer}, {\em Stochastic global
  optimization methods part {I}: Clustering methods}, Mathematical Programming,
  39 (1987), pp.~27--56.

\bibitem{RinnooyKan1987a}
\leavevmode\vrule height 2pt depth -1.6pt width 23pt, {\em Stochastic global
  optimization methods part {II}: Multi level methods}, Mathematical
  Programming, 39 (1987), pp.~57--78.

\bibitem{Rios2013}
{\sc L.~M. Rios and N.~V. Sahinidis}, {\em Derivative-free optimization: a
  review of algorithms and comparison of software implementations}, Journal of
  Global Optimization, 56 (2013), pp.~1247 -- 1293.

\bibitem{TR-Miniproject-Lindon}
{\sc L.~Roberts}, {\em Derivative-free optimisation for data fitting}, Tech.
  Rep. InFoMM CDT Report, University of Oxford, 2016.
\newblock
  \url{http://people.maths.ox.ac.uk/robertsl/docs/DFO_MiniprojectReport_updateNov18.pdf}.

\bibitem{shahriari2016taking}
{\sc B.~Shahriari, K.~Swersky, Z.~Wang, R.~P. Adams, and N.~de~Freitas}, {\em
  Taking the human out of the loop: A review of {B}ayesian optimization},
  Proceedings of the IEEE, 104 (2016), pp.~148--175.

\bibitem{Snoek2012}
{\sc J.~Snoek, H.~Larochelle, and R.~P. Adams}, {\em {Practical Bayesian
  Optimization of Machine Learning Algorithms.}}, in Adv. Neural Inf. Process.
  Syst., 2012, pp.~2951--2959.

\bibitem{gpyopt2016}
{\sc {The GPyOpt authors}}, {\em {GPyOpt}: A {Bayesian} optimization framework
  in {Python}}.
\newblock \url{http://github.com/SheffieldML/GPyOpt}, 2016.

\bibitem{Torn1994}
{\sc A.~T{\"{o}}rn and S.~Viitanen}, {\em Topographical global optimization
  using pre-sampled points}, Journal of Global Optimization, 5 (1994),
  pp.~267--276.

\bibitem{MHWright}
{\sc M.~H. Wright}, {\em Direct search methods: once scorned, now respectable},
  in Numerical Analysis 1995 (Proceedings of the 1995 Dundee Biennal Conference
  in Numerical Analysis), D.~Griffiths and G.~Watson, eds., 1996, pp.~191--208.

\end{thebibliography}
\appendix

\section{Test Problems} \label{app_test_problems}
\begin{table}[H]
	\centering
	\small{
	\begin{tabular}{rlccc} 
		\hline\noalign{\smallskip}
		\textbf{\#} & \textbf{Problem} & \textbf{Dimension} & \textbf{Bounds} & \textbf{Global minimum} \\ \noalign{\smallskip}\hline\noalign{\smallskip} 
		1 & Ackley & 10 & $-30 \leq x_i \leq 30$ & $0$ \\
		2 & Aluffi-Pentini & 2 & $-10 \leq x_i \leq 10$ & $-0.3523861$ \\
		3 & Becker and Lago & 2 & $-10 \leq x_i \leq 10$ & $0$ \\
		4 & Bohachevsky 1 & 2 & $-50 \leq x_i \leq 50$ & $0$ \\
		5 & Bohachevsky 2 & 2 & $-50 \leq x_i \leq 50$ & $0$ \\
		6 & Branin & 2 & $[-5, 0] \leq \bx \leq [10, 15]$ & $3.926991$ \\
		7 & Camel 3 & 2 & $-5 \leq x_i \leq 5$ & $0$ \\
		8 & Camel 6 & 2 & $-5 \leq x_i \leq 5$ & $-1.031628$ \\
		9 & Cosine Mixture & 4 & $-1 \leq x_i \leq 1$ & $-0.4$ \\
		10 & Dekkers and Aarts & 2 & $-20 \leq x_i \leq 20$ & $-2.477289\times 10^{4}$ \\
		11 & Easom & 2 & $-10 \leq x_i \leq 10$ & $-1$ \\
		12 & Epistatic Michalewicz & 5 & $0 \leq x_i \leq \pi$ & $-4.687658$ \\
		13 & Exponential* & 40 & $-1 \leq x_i \leq 1$ & $-1$ \\
		14 & Goldstein and Price & 2 & $-2 \leq x_i \leq 2$ & $3$ \\
		15 & Griewank* & 5 & $-600 \leq x_i \leq 600$ & $0$ \\
		16 & Gulf Research & 3 & $[0.1, 0, 0] \leq \bx \leq [100, 25.6, 5]$ & $0$ \\
		17 & Hartman 3 & 3 & $0 \leq x_i \leq 1$ & $-3.862782$ \\
		18 & Hartman 6 & 6 & $0 \leq x_i \leq 1$ & $-3.322368$ \\
		19 & Helical Valley & 3 & $-10 \leq x_i \leq 10$ & $0$ \\
		20 & Hosaki & 2 & $[0, 0] \leq \bx \leq [5, 6]$ & $-2.345812$ \\
		21 & Kowalik & 4 & $0 \leq x_i \leq 0.42$ & $3.074871\times 10^{-4}$ \\
		22 & Levy and Montalvo 1 & 3 & $-10 \leq x_i \leq 10$ & $0$ \\
		23 & Levy and Montalvo 2 & 10 & $-5 \leq x_i \leq 5$ & $0$ \\
		24 & McCormick & 2 & $[-1.5, -3] \leq \bx \leq [4, 3]$ & $-1.913223$ \\
		25 & Meyer and Roth & 3 & $-20 \leq x_i \leq 20$ & $4.355269\times 10^{-5}$ \\
		26 & Miele and Cantrell & 4 & $-1 \leq x_i \leq 1$ & $0$ \\
		27 & Modified Langerman & 10 & $0 \leq x_i \leq 10$ & $-0.965$ \\
		28 & Modified Rosenbrock & 2 & $-5 \leq x_i \leq 5$ & $0$ \\
		29 & Multi-Gaussian & 2 & $-2 \leq x_i \leq 2$ & $-1.296954$ \\
		30 & Neumaier 2 & 4 & $0 \leq x_i \leq 4$ & $0$ \\
		31 & Neumaier 3* & 30 & $-900 \leq x_i \leq 900$ & $-4930$ \\
		32 & Odd Square & 20 & $-15 \leq x_i \leq 15$ & $-1$ \\
		33 & Paviani & 10 & $2.001 \leq x_i \leq 9.999$ & $-45.77845$ \\
		34 & Periodic & 2 & $-10 \leq x_i \leq 10$ & $0.9$ \\
		35 & Powell Quadratic & 4 & $-10 \leq x_i \leq 10$ & $0$ \\
		36 & Price Transistor Monitoring & 9 & $-10 \leq x_i \leq 10$ & $0$ \\
		37 & Rastrigin* & 30 & $-5.12 \leq x_i \leq 5.12$ & $0$ \\
		38 & Rosenbrock* & 50 & $-30 \leq x_i \leq 30$ & $0$ \\
		39 & Salomon* & 50 & $-100 \leq x_i \leq 100$ & $0$ \\
		40 & Schaffer 1 & 2 & $-100 \leq x_i \leq 100$ & $0$ \\
		41 & Schaffer 2 & 2 & $-100 \leq x_i \leq 100$ & $0$ \\
		42 & Schwefel* & 40 & $-500 \leq x_i \leq 500$ & $-1.675932\times 10^{4}$ \\
		43 & Shekel 10 & 4 & $0 \leq x_i \leq 10$ & $-10.53628$ \\
		44 & Shekel 5 & 4 & $0 \leq x_i \leq 10$ & $-10.15320$ \\
		45 & Shekel 7 & 4 & $0 \leq x_i \leq 10$ & $-10.40282$ \\
		46 & Shekel Foxholes & 10 & $0 \leq x_i \leq 10$ & $-10.20879$ \\
		47 & Shubert & 2 & $-10 \leq x_i \leq 10$ & $-186.7309$ \\
		48 & Sinusoidal & 20 & $0 \leq x_i \leq 180$ & $-3.5$ \\
		49 & Storn Tchebyshev & 9 & $-256 \leq x_i \leq 256$ & $0$ \\
		50 & Wood & 4 & $-10 \leq x_i \leq 10$ & $0$ \\
		\noalign{\smallskip}\hline
	\end{tabular}}
	\caption{Final set of global optimization test problems, taken from \cite{Ali2005}. Problems marked with * had $n=10$ originally, but had their dimension changed here.}
	\label{tab_final_global_problems}
\end{table}

\section{Further Numerical Results} \label{app_pybobyqa_extra_results}

\subsection{Selecting best Py-BOBYQA variant: low accuracy and budget}

\begin{figure}[H]
	\centering
	\begin{subfigure}[b]{0.48\textwidth}
		\includegraphics[width=\textwidth]{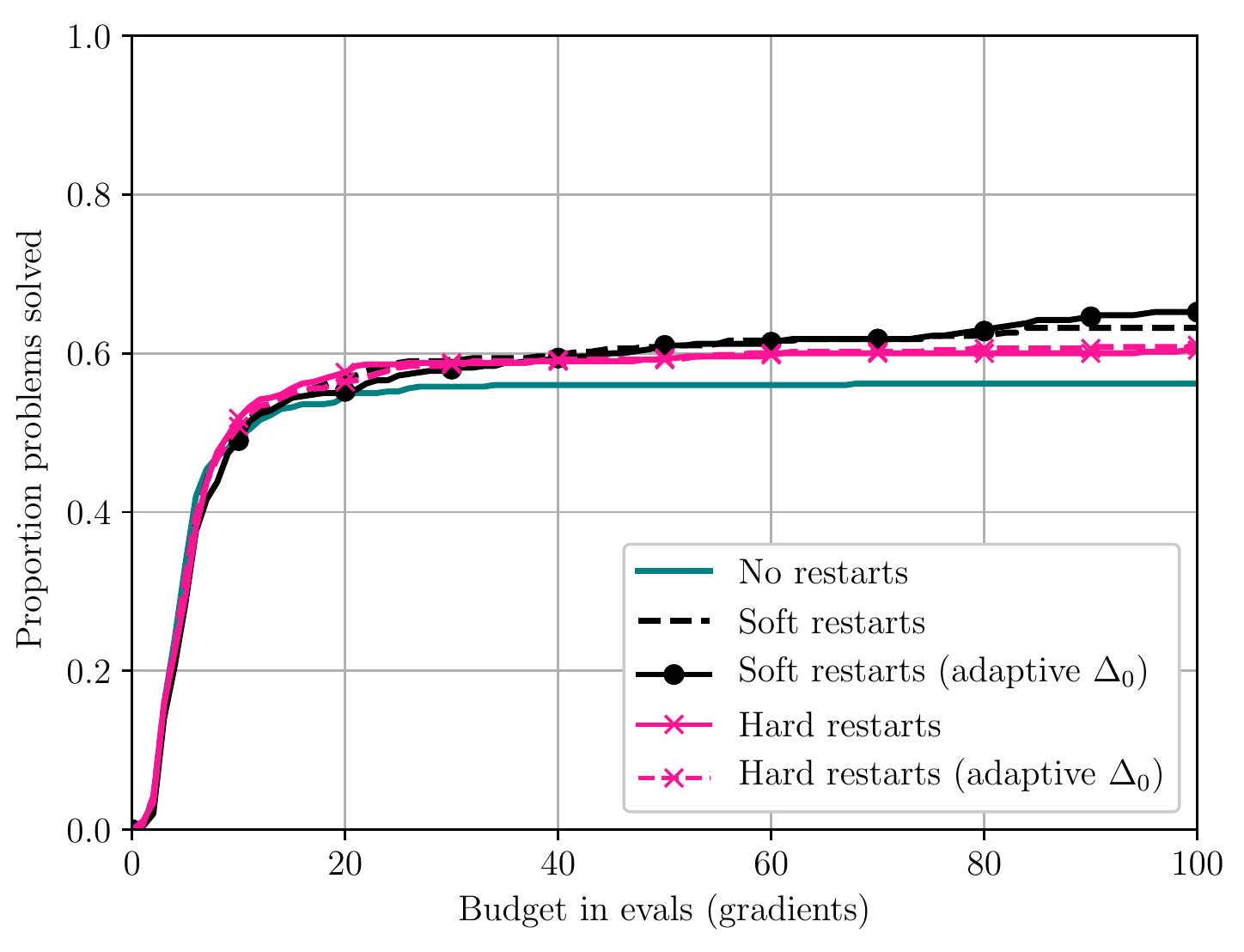}
		\caption{Comparison of Py-BOBYQA ($p=2n+1$) variants  on smooth GO test set}
	\end{subfigure}
	~
	\begin{subfigure}[b]{0.48\textwidth}
		\includegraphics[width=\textwidth]{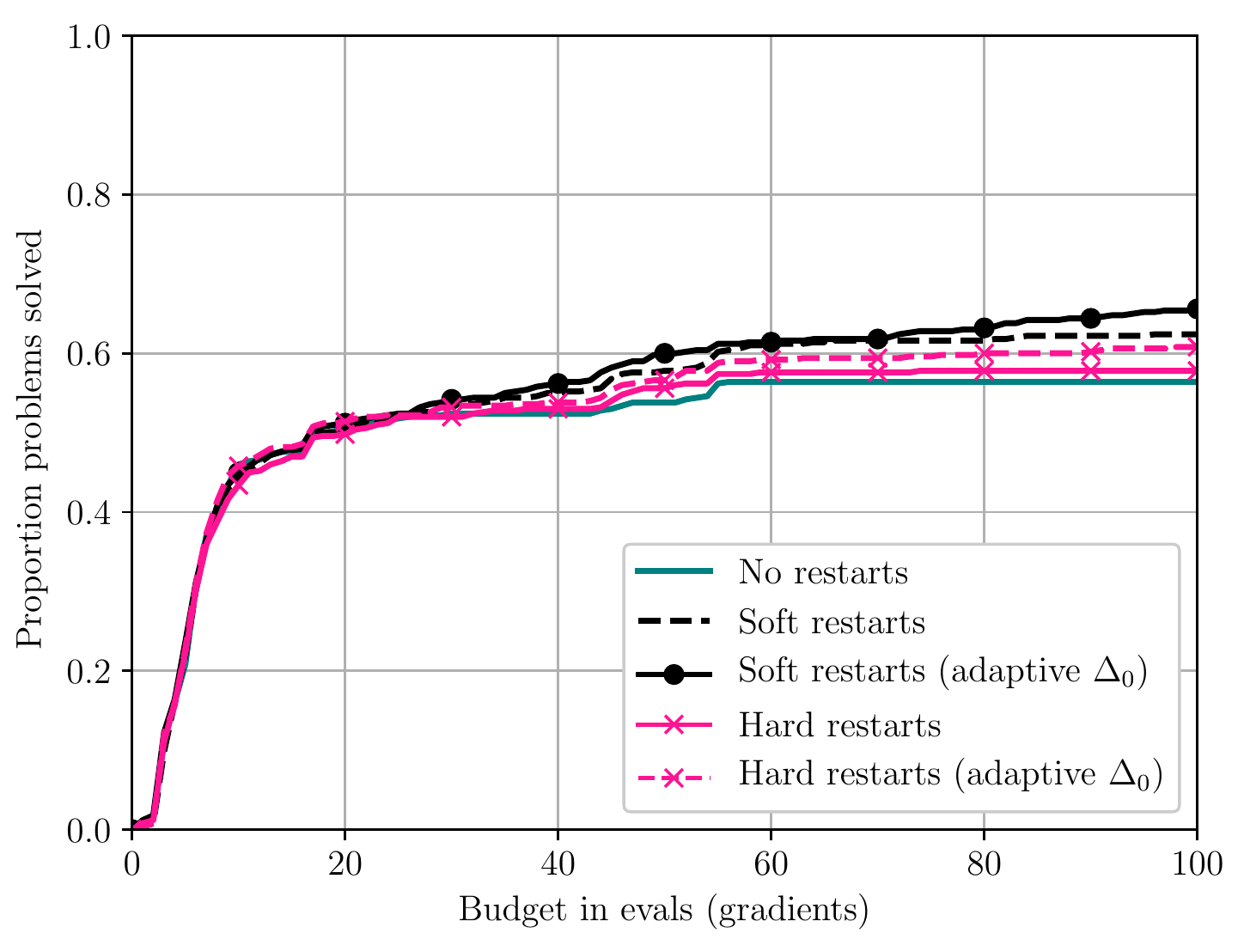}
		\caption{Comparison of Py-BOBYQA ($p=\mathcal{O}(n^2)$) variants  on smooth GO test set}
	\end{subfigure}
	\\
	\begin{subfigure}[b]{0.48\textwidth}
		\includegraphics[width=\textwidth]{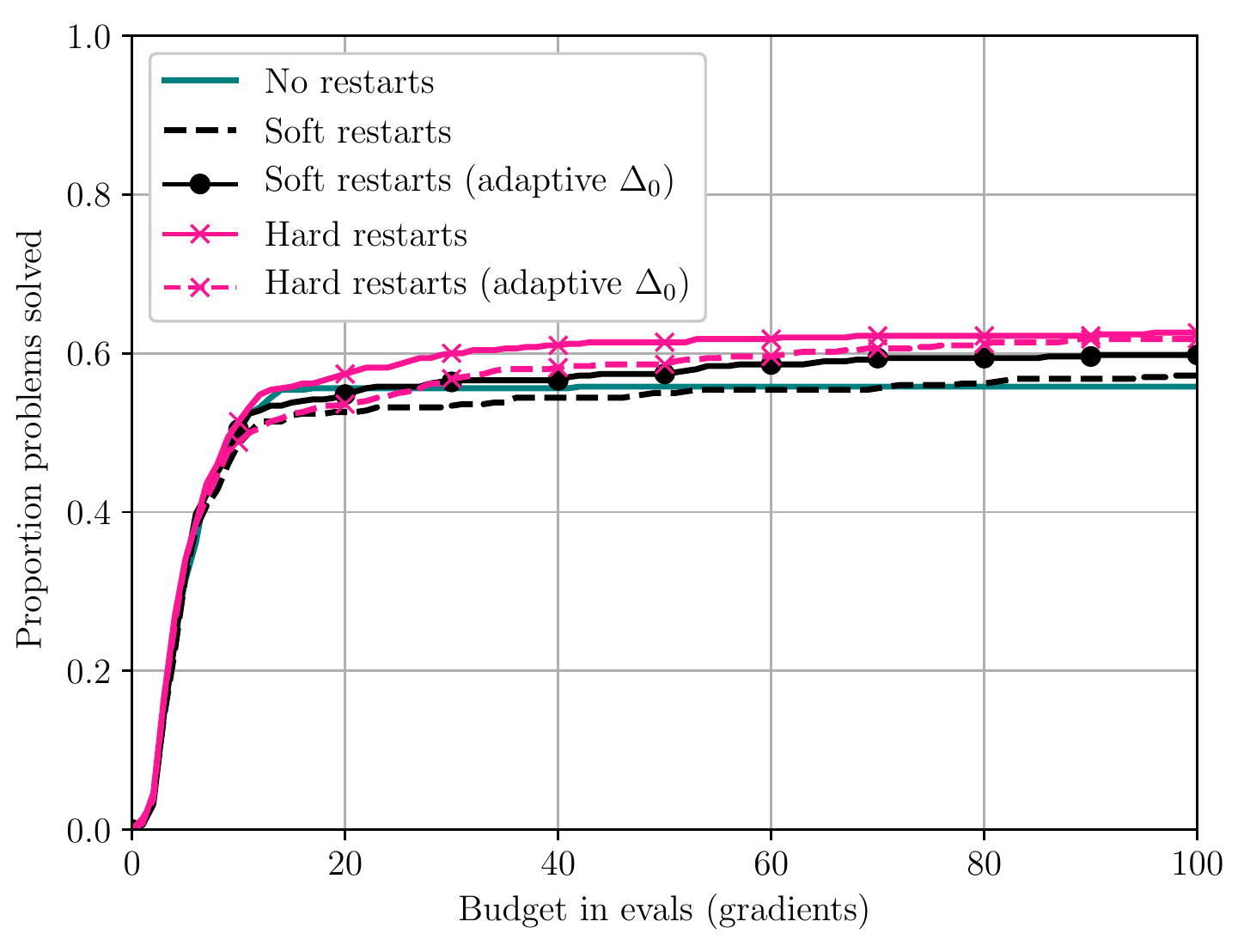}
		\caption{Comparison of Py-BOBYQA ($p=2n+1$) variants  on GO test set with multiplicative noise}
	\end{subfigure}
	~
	\begin{subfigure}[b]{0.48\textwidth}
		\includegraphics[width=\textwidth]{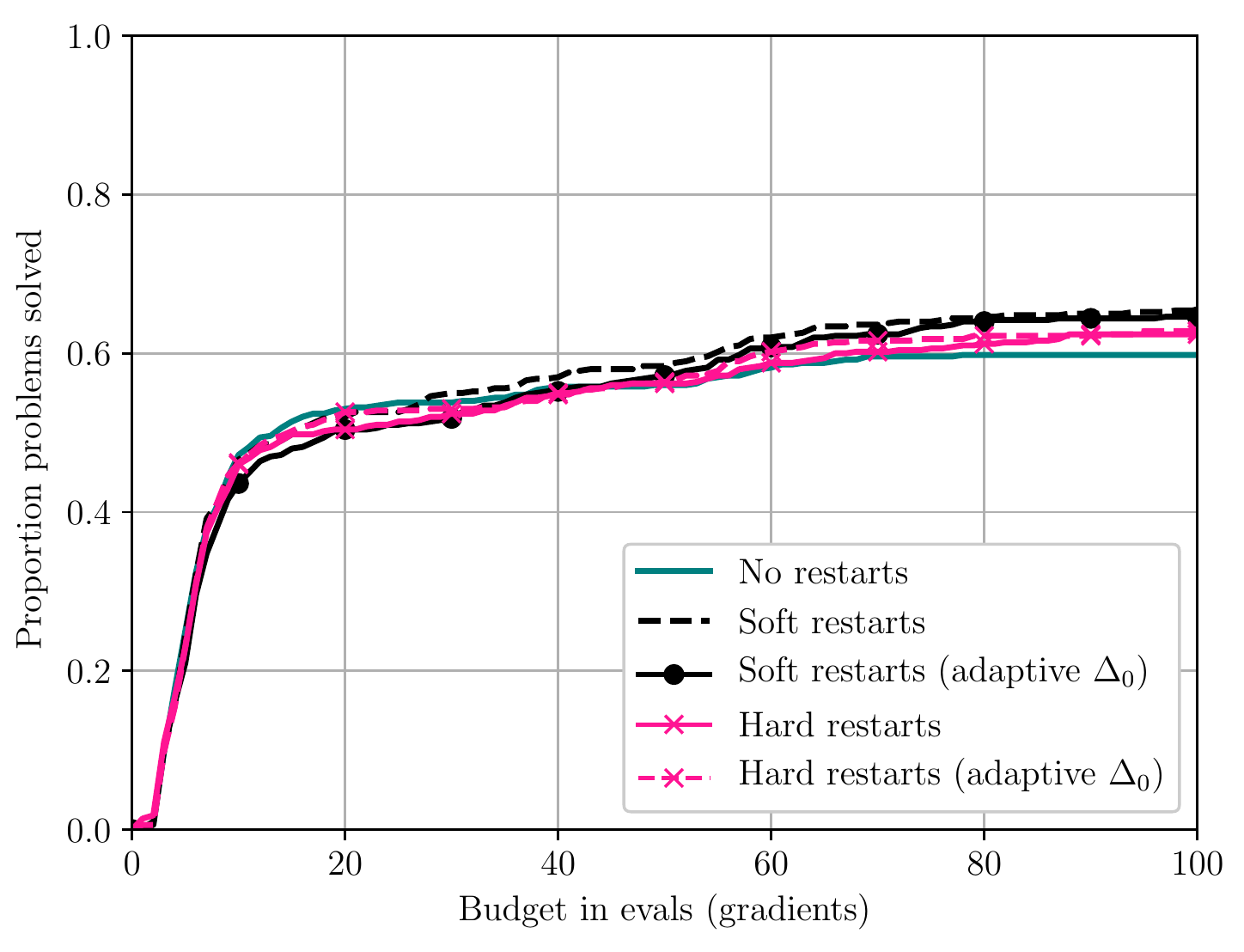}
		\caption{Comparison of Py-BOBYQA ($p=\mathcal{O}(n^2)$) variants  on GO test set with multiplicative noise}
	\end{subfigure}
	\\
	\begin{subfigure}[b]{0.48\textwidth}
		\includegraphics[width=\textwidth]{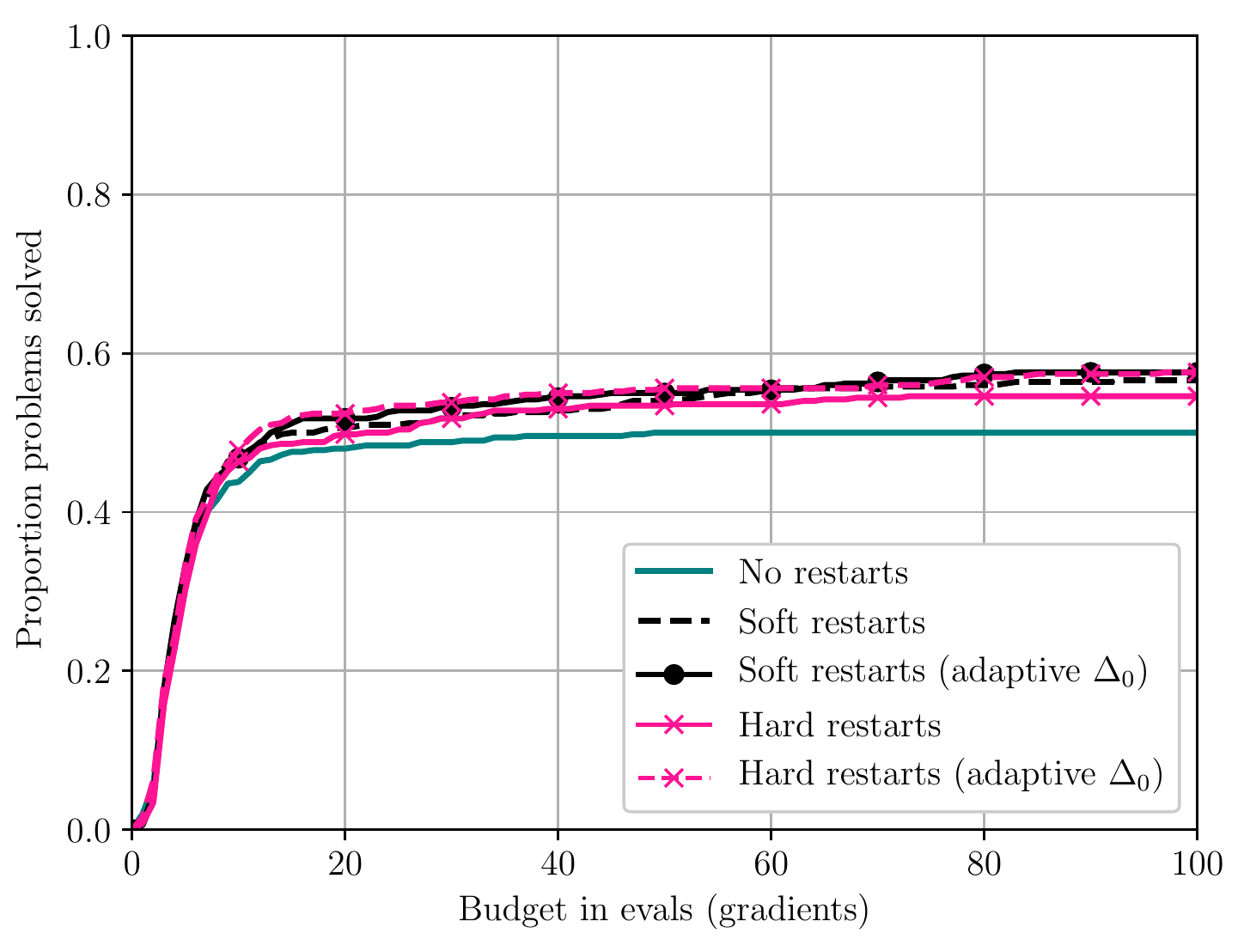}
		\caption{Comparison of Py-BOBYQA ($p=2n+1$) variants  on GO test set with additive noise}
	\end{subfigure}
	~
	\begin{subfigure}[b]{0.48\textwidth}
		\includegraphics[width=\textwidth]{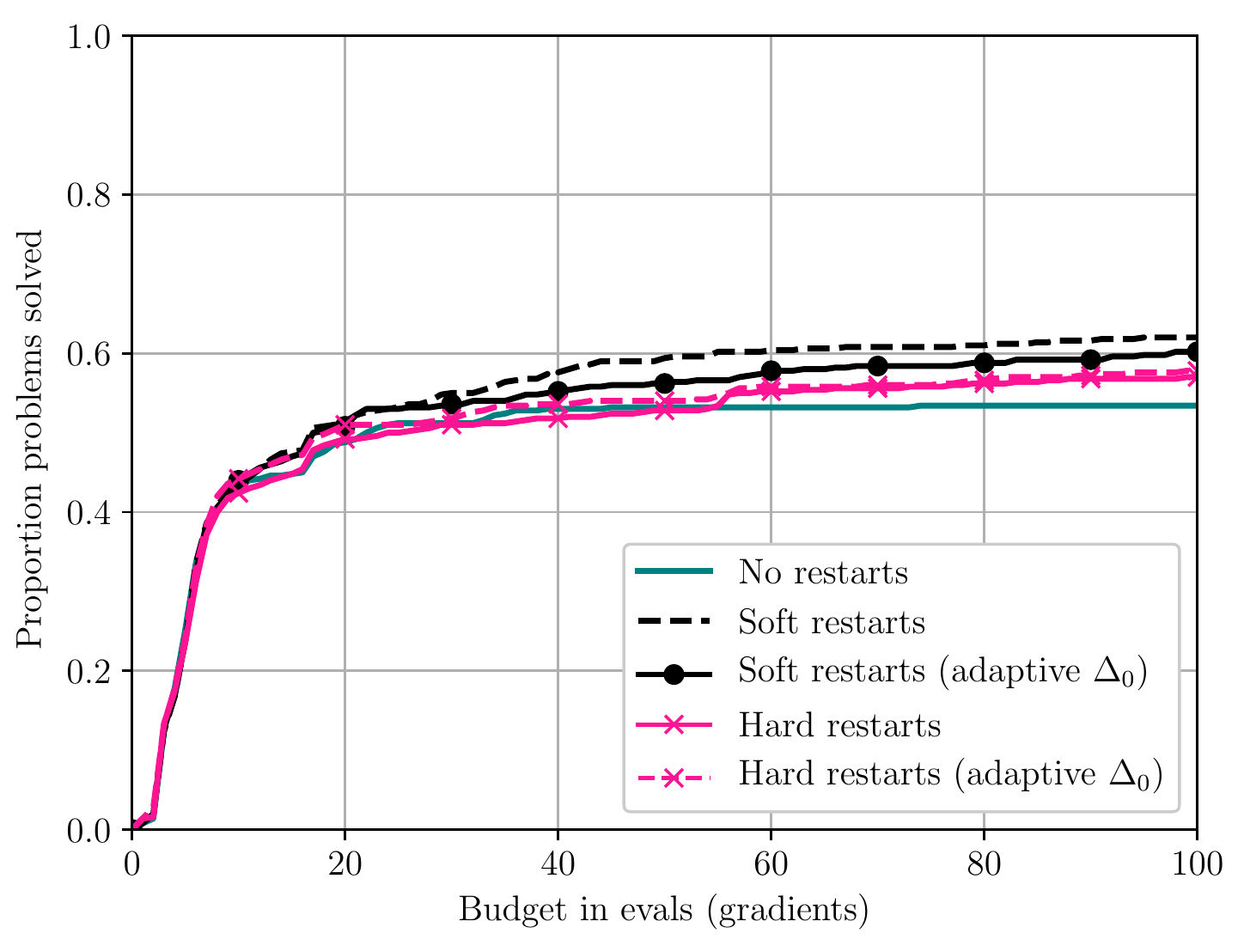}
		\caption{Comparison of Py-BOBYQA ($p=\mathcal{O}(n^2)$) variants  on GO test set with additive noise}
	\end{subfigure}
	\caption{Comparison of Py-BOBYQA variants on GO test set---low accuracy ($\tau=10^2$) and budget ($10^2 (n+1)$ evaluations).}
	\label{pybobyqa-low}
\end{figure}

\subsection{Py-BOBYQA Variance Testing: $2n+1$ interpolation points}
\begin{figure}[H]
	\centering
	\begin{subfigure}[b]{0.48\textwidth}
		\includegraphics[width=\textwidth]{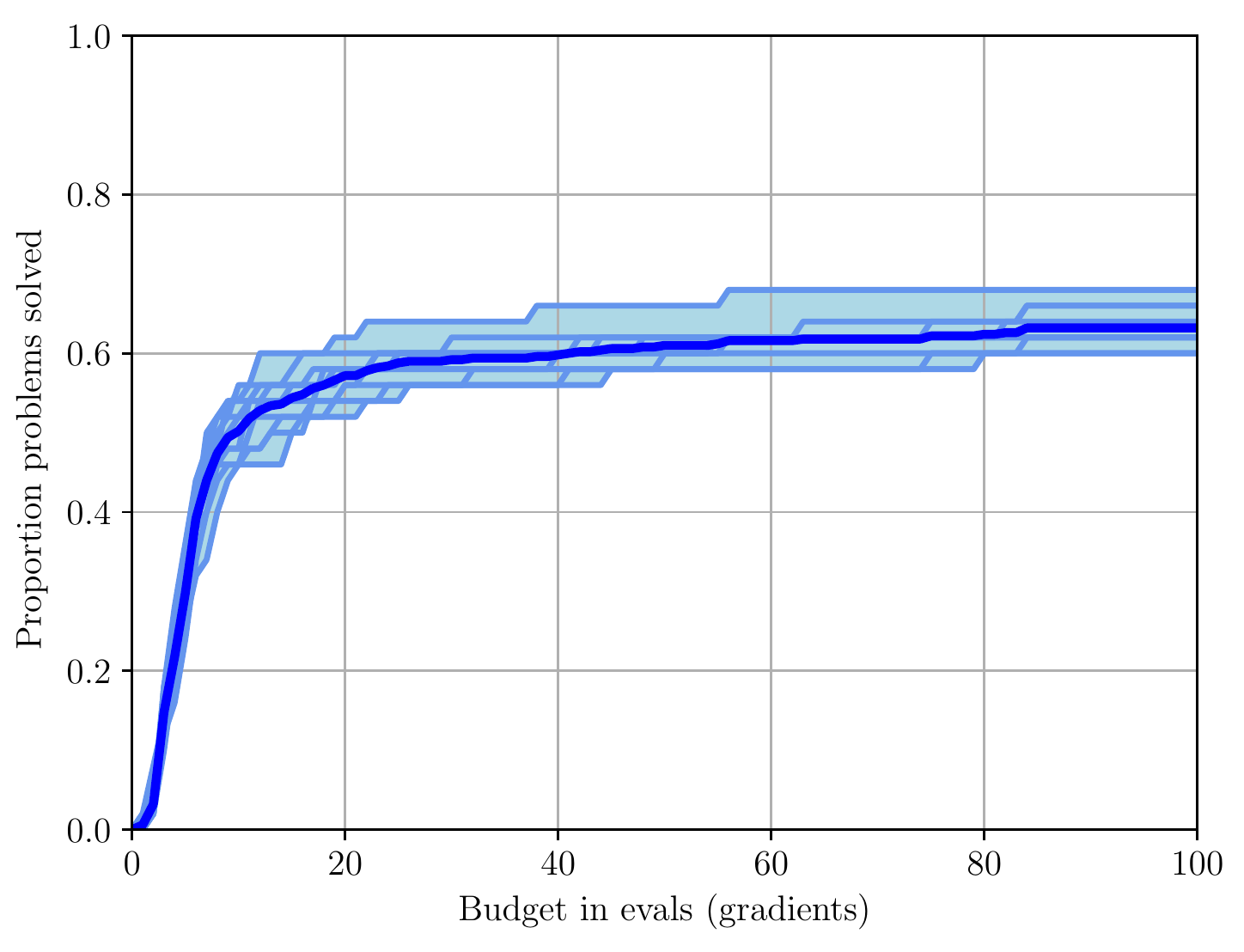}
		\caption{Fixed restart radius, low accuracy/budget.}
	\end{subfigure}
	~
	\begin{subfigure}[b]{0.48\textwidth}
		\includegraphics[width=\textwidth]{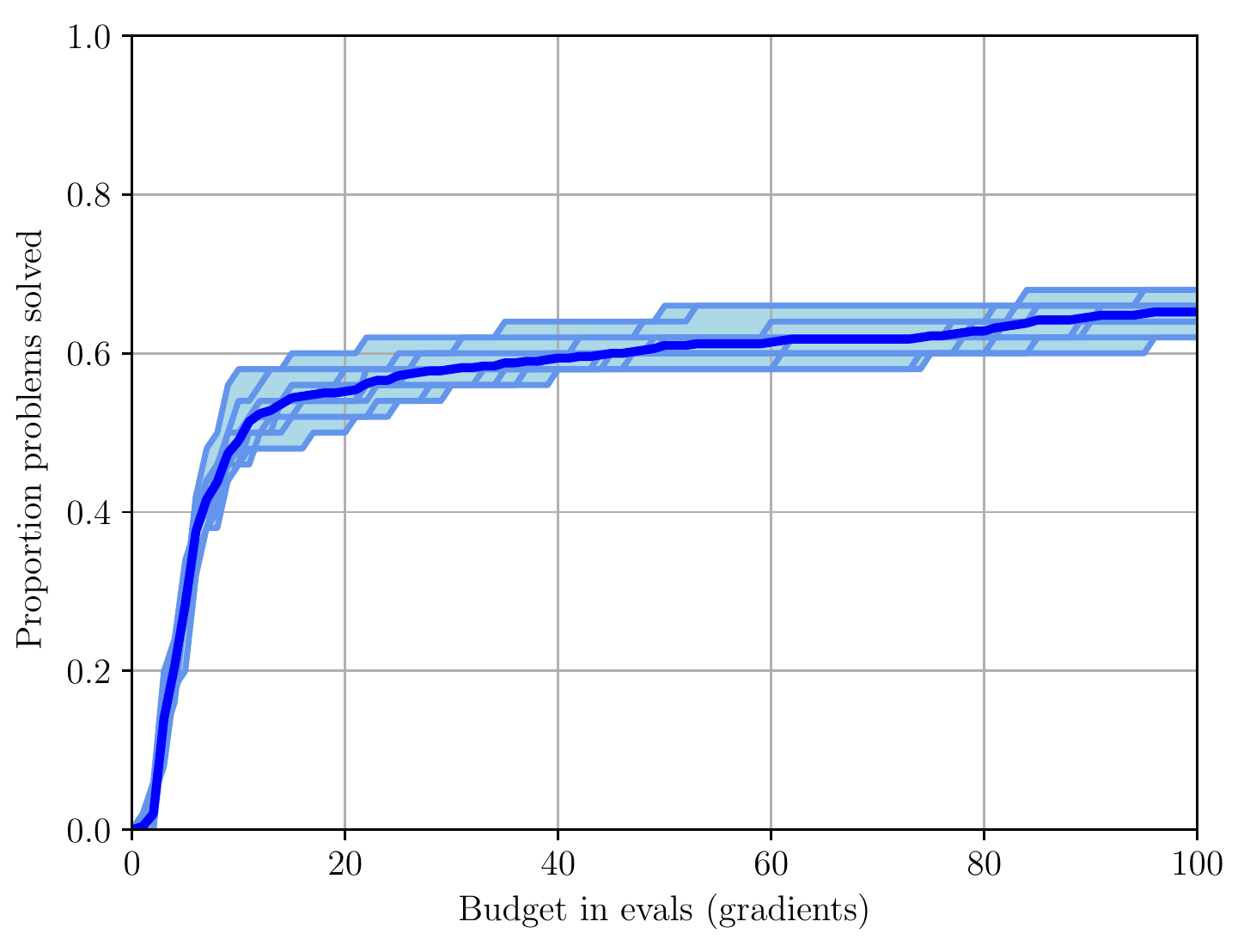}
		\caption{Adaptive restart radius, low accuracy/budget.}
	\end{subfigure}
	\\
	\begin{subfigure}[b]{0.48\textwidth}
		\includegraphics[width=\textwidth]{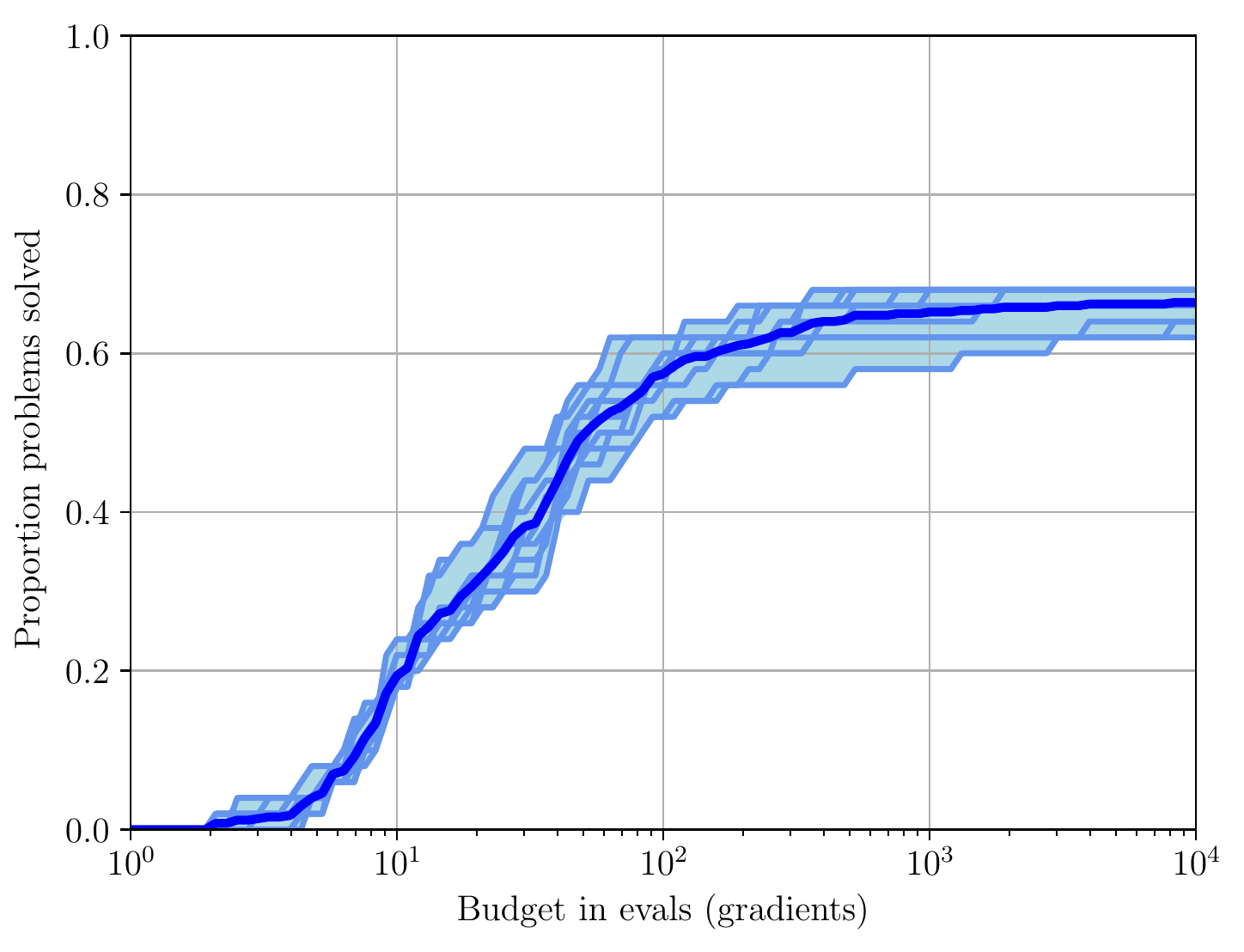}
		\caption{Fixed restart radius, high accuracy/budget.}
	\end{subfigure}
	~
	\begin{subfigure}[b]{0.48\textwidth}
		\includegraphics[width=\textwidth]{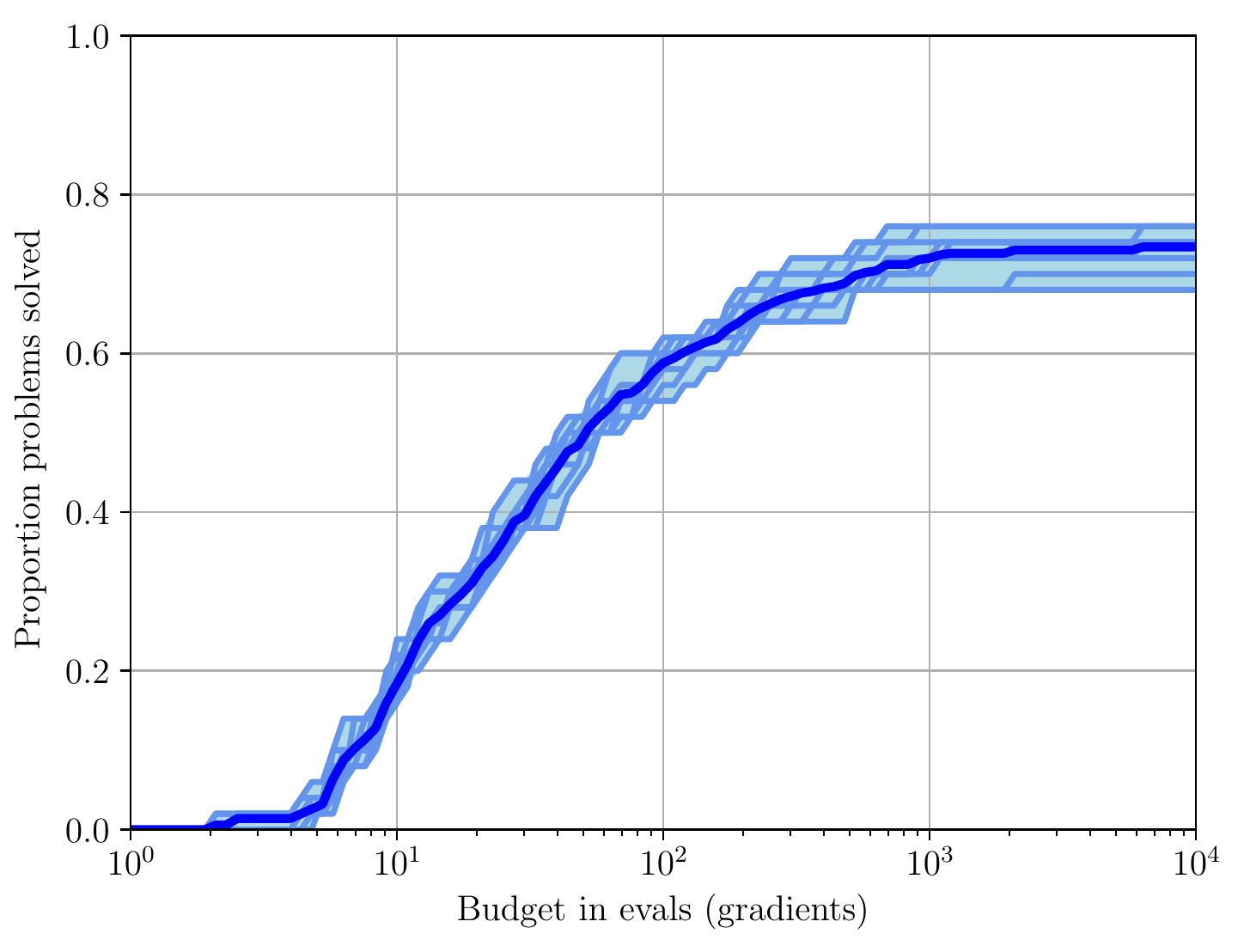}
		\caption{Adaptive restart radius, high accuracy/budget.}
	\end{subfigure}
	\caption{Demonstration of impact of random starting point for Py-BOBYQA ($p=2n+1$ with soft restarts); compare with \figref{fig_pybobyqa_variance_nsq} for $p=\bigO(n^2)$ points. Low accuracy and budget are $\tau=10^{-2}$ and $10^2 (n+1)$ evaluations respectively; high accuracy and budget are $\tau=10^{-5}$ and $10^4 (n+1)$ evaluations respectively. The dark lines are the average of all runs.}
	\label{fig_pybobyqa_variance_2np1}
\end{figure}

\subsection{Selecting the best Bayesian/surrogate solver: noisy objectives}
\begin{figure}[H]
	\centering
	\begin{subfigure}[b]{0.48\textwidth}
		\includegraphics[width=\textwidth]{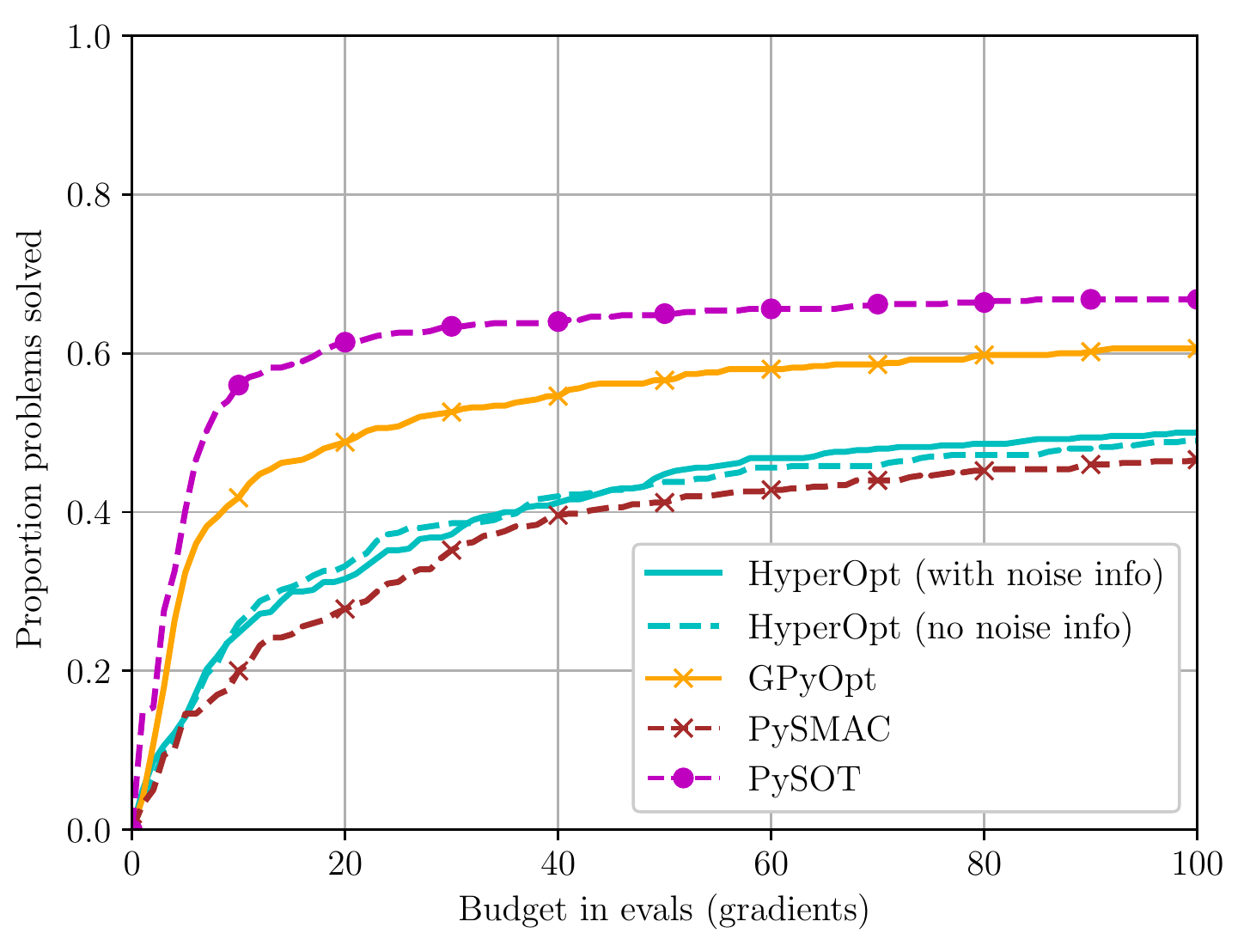}
		\caption{Low accuracy ($\tau=10^{-2}$), multiplicative noise.}
	\end{subfigure}
	~
	\begin{subfigure}[b]{0.48\textwidth}
		\includegraphics[width=\textwidth]{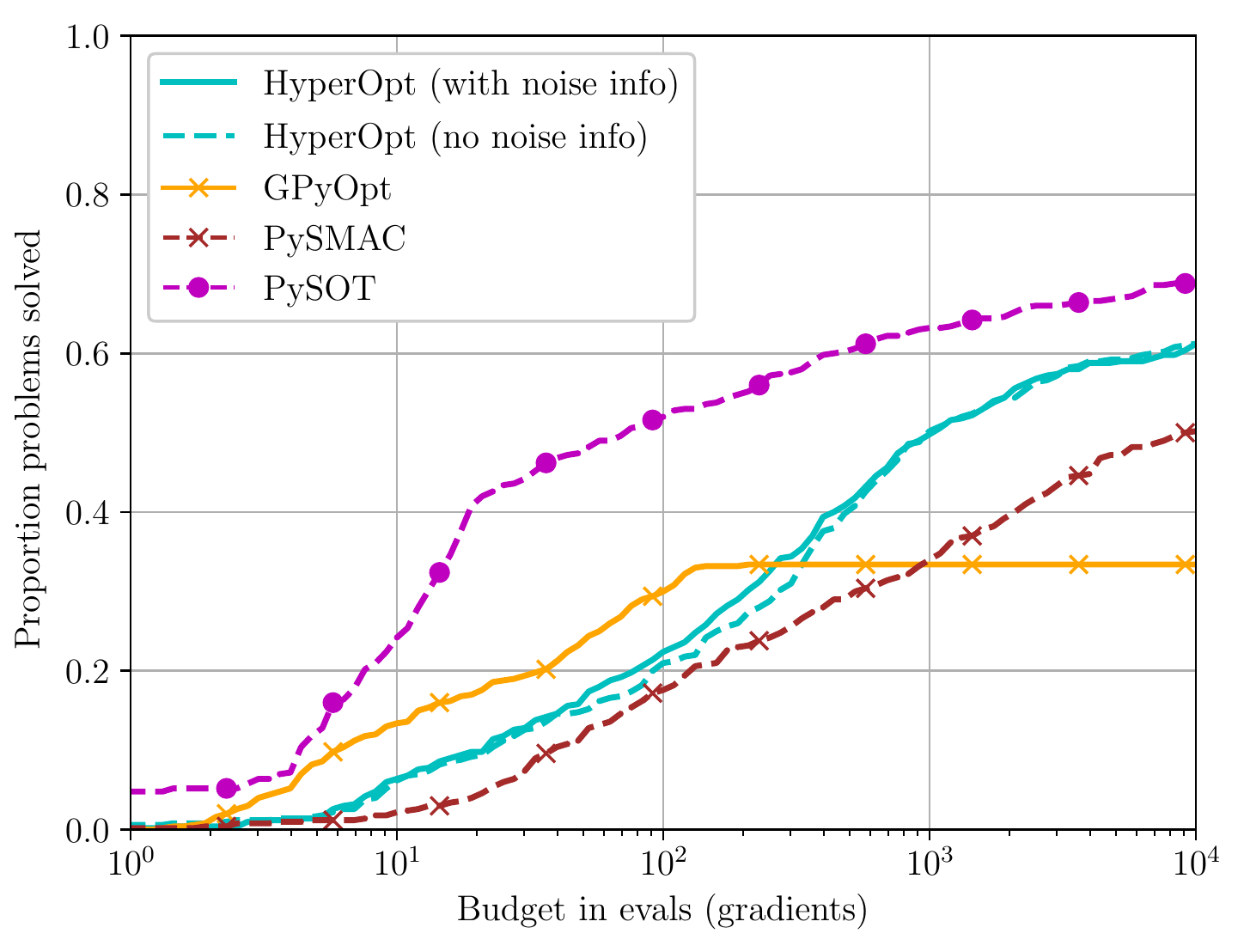}
		\caption{High accuracy ($\tau=10^{-5}$), multiplicative noise.}
	\end{subfigure}
	\\
	\begin{subfigure}[b]{0.48\textwidth}
		\includegraphics[width=\textwidth]{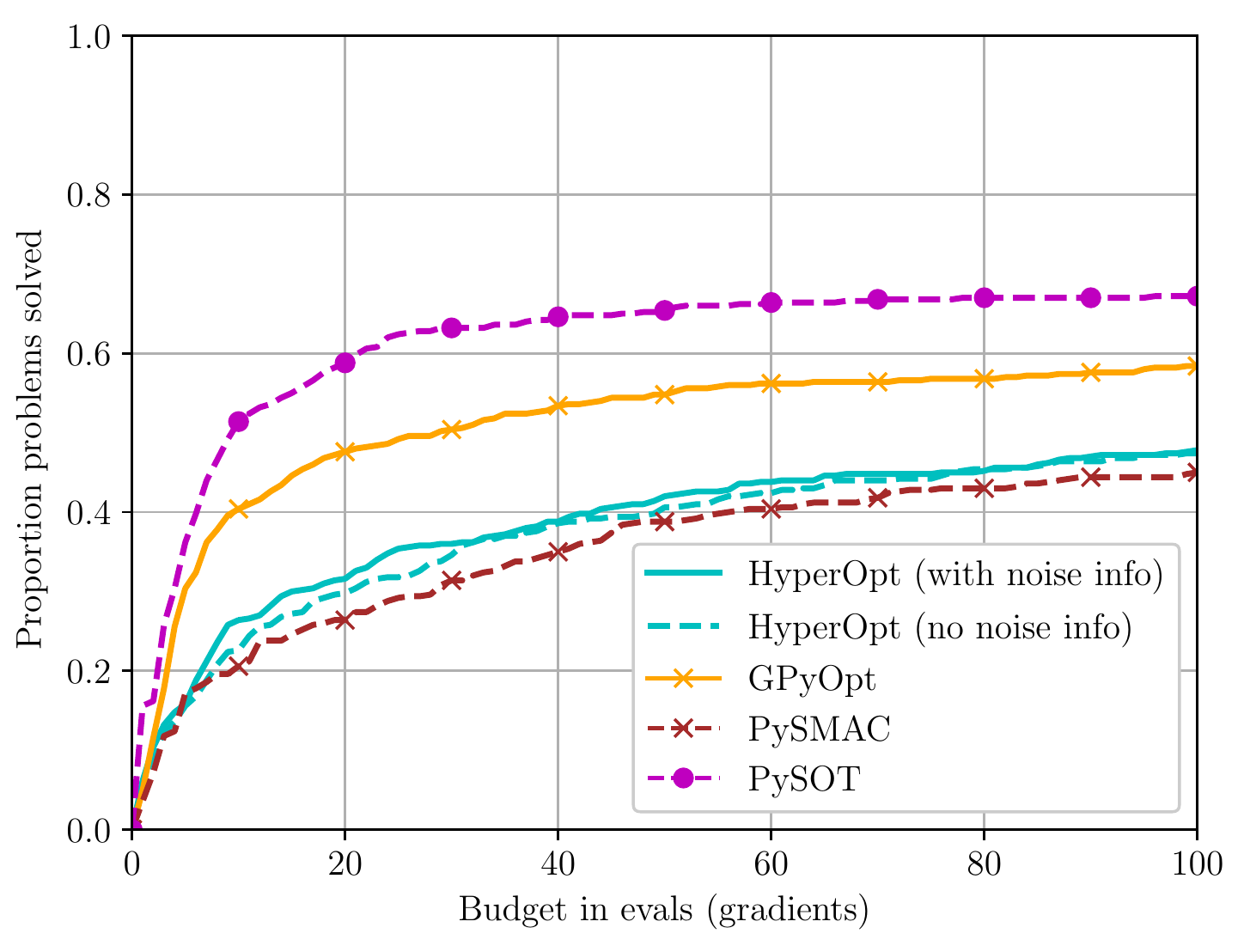}
		\caption{Low accuracy ($\tau=10^{-2}$), additive noise.}
	\end{subfigure}
	~
	\begin{subfigure}[b]{0.48\textwidth}
		\includegraphics[width=\textwidth]{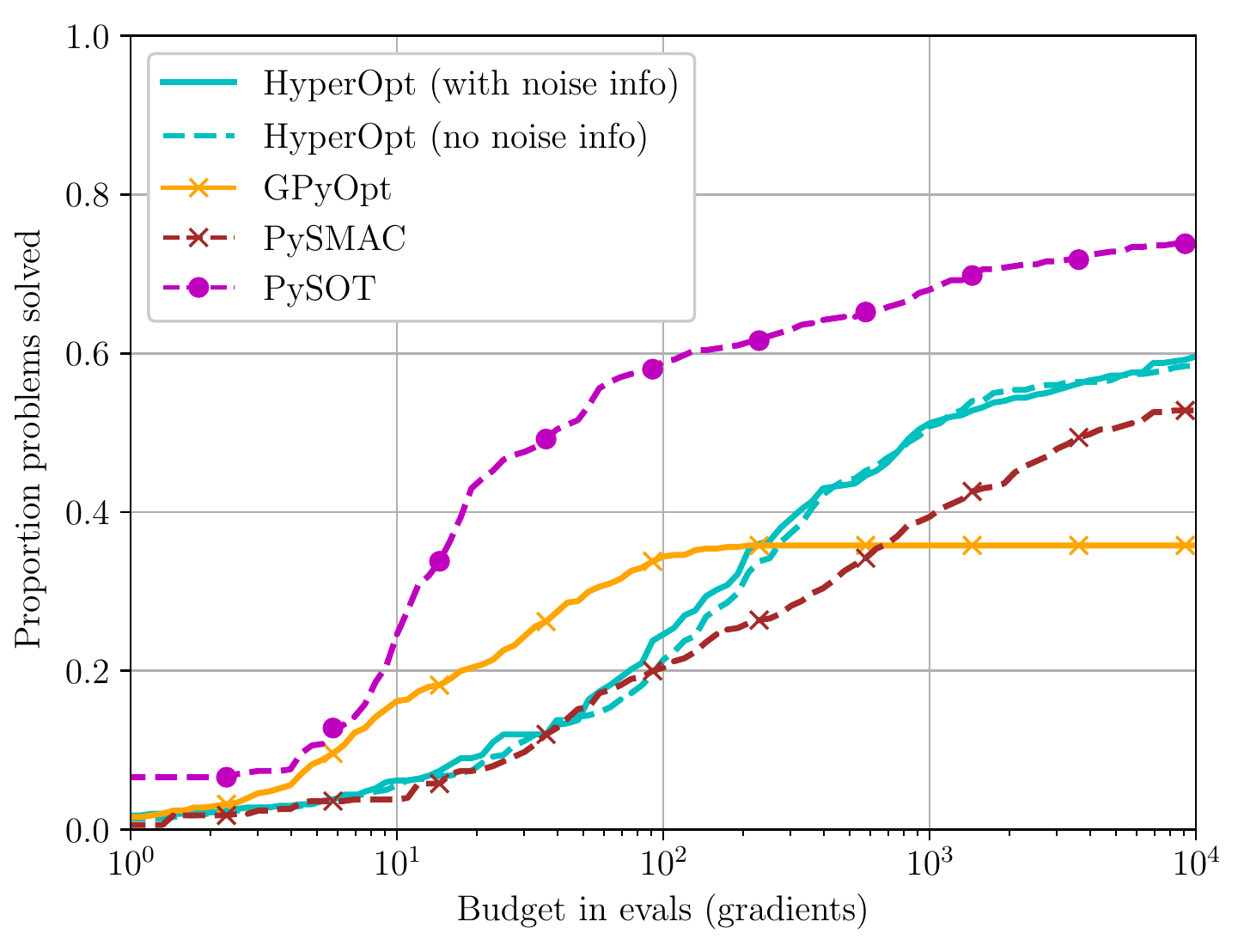}
		\caption{High accuracy ($\tau=10^{-5}$), additive noise.}
	\end{subfigure}
	\caption{Comparison of Bayesian and surrogate solvers on the GO test set for noisy problems; see \figref{BayesianComparison-GO} for the smooth case.}
	\label{BayesianComparison-GO-noisy}
\end{figure}


\end{document}